\documentclass[11pt, reqno]{amsart} 
\usepackage[utf8]{inputenc}
\usepackage [english]{babel}
\usepackage{a4wide}
\usepackage{amssymb}
\usepackage{amsmath}
\usepackage{amsthm}
\usepackage{amsbsy}
\usepackage{mathrsfs}
\usepackage{amsfonts}
\usepackage{mathtools}
\usepackage{verbatim}
\usepackage{tikz}

\usepackage{bbm}
\usepackage{tikz-cd}
\usetikzlibrary{cd}
\usepackage{amscd}
\usepackage{hyperref}
\usepackage{xcolor}
\hypersetup{
	colorlinks,
	linkcolor={blue!80!black},
	citecolor={red!90!black},
	urlcolor={black}
}
\usepackage[normalem]{ulem}
\usepackage{cleveref}   
\usepackage[all]{xy}
\usepackage{xcolor}
\theoremstyle{plain}
\usepackage[margin=1 in, a4paper]{geometry}
\usepackage{relsize}
\usepackage{lipsum}
\usepackage[autostyle]{csquotes}
\usepackage{relsize}
\usepackage{enumitem}  
\usepackage{titlesec}
\usepackage{setspace}
\usepackage{adjustbox}
\usepackage{crossreftools}
\usepackage{titlesec}
\allowdisplaybreaks

\titleformat{\section}
{\centering\large\sc}{\thesection.}{1em}{}

\makeatletter

\def\l@section{\@tocline{1}{0pt}{0pc}{1em}{}}

\def\l@subsection{\@tocline{2}{0pt}{1em}{3em}{}}

\makeatother

\newcommand{\Sym}{\mathrm{Sym}}

\newcommand{\F}{{\mathbb{F}}}
\newcommand{\Z}{{\mathbb{Z}}}
\newcommand{\Q}{{\mathbb{Q}}} 

\newcommand{\N}{{\mathbb{N}}}    
    
\newcommand{\C}{{\mathbb{C}}}

\newcommand{\p}{{\mathfrak{p}}}

\newcommand{\FF}{{\mathcal{F}}}

\newcommand{\OO}{{\mathcal{O}}}

\newcommand{\lra}{\longrightarrow}

\newcommand{\Gal}{\mathrm{Gal}}
\newcommand{\Ker}{\mathrm{Ker}}
\newcommand{\Coker}{\mathrm{Coker}}
\newcommand{\tor}{\mathrm{tor}}

\newcommand{\Hom}{\mathrm{Hom}}
\newcommand{\ord}{\mathrm{ord}}

\newcommand{\rank}{\mathrm{rank}}
\newcommand{\rk}{\mathrm{rk}}
\newcommand{\Frob}{\mathrm{Frob}}

\newcommand{\n}{\noindent}

\newcommand{\str}{\mathrm{str}}

\newcommand{\BK}{\mathrm{BK}}
 \newcommand{\genlegendre}[4]{%
  \genfrac{(}{)}{}{#1}{#3}{#4}%
  \if\relax\detokenize{#2}\relax\else_{\!#2}\fi
}

\newcommand{\ur}{\mathrm{ur}}

\newcommand{\pst}{\mathrm{pst}}

\makeatletter

\renewcommand\subsection{\@startsection{subsection}{2}{\z@}%
  {-2.25ex\@plus -1ex \@minus -.2ex}%
  {-1em}%
  {\normalfont\bfseries\raggedright}}
\makeatother

\usepackage[OT2,T1]{fontenc}
\DeclareSymbolFont{cyrletters}{OT2}{wncyr}{m}{n}
\DeclareMathSymbol{\Sha}{\mathalpha}{cyrletters}{"58}

\newcommand{\G}{\mathcal{G}}
\newcommand{\dR}{\mathrm{dR}}
\newcommand{\crys}{\mathrm{crys}}
\newcommand{\st}{\mathrm{st}}
\newcommand{\WD}{\mathrm{WD}}

\newtheorem{theorem}{Theorem}[section]
\newtheorem*{theorem*}{Theorem}
\newtheorem{corollary}[theorem]{Corollary}
\newtheorem{lemma}[theorem]{Lemma}
\newtheorem*{lemma*}{Lemma}
\newtheorem{proposition}[theorem]{Proposition}
\newtheorem{conjecture}[theorem]{Conjecture}

\theoremstyle{definition}

\newtheorem{remark}[theorem]{Remark}
\newtheorem{definition}[theorem]{Definition}

\newtheorem{example}[theorem]{Example}

\newtheorem*{heegnerhyp*}{Heegner Hypothesis[Heg]}
\numberwithin{figure}{section}			
\numberwithin{table}{section}				

\newtheorem{maintheorem}{Theorem}

\begin{document}
	\title[$p$-\MakeLowercase{parity conjecture for}  P\MakeLowercase{anchishkin} G\MakeLowercase{alois  representations}]{$p$-\MakeLowercase{parity conjecture for}  P\MakeLowercase{anchishkin} G\MakeLowercase{alois  representations}}
\author{S\MakeLowercase{ubham} G\MakeLowercase{arai},\ S\MakeLowercase{omnath} J\MakeLowercase{ha},\ \MakeLowercase{and} S\MakeLowercase{udhanshu} S\MakeLowercase{hekhar}} 
	\address{Department of Mathematics and Statistics, IIT Kanpur,  Kanpur-208016, India}\email{subhamg23@iitk.ac.in, jhasom@iitk.ac.in, sudhansh@iitk.ac.in}
\begin{abstract}

Let $p$ be an odd prime, and let $V_1$ and $V_2$ be symplectic self-dual $p$-adic de Rham representations of $G_{\mathbb{Q}}$ that are pure (of weight $-1$) and satisfy the Panchishkin condition at $p$. Suppose that $V_1$ and $V_2$ are residually isomorphic, compatibly with their Panchishkin filtrations.  We prove the relative $p$-parity conjecture relating the ratio of the global root numbers of $V_1$ and $V_2$ to the difference of their Bloch-Kato Selmer dimensions. We have  no restriction on the dimension of the Galois representation, the Hodge-Tate weights or the monodromy operator at $p$, and the residual isomorphism need not be symplectic. The proof expresses ratios of local root numbers in terms of local invariants over the cyclotomic $\mathbb{Z}_p$-extension: a new invariant introduced at $p$, and local Iwasawa $\lambda$-invariants at $\ell\neq p$. Both are compatible with the arithmetic local constants of Mazur-Rubin and Nekov\'a\v{r}. A key input, proved using Sen theory, is that the Bloch-Kato and strict local conditions at $p$ coincide over the cyclotomic $\mathbb{Z}_p$-extension, generalising results of Coates-Greenberg and Nekov\'a\v{r}. As applications, we obtain new cases of the $p$-parity conjecture, including cases of rank greater than two: Artin twists of CM Hida families, orthogonal Artin twists of ordinary modular forms, symmetric cubes along Hida families, and families of abelian surfaces, together with unconditional numerical examples.
\end{abstract}
\maketitle
\tableofcontents

\section{Introduction}
\subsection{$p$-parity and relative $p$-parity conjectures} We fix an odd prime $p$.  For a fix algebraic closure $\bar{\Q}$ of $\Q$ and define the embedding $\iota_\infty : \bar{\Q} \hookrightarrow \C$ and for any finite prime $\ell$ of $\Q$, $\iota_\ell: \bar{\Q} \hookrightarrow \bar{\Q}_\ell$.
 Let  $M$ be a pure motive over a number field $L$ with coefficients in  a number field $F$. For a prime $\p\mid p$ of $F$, let $V=M_\p$ be a finite dimensional $F_\p$ vector space equipped with a continuous action of $G_L=\Gal(\bar{\Q}/L)$. Let $S_\BK(V/L)$ be the Bloch-Kato Selmer group associated to $V$ over $L$ (see \eqref{global_BK_selmer_group}). Let $L(M,s)$ be the complex $L$-function associated to $M$. It is expected to have a meromorphic continuation to $\C$ and the completed $L$-function $\tilde{L}(M,s)$ is conjectured to  satisfy a functional equation (see \cite[($C_{FE}$)]{Ne3} for a precise statement). Let $\chi_p: G_\Q \rightarrow \Z_p^\times$ denote the $p$-adic cyclotomic character and for any $r\in \N$, put $V^*(r):= \Hom_{F_\p} (V, F_\p(r))$,  with $F_\p(r):= F_\p \otimes \chi_p^r$. Then the Bloch-Kato conjecture states that (see \cite{bk, fp}),
\begin{equation}\label{BK_conjecture}
    \ord_{s=0}L(M,s)= \rank_{F_\p}\ S_\BK(V^*(1)/L) - \rank_{F_\p}\ H^0(L, V^*(1)).
\end{equation}
Assume that $V$ is a symplectic self-dual representation (see Definition \ref{ssd}). Then $V=V^*(1)$ and the sign of the conjectural functional equation of $\tilde L(M,s)$   determines the parity of $\ord_{s=0}\ L(M,s)$. In fact, more generally,  for a symplectic self-dual $p$-adic representation $V$ of $G_L$, which is unramified at all but finitely many primes and de Rham at $p$, we can define the global root number $W(V)$ (see definition \eqref{global_root}). This global root number $W(V)$  conjecturally coincides with the sign of the expected  functional equation of the completed complex $L$-function associated to $V$. 
 Now we can define the mod $2$ analogue of the Bloch-Kato conjecture \eqref{BK_conjecture} as the $p$-parity conjecture.
\begin{conjecture}\label{p-parity} ($p$-parity conjecture) Let $V$ be a  symplectic  self-dual $p$-adic representation  of $G_L$, which is unramified at all but finitely many primes and de Rham at $p$. Then we have 
\begin{equation*}
    W(V)= (-1)^{\rank_{F_\p}\ S_\BK(V/L)-\rank_{F_\p}\ H^0(L, V)}. \tag*{\qed}
\end{equation*}
\end{conjecture}
Note that if $V$ is symplectic self-dual and pure (of weight $-1$), then  $H^0(L, V)=0$.  Conjecture \ref{p-parity} implies the following relative version for two $p$-adic representations $V_1$ and $V_2$ of $G_L$.

\begin{conjecture}\label{relative_p-parity}(Relative $p$-parity conjecture) 
  Let $V_1$ and $V_2$ be   two symplectic  self-dual $p$-adic representations  of $G_L$, that are unramified at all but finitely many primes and de Rham at $p$ and pure (of weight $-1$). Then we have 
\begin{equation*}
    \frac{W(V_1)}{W(V_2)}= (-1)^{\rank_{F_\p}\ S_\BK(V_1/L)- \rank_{F_\p}\ S_\BK(V_2/L)}. \tag*{\qed}
\end{equation*}
\end{conjecture}



\subsection{Previous works} We now recall some known results concerning Conjectures~\ref{p-parity} and~\ref{relative_p-parity}. The earliest results in this direction go back to Birch and Stephens, who established $p$-parity results for  elliptic curves over $\Q$ \cite{BS} with $j$-invariant equal to $0$ and $1728$.  Greenberg and Guo discussed $p$-parity results for CM elliptic curves over $\Q$ (see \cite{Gr6, Gu2}). Subsequently, Monsky proved the $2$-parity conjecture for all elliptic curves over $\Q$ \cite{Mo} and also discussed $p$-parity for elliptic curves with a degree $p$-isogeny. He used results of Kolyvagin, Cassels and Kramer. Significant progress on the $p$-parity conjecture was made by Nekov\'{a}\v{r}, who established the $p$-parity conjecture for elliptic curves over $\Q$ that have good ordinary reduction at $p$ or multiplicative (ordinary) reduction at $p$, or additive reduction at $p$ that becomes ordinary or multiplicative over an extension.  For a  prime $p$ of supersingular reduction for an elliptic curve over $\Q$,  Kobayashi introduced the theory of plus/minus Selmer groups, which was used  in Kim's proof of the $p$-parity conjecture for elliptic curves over $\Q$ at a supersingular prime $p>3$ \cite{Ki}. 

The above mentioned works do not cover the case of  an elliptic curve over $\Q$ with good supersingular reduction at $p=3$ or additive reduction at an odd $p$ that becomes good supersingular over an extension of $\Q$. A different  approach was developed by Tim and Vladimir Dokchitser, based on \emph{regulator constants} and relations between permutation representations of finite groups. Their powerful method  led to $p$-parity results for the remaining cases over $\Q$ and also applies for elliptic curves over certain number fields and their Artin twists \cite{Dokchister-Dokchister_parity-conjecture,DD3}.  Another contribution is due to \v{C}esnavi\v{c}ius  \cite{Ce}, who completed the proof of the $p$-parity conjecture for elliptic curves admitting a rational $p$-isogeny over an arbitrary number field, by removing certain restrictions on reduction types at primes above $p$, which were required for Dokchitser-Dokchitser. It is also shown in  \cite{Ce} that if an elliptic curve $E$ has CM over a number field $K$, then  the $p$-parity conjecture holds for $E/K$. 

Combining Nekov\'{a}\v{r}'s later results and techniques, overall his work on $p$-parity conjecture also covers all non-CM and most CM elliptic curves over any totally real number field, most  Hilbert modular forms of parallel weight two with trivial central character,  and many cases of abelian varieties with real multiplication. In his work on the $p$-parity conjecture for elliptic curves or Hilbert modular forms, Nekov\'{a}\v{r} developed a method where,  (under certain assumptions), he first established a relative $p$-parity result (see Conjecture \ref{relative_p-parity})  for two (or a family of) congruent mod-$p$ Galois representations  and then specializes to a situation where the $p$-parity result is known for one of the member. A central ingredient in this approach is the theory of \emph{arithmetic local constants}, certain local invariants introduced by Mazur and Rubin \cite{MR}, which describe the variation of Selmer ranks modulo $2$.  Nekov\'{a}\v{r} investigated the compatibility of \emph{arithmetic local constants} with the corresponding local root numbers, proving this compatibility at places $\ell\neq p$. The case $\ell=p$ is substantially more delicate and requires additional hypotheses and techniques from $p$-adic Hodge theory. By establishing the required compatibility in progressively more general situations, he obtained the relative $p$-parity results. Further, combining relative parity results with known cases of the $p$-parity conjecture, Nekov\'{a}\v{r} obtained far-reaching results for elliptic curves over totally real fields, Hilbert modular forms of parallel weight two, and abelian varieties with real multiplication. Moreover, for cuspidal Hilbert modular forms of weight $(2k_1,\dots 2k_d, 2)$ with trivial central character, some cases of $p$-parity conjecture were established by Johansson-Newton \cite{JN}.

In a recent preprint,  Burungale-Kobayashi-Nakamura-Ota  \cite{bkno} proved the relative $p$-parity conjecture for two residually isomorphic symplectic self-dual $p$-adic representations of $G_\Q$ of rank two.
\subsection{Aim of the article} In this article, we aim to establish new cases of the relative $p$-parity and $p$-parity conjectures for a broader class of $p$-adic Galois representations, with particular emphasis on cases not covered by existing results, including those in rank greater than two. We also discuss to what extent Nekovář’s methods extend to this more general setting and clarify the limitations of his techniques.

Nekov\'{a}\v{r}'s  relative $p$-parity results in \cite{Ne6} for  pairs of symplectically self-dual $p$-adic Galois representations with isomorphic residual representations, assume that their restrictions to the decomposition groups at primes above $p$ become Barsotti–Tate over finite abelian extensions of the corresponding local fields.  Note that potentially Barsotti–Tate representations are potentially crystalline, with Hodge-Tate weights in $\{ 0,1\}$ and Weil–Deligne monodromy operator vanishes $N=0$. Our results will cover the cases of symplectic self-dual, Panchishkin and de Rham Galois representations of arbitrary Hodge-Tate weights and we do not assume the vanishing of monodromy operator.   In a different direction, his deformation-theoretic result \cite{Ne3} on  relative $p$-parity for a family of pure, symplectically self-dual Galois representations satisfying the Panchishkin condition, assumes that the family is defined over a complete Noetherian local integral domain of Krull dimension two, that is, the family is defined over an irreducible base of dimension two. In this work we consider an arbitrary pair of $p$-adic Galois representations which a priori need not be members of such a family. In particular, we also include families of Galois representations parameterised by arbitrary complete Noetherian local rings which are not necessarily domains that is families over non irreducible bases are also included.

In a different direction, the $p$-parity conjecture has also been studied using Iwasawa-theoretic methods by Greenberg in  \cite{Greenberg_projectivemodule_modularrepresentation} and Coates-Fukaya-Kato-Sujatha \cite{CFKS}. Greenberg investigated the relationship between the parity of Selmer ranks and global root numbers for Artin twists of elliptic curves with good ordinary reduction at $p$. Using the structure of Selmer groups over cyclotomic extensions and modular representation theory, he established, under the vanishing of Iwasawa $\mu$ invariant, that the parity conjecture is preserved under congruences between self-dual Artin representations. On the other hand, \cite{CFKS} studied Selmer groups over certain noncommutative $p$-adic Lie extensions and established parity results for twists of elliptic curves by orthogonal Artin representations of the $p$-adic Lie extension under certain additional Iwasawa-theoretic hypotheses. They also proved the $p$-parity conjecture for certain abelian varieties admitting suitable $p$-power isogenies, including elliptic curves with a rational $p$-isogeny under specified local conditions.

Building on and extending the work of Nekov\'{a}\v{r} and the results in \cite{Greenberg_projectivemodule_modularrepresentation,CFKS}, in this article, we will study the relative $p$-parity conjecture for a more general class of $p$-adic Galois representations that are unramified outside finitely many primes, satisfy the de Rham and  the Panchishkin condition at $p$. We will then apply these results to establish new cases of the $p$-parity conjecture.
\subsection{Main results of the relative $p$-parity conjecture} Now we fix some notation and state our main results. Some of the terminology involved are introduced and explained in  detail in later sections.

We say a finite dimensional $p$-adic representation $V$ of $G_\Q$ satisfies \hyperref[Hyp:Mov]{(Mov)} if $V$ is unramified outside finitely many primes, de Rham at $p$ and its coefficient field is a finite extension  of $\Q_p$.
Let $V_1$ and $V_2$ be two  $p$-adic Galois representations of $G_{\Q}$ satisfying \hyperref[Hyp:Mov]{(Mov)} with coefficient field $\FF$. Assume that $V_1$ and $V_2$ are symplectic self-dual representations (Definition \ref{ssd}) and satisfy the Panchishkin condition (Definition \ref{panchishkin})  at $p$. 
We denote the ring of integers of $\mathcal{F}$  by $\mathcal{O}$ and let  $\pi$  be a uniformizer of $\OO$ and   $\mathfrak{f}$ be the  residue field of $\OO$. For $m\in \{1,2\}$, let $T_m$ be a symplectic self-dual $G_{\Q}$-stable  $\OO$-lattice (see Remark \ref{self-dual_lattice}) of $V_m$ and put $A_m:= V_m/T_m$. Let $A_m^+$ be the image of $V_m^+$ under the natural map $V_m\lra A_m$ and set $A_m^-:= A_m/A_m^+$. For any $\mathcal{O}$-module $Y$, let $Y[\pi]$ denote the $\pi$-torsion submodule of $Y$. For a discrete or a compact $\Z_p$-module $Y$, let $Y^\vee$ denote the Pontryagin dual of $Y$.

We assume that $A_1[\pi] \cong A_2[\pi]$  as $G_\Q$-modules. 
 Let $\Q_\infty$ be the cyclotomic $\mathbb{Z}_p$-extension of $\Q$ with $\Gamma:= \Gal(\Q_\infty /\Q)\cong \mathbb{Z}_p$ and let $\Lambda:= \OO[[\Gamma]]\cong  \OO[[X]]$ denote the associated Iwasawa algebra. Let $\Q_{\infty,p}$ denotes the completion of $\Q_\infty$ at the unique prime above $p$. Define 
\begin{equation*}
     H^\str (\Q_{\infty,p}, A_m[\pi]):= \Ker\left( H^1(\Q_{\infty,{p}}, A_m[\pi]) \lra H^1(\Q_{\infty,{p}}, A_m^-) \right).
\end{equation*}
Now we introduce the local  invariants at $p$,
\begin{align*}
   \delta(A_m,p):=  \dim_\mathfrak{f}\ \frac{ H^\str (\Q_{\infty,p}, A_m[\pi])}{ H^\str (\Q_{\infty,p}, A_1[\pi]) \cap  H^\str (\Q_{\infty,p}, A_2[\pi])},
\end{align*}
where the intersection is taken via the identification $H^1(\Q_{\infty,{p}}, A_1[\pi])=H^1(\Q_{\infty,{p}}, A_2[\pi])$.

On the other hand,  for a finite prime $\ell\neq p$, we take the local invariant as  the local Iwasawa $\lambda$-invariant defined as follows:
\begin{align*}
   \delta(A_m,\ell):=  \rk_\mathcal{O}\
       \Big( \prod_{\eta\mid \ell} H^1(\Q_{\infty,\eta} , A_m)^\vee \Big).
\end{align*}

Towards the relative $p$-parity conjecture, we have the following comparison of the global root numbers in terms of the local invariants $\delta(-,\ell)$:

\begin{theorem}\label{global_root_lambda}
     Let $V_1$ and $V_2$ be two symplectic self-dual $p$-adic Galois representations of $G_{\Q}$ satisfying \hyperref[Hyp:Mov]{(Mov)} such that $A_1[\pi] \cong A_2[\pi]$  as $G_\Q$-modules. Assume that $V_1$ and $V_2$ satisfy the Panchishkin condition  at $p$ and are pure (of weight $-1$) (Definition \ref{purity}).  We further assume that the isomorphism $A_1[\pi] \cong A_2[\pi]$ restricted to $G_{\Q_p}$  induces an isomorphism $A_1^-[\pi] \cong A_2^-[\pi]$. Then the ratio of the global root numbers of $V_1$ and $V_2$ is given by  
 \begin{equation*}
   \frac{W(V_1)}{W(V_2)}=(-1)^{\sum_{{\ell } }\  (\delta(A_1,\ell)- \delta(A_2,\ell))},  
\end{equation*}
where the sum is over   finite primes $\ell$ dividing $p$ and the primes where  $V_1$ or $V_2$ is ramified. \qed
\end{theorem}

Let $S_\BK(A_m/\Q_\infty)$ be the Bloch-Kato Selmer group of $A_m$ over $\Q_\infty$. Then 
recall that $S_\BK(A_m/\Q_\infty)^\vee$  is a finitely generated $\Lambda$-module and we denote the Iwasawa $\mu$-invariant of $S_\BK(A_m/\Q_\infty)^\vee$ by $\mu(S_\BK(A_m/\Q_\infty)^\vee)$.
Our main result on the relative $p$-parity conjecture (Conjecture \ref{relative_p-parity}) is the following:

\begin{theorem}\label{relative_p_parity_over_Q}
         Let $V_1$, $V_2$ be two $p$-adic representations of $G_\Q$ satisfying \hyperref[Hyp:Mov]{(Mov)}. Assume that they are  symplectic self-dual, pure (of weight $-1$),  and satisfy the Panchishkin condition at $p$ and $A_1[\pi]\cong A_2[\pi]$ as $G_\Q$-modules. We make the following additional assumptions.
     \begin{enumerate}
         \item $S_\BK(A_1/\Q_\infty)^\vee$ is $\Lambda$-torsion and $\mu(S_\BK(A_1/\Q_\infty)^\vee)$ is equal to $0$.
    \item    The given isomorphism $A_1[\pi] \cong A_2[\pi]$ when restricted to $G_{\Q_p}$  induces an isomorphism $A_1^-[\pi] \cong A_2^-[\pi]$. 
     \end{enumerate}
     Then \begin{equation*}
         \frac{W(V_1)}{W(V_2)}= (-1)^{\mathrm{rk}_\mathcal{O}\ S_{\BK}(A_1/\Q)^\vee- \mathrm{rk}_\mathcal{O}\ S_{\BK}(A_2/\Q)^\vee}.\tag*{\qed}
     \end{equation*} 
\end{theorem}

Note that $\rank_\OO\ S_{\BK} (A/\Q)^\vee = \rank_\FF\ S_{\BK} (V/\Q)$.
In \S \ref{arithmetic_examples}, we give several examples of  families of $p$-adic Galois representations 
   that satisfy the relative $p$-parity conjecture stated in Theorem \ref{relative_p_parity_over_Q}. 
   These examples include   congruent families (including Hida families) of elliptic cuspidal newforms twisted by Artin representations (\S \ref{example_modular-forms}), odd symmetric powers of congruent elliptic curves (\S \ref{symmetric_power_elliptic}), tensor product of symmetric powers of congruent elliptic curves (\S \ref{tensor_product_symmetric_powers}). We also prove a relative $p$-parity result for a given $p$-adic representation twisted by residually isomorphic Artin representations (\S \ref{twist_by_artin_representation}). In particular, this establishes the relative $p$-parity conjecture for a $p$-adic Galois representation associated to a given elliptic cuspidal newforms twisted by two congruent Artin representations,  which generalizes a result of Greenberg \cite{Greenberg_projectivemodule_modularrepresentation}. 
Our results on   relative $p$-parity and $p$-parity also cover  Galois representations exhibiting exceptional zero
 phenomenon (see \cite{Gr3}, \cite{Greenberg_Stevens})(for example, modular abelian varieties with split multiplicative reduction at primes above $p$).
\subsection{Proof ideas, novelty and difficulty}
Now we discuss some key ideas, novelty and difficulty in the main results described above.

To  compare  the global root numbers for congruent $p$-adic Galois representations (Theorem \ref{global_root_lambda}), we compare the local root numbers and patch them to obtain the global root number comparison. The most subtle case here is the local root number comparison at the prime $p$. Under the hypothesis that the residual representations of the  two $p$-adic representations are residually symplectically isomorphic, Nekov\'{a}\v{r} \cite{Ne, Ne6} uses the {\it arithmetic local constant} defined by Mazur-Rubin 
and under some additional assumptions (see \cite{Ne, Ne6}) relates  it with the ratio of local root numbers at $p$. 
On the other hand,  in \cite{Ne3}, he has compared the ratio of local root numbers at $p$ of two Galois representations which come from a family defined over a complete Noetherian local {\it integral domain} of Krull dimension two. Further in their study of relative $p$-parity result, the authors of  \cite{bkno}  mainly focus on $\ell=p$ case and consider  rank two $p$-adic Galois representations.

To circumvent these restrictions for a general $p$-adic de Rham, Panchishkin Galois representation of arbitrary rank, we climb up the cyclotomic $\Z_p$ tower and introduce a new invariant, $\delta(-,p)$ \eqref{delta_at_p} and relate it with the ratio of local root numbers at $p$. Moreover, we also show the compatibility of $\delta(-,p)$ with the {\it arithmetic local constant} in Proposition \ref{leqp}. For a prime $\ell\neq p$, we also relate the difference of local Iwasawa invariants $\delta(-,\ell)$ with the ratio of local root numbers in our general setting.  We also prove the compatibility results for {\it arithmetic local constant} and our local invariant $\delta(-,\ell)$  in Lemma \ref{lnp}.
We mention that for the local root number comparison at a   finite prime (including $p$) or at the infinite prime,  we do not assume that the isomorphism of the residual Galois representations preserves the symplectic pairing. See Example \ref{example_isomorphic_residual_but_not_symplectic_similitude}, Lemma \ref{sym} and Remark \ref{remark 4.15} for a detailed discussion.

Our local invariants are defined over cyclotomic $\Z_p$-extension. To
 compare them with root numbers and Selmer ranks over number fields, we need a careful study of the structure of $(\varphi, N)$-module associated to $p$-adic Galois representations. We use  Sen's theory to study the $\Gamma$-module structure of $H^0(L_v(\mu_{p^\infty}), V^-)$, where $L_v$ is a finite extension of $\Q_p$ (see Lemma \ref{Sen_theory}, Proposition \ref{trivial_action_gamm_v} and Corollary \ref{mod_2_Qcyc_Q}). In particular, Corollary \ref{mod_2_Qcyc_Q} plays a crucial role in comparing root numbers with local invariants $\delta(-,p)$. 
 
 In an influential work on Kummer theory for abelian varieties, Coates and Greenberg \cite{CG} compares the local Kummer map at prime $p$ with strict Selmer local condition at $p$ over deeply ramified extensions of $\Q_p$. Using a modified version of  Cartesian local condition at $p$ defined in \cite{MR1}, we extend the result in \cite{CG} to a more general symplectic self-dual $p$-adic Galois representation. More precisely, let $H^{\str}(\Q_{\infty,p},A)$ (resp. $H^\BK(\Q_{\infty,p},A)$)   be the strict (resp. Bloch-Kato) Selmer local condition of $A=V/T$ over $\Q_\infty$ at $p$ defined in \S \ref{Selmer Groups}.   Our main result on this is the following.

\begin{theorem}
   Let $V$ be a self-dual symplectic $p$-adic Galois representation of $G_{\Q_p}$ which is de Rham. Also assume that $V$ satisfies the Panchishkin condition and is pure (of weight $-1$). 
    Then 
    \begin{equation*}
        H^{\str}(\Q_{\infty,p},A)=H^\BK(\Q_{\infty,p},A). \tag*{\qed} 
    \end{equation*}
 \end{theorem}

As a corollary of this, we deduce the equality  of Bloch-Kato and strict Selmer groups over $\Q_\infty$, $S_\BK(A/\Q_\infty)= S_\str(A/\Q_\infty)$ (see Theorem \ref{BK_str_equality_over_cyclotomic}). Due to this, we can work with strict Selmer group.   We  extend  techniques of    Greenberg-Vatsal (see \cite{Greenberg_Vatsal}) to express the difference of coranks of strict Selmer groups   over cyclotomic $\Z_p$-extension with sum of our local invariants $\delta(-,\ell)$. At this stage, the Iwasawa theoretic hypotheses that $S_\BK(A/\Q_\infty)^\vee$ is $\Lambda$-torsion and $\mu(S_\BK(A/\Q_\infty)^\vee)=0$ appears to be essential. Then extending techniques from Greenberg, Guo \cite{Greenberg1997,Guo_Li},  we compare the parity of coranks of strict Selmer groups over $\Q_\infty$ and $\Q$. This enables us to prove the relative $p$-parity result (see Theorem \ref{relative_p_parity_over_Q}).

Note that in \S \ref{arithmetic_examples},  we discuss the prevalence of $p$-adic Galois representations in literature satisfying the hypotheses of our relative $p$-parity result in Theorem \ref{relative_p_parity_over_Q}.  Moreover, in \S \ref{numerical_examples}, we have   numerical examples where we verify the hypotheses of Theorem \ref{relative_p_parity_over_Q}.

\subsection{Application to the $p$-parity conjecture} 
Having established  relative $p$-parity results, we follow the strategy considered by Nekov\'{a}\v{r} and others to specialize to a situation where the $p$-parity is known for one of the member. This enables us to prove, under suitable hypotheses, several new cases of the $p$-parity conjecture in higher rank ($>2$) cases. 

Castella \cite[Theorem C]{Castella_Heegner_points} has established, under certain hypotheses (see (Ca1)-(Ca4) in Theorem \ref{castella_main_thm}), Bloch-Kato conjecture for CM elliptic curves twisted by two dimensional Artin representation, when the complex  $L$-value of the CM elliptic curve twisted by the Artin representation is non-zero. Using this, we prove a $p$-parity result for a family of  modular forms varying over a CM Hida family twisted by a two dimensional Artin representation. More precisely, we have the following theorem:
\begin{maintheorem}\label{castella_main_thm}
    Let $p$ be an odd prime. Suppose that $E/\Q$ is a $p$-ordinary elliptic curve of conductor $M$ and $E$ has CM by an imaginary quadratic field $L$. Let $K$ be a Galois number field and $\sigma$ be a two-dimensional, odd, irreducible orthogonal Artin representation of $\Gal(K/\Q)$. Assume that 
\begin{enumerate}
    \item[(Ca1)] $p\nmid 6NM$ and $\gcd(N,M)=1$, where $N$ is the Artin conductor of $\sigma$.
    \item[(Ca2)] The prime $p$ splits as  $p=\mathfrak{p}\bar{\mathfrak{p}}$ in the ring of integers of $L$.
    \item[(Ca3)] Every prime $q|N$ either splits or ramifies in $L$, with $q^2\nmid N$ in the latter case.
    \item[(Ca4)] $\sigma(\mathrm{Frob}_p)$ has distinct eigenvalues modulo $\pi$.
    \item[(5)] The Iwasawa $\mu$ invariant  $\mu(S_\BK(E\otimes\sigma/\Q_\infty)^\vee)$ is equal to $0$.
\end{enumerate}
    If $L(E\otimes \sigma/\Q,1)\neq 0$, then for every member $f_n$ in the Hida family $\{f_n\}_{n}$ associated to $E$, with weight $k_n\equiv 2 \pmod{2(p-1)}$, the $p$-parity conjecture holds for the $p$-adic Galois representation $V_{f_n}(\sigma):= V_{f_n}(1-\frac{k_n}{2}) \otimes \sigma$ associated to   $f_n\otimes\sigma$, i.e., 
    \begin{equation*}
        W(V_{f_n}(\sigma))=(-1)^{\rk _\mathcal{O}\ S_{\BK}(A_{f_n}(1-\frac{k_n}{2})\otimes\sigma/\Q)^\vee}. \tag*{\qed}
    \end{equation*}
\end{maintheorem}

On the other hand, let $f\in S_2(\Gamma_0(M))$ be a non-CM $p$-ordinary  cuspidal newform.  Let $K$ be a totally real Galois number field and assume that the order of vanishing of the complex $L$-function of $f$ over $K$, $L(f/K,s)$ at $s=1$ is at most $1$. In this setting, Shou-Wu Zhang \cite{Shouwu_Zhang_Heights_Of_Heggner_points}  has extended the Gross-Zagier and Kolyvagin's result for  abelian varieties associated to $f$. Using the result in \cite{Shouwu_Zhang_Heights_Of_Heggner_points}, we  prove the following $p$-parity result for every member of  the Artin twist of the Hida family passing through $f$.

\begin{maintheorem}\label{zhang_p_parity}
   Let $p$ be an odd prime. Let $K$ be a totally real Galois number field and $f \in S_2(\Gamma_0(M))$ be a $p$-ordinary newform   with trivial nebentypus. For a prime $\lambda$  of the Hecke field $\Q_f :=\Q(\{a_n(f)\}_{n\geq 1})$ dividing $p$,  let $V_{f,\lambda}$ be the associated $\lambda$-adic Galois representation. Let $\sigma$ be an orthogonal, irreducible Artin representation of $\Gal(K/\Q)$. 
   Further assume that 
  \begin{enumerate}
      \item  $\mathrm{ord}_v(M)=1$ for at least one finite place $v$ of $K$. 
      \item The order of vanishing of the complex $L$-function of $f$ over $K$, $L(f/K,s)$ at $s=1$ is at most $1$.
      \item $S_\BK(A_{f,\lambda}(\sigma)/\Q_\infty)^\vee$ is $\Lambda$-torsion and $\mu(S_\BK(A_{f,\lambda}(\sigma)/\Q_\infty)^\vee)$ is equal to $0$ for every $\lambda\mid p$.
      \end{enumerate} 
     Then for every member $f_n$ of the  Hida family associated to $f$ with weight $ k_n\equiv 2 \pmod{2(p-1)}$ and for every $\lambda|p$, the $p$-parity conjecture holds for   $f_n\otimes \sigma$,  i.e., 
     \begin{equation*}
         W(V_{f_n,\lambda}(\sigma)) =(-1)^{\rk_\mathcal{O}\ S_{\BK}(A_{f_n,\lambda}(1-\frac{k_n}{2})\otimes \sigma/\Q)^\vee }. \tag*{\qed}
     \end{equation*} 
    \end{maintheorem}

Some important cases of $p$-parity conjecture for abelian varieties twisted by Artin representation are discussed  in \cite[Theorem 4.5]{Dokchister-Dokchister_parity-conjecture}.  
Let $A$ be an abelian variety over $\Q$ and $K/\Q$ is a Galois extension of number fields such that the Sylow $p$-subgroup $P$ of $G=\Gal(K/\Q)$ is normal and $G/P$ is abelian. In \cite{Dokchister-Dokchister_parity-conjecture}, they  first establish a key technical result (called {\it Hypothesis 4.1} in \cite{Dokchister-Dokchister_parity-conjecture}) and  using this, they show  if the $p$-parity conjecture is true for $A/\Q$ and $A/L$ for any quadratic extension $\Q \subset L \subset K$, then the $p$-parity conjecture holds  for $A$  twisted by any orthogonal Artin representations of $G$. 

Using our relative $p$-parity result (Theorem \ref{relative_p_parity_over_Q}), we prove an analogue of  \cite[Hypothesis 4.1]{Dokchister-Dokchister_parity-conjecture} for more general $p$-adic representations in Lemma \ref{hypothesis_4.1} and combining it with  results of Nekov\'{a}\v{r} on the $p$-parity conjecture for modular forms,  we prove the following general $p$-parity result for Artin twists of modular forms.

\begin{maintheorem}\label{P_parity_as_DD}
     Let $p$ be an odd prime.  Suppose that $K/\Q$ is a Galois extension of number fields such that the Sylow $p$-subgroup $P$ of $G=\Gal(K/\Q)$ is normal and $G/P$ is abelian.  Let $f \in S_{2k}(\Gamma_0(M))$ be a $p$-ordinary cuspidal newform with even weight $2k$ and trivial nebentypus  such that $S_\BK(A_f(1-k)/K_\infty)^\vee$ is $\Lambda$-torsion and $\mu(S_\BK(A_f(1-k)/K_\infty)^\vee)$ is equal to $0$. Then for any   irreducible orthogonal Artin representation $\sigma$ of $G$ the $p$-parity conjecture holds for $f\otimes \sigma$ over $\Q$, i.e.,
     \begin{equation*}
         W(V_f(\sigma))= (-1)^{\rk_\OO\ S_\BK(A_f(\sigma)/\Q)^\vee}. \tag*{\qed}
     \end{equation*}
\end{maintheorem}

Let $f\in S_2(\Gamma_0(N))$ be a $p$-ordinary elliptic cuspidal newform such that the $L$ value $L(f,1)\neq 0$ and the central critical complex $L$ value,  $L(\Sym^3 f, 2)\neq 0$. In this setting,
 H. Wang \cite[Theorem 5.8]{wa} has shown, under suitable hypotheses (see (Wa1)-(Wa5) in Theorem \ref{wang_result_p_parity}), the Bloch-Kato conjecture for $\Sym^3 f$.  Using this,  we prove a $p$-parity result for the symmetric cubes of (higher weight) cuspidal newforms in a Hida family passing through $f$. 

\begin{maintheorem}\label{wang_result_p_parity}
    Let $p\geq 5$ be a prime. Suppose that $f\in S_2(\Gamma_0(N))$ is a $p$-ordinary normalised newform with trivial nebentypus. We assume that $N$ admits a factorization $N=N^+N^-$ such that $\gcd(N^+, N^-)=1$ and $N^-$ is square-free and is a product of an odd number of primes.  We make the following additional assumptions.
    \begin{enumerate}
        \item[(Wa1)] The residual representation $\bar{\rho}_f$ is absolutely irreducible.
        \item[(Wa2)] The $\mathbb{T}_{\mathfrak{m}_f}$ module $S_2^B(N^+, \OO)_{\mathfrak{m}_f}$ is free of rank $1$ (for definition see \cite{wa}). 
        \item[(Wa3)] The residual representation $\bar{\rho}_f$ is surjective. 
        \item[(Wa4)] The $L$ value $L(f,1)$ is non-vanishing.
         \item[(Wa5)] The central critical $L$ value $L(\Sym^3 f, 2)$ is non-zero.
        \item[(6)] For $V= (\Sym^{3}V_f)(-1)$, $S_\BK(A/\Q_\infty)^\vee$ is $\Lambda$-torsion and $\mu(S_\BK(A/\Q_\infty)^\vee)=0$.
       
    \end{enumerate}
Then for every member $f_n$ of   the Hida family $\{f_n\}_n$ associated to $f$, with weight $k_n \equiv 2 \pmod{2(p-1)}$, the $p$-parity conjecture holds for $V_n:= (\Sym^{3}V_{f_n}(1-\frac{k_n}{2}))(-1)$, i.e., 
\begin{equation*}
    W(V_n) = (-1)^{\rk_\OO\ S_\BK(A_n/\Q)^\vee}. \tag*{\qed}
\end{equation*}
\end{maintheorem}

Further, we have also discussed  numerical examples in \S \ref{numerical_examples}, where all the hypotheses are satisfied and  the $p$-parity conjecture holds unconditionally.





\textbf{Structure of the article:} In \S \ref{preliminaries}, we   discuss some preliminary setup, definitions and notations. Then in \S \ref{Kummer_theory},  we have shown the equality of the local  Bloch-Kato  and strict  Selmer conditions at  $p$ over the cyclotomic $\Z_p$-extension. Next in \S \ref{cmparison_of_root_number}, we compare the local root numbers of two $p$-adic Galois representations at both finite and Archimedean places  and patching them, we obtain the relation between the ratio of the global root numbers with the local invariants $\delta(-,\ell)$.
Then in \S \ref{Selmer Groups}, we define the Selmer groups and establish the equality of the (global) Bloch-Kato and strict Selmer groups over the cyclotomic $\Z_p$-extension. In  \S \ref{relative_p_parity_conjecture}, we establish our main relative $p$-parity result and give several families of $p$-adic Galois representations, where relative $p$-parity conjecture holds. Finally, in \S \ref{some_new_cases_p_parity}, using relative $p$-parity result, we prove some new cases of $p$-parity conjecture. We also exhibit  numerical examples of families of $p$-adic representations such that the $p$-parity conjecture holds unconditionally for them.

\textbf{Acknowledgements:} We thank   C. S. Rajan for his suggestions and encouragement. The first named author thanks Abhishek for several helpful discussions and acknowledges the support of  IIT Kanpur institute fellowship for PhD students.

\section{Preliminaries}\label{preliminaries} 
We fix an odd prime $p$.  For a fix algebraic closure $\bar{\Q}$ of $\Q$ and define the embedding $\iota_\infty : \bar{\Q} \hookrightarrow \C$ and for any finite prime $\ell$ of $\Q$, $\iota_\ell: \bar{\Q} \hookrightarrow \bar{\Q}_\ell$. For any field $L$, $G_L$ will denote the absolute Galois group of $L$. Note that $G_{\Q_\ell}$ can be identified as the decomposition subgroup of $G_\Q$ at $\ell$. Let $I_\ell$ be the inertia subgroup of $G_{\Q_\ell}$ and denote by $\Frob_\ell$, the arithmetic Frobenius at $\ell$. For each finite prime $\ell$, $\F_\ell$ be the finite field of $\ell$ elements.  Let $\chi_p: G_\Q \rightarrow \Z_p^\times$ denote the $p$-adic cyclotomic character and $\omega_p: G_\Q \longrightarrow \mathbb{F}_p^\times \hookrightarrow \Z_p^\times$ be the Teichm\"{u}ller character. Let $\C_p$ denote the $p$-adic completion of $\bar{\Q}_p$ and we put $\C_p(i):= \C_p \otimes_{\Q_p} \Q_p(\chi_p^i)$.

Throughout this article, $\FF$ will denote the $p$-adic completion of a number field $F$ at a prime $\pi|p$.
  The ring of integers of $\mathcal{F}$ will be denoted by $\mathcal{O}$ and we will use the same notation $\pi$ to denote a uniformizer of $\OO$. Further,  $\mathfrak{f}$ will denote the  residue field of $\OO$.

 For a discrete or a compact $\Z_p$-module $X$, let $X^\vee$ denotes the Pontryagin dual of $X$.
  If $Y$ is a finitely generated $\mathcal{O}$-module, define $\rk_\mathcal{O}\, Y$ to be the $\mathcal{O}$-rank of $Y$ and the $\mathcal{O}$-corank of  $Y$ to be the $\mathcal{O}$-rank of $Y^\vee$.
  For an $\OO$-module $X$, $X_\mathrm{div}$ denotes the maximal $\pi$-divisible submodule of $X$. 
  For a finitely generated $\OO$-module $X$, $X_\tor$ will denote the torsion submodule of $X$.
 For any $\mathcal{O}$-module $Y$, let $Y[\pi]$ denotes the $\pi$-torsion submodule of $Y$.
\subsection{Theory of $p$-adic  representation} 
Recall the rings of $p$-adic periods
$$B_{\mathrm{crys}}\subset B_{\mathrm{st}}\subset B_{\mathrm{dR}},$$
constructed by \cite{Fontaine_94}, where the inclusion $B_{\mathrm{st}}\subset B_{\mathrm{dR}}$ depends on the choice of a branch of the $p$-adic logarithm, which we fix. Each of the above rings is a $\mathbb{Q}_p$-algebra equipped with a continuous action of $G_{\mathbb{Q}_p}$, compatible with the inclusions, along with additional structures: a Frobenius $\varphi$ on
$B_{\mathrm{cris}}$ and $B_{\mathrm{st}}$, a monodromy operator $N$ on $B_{\mathrm{st}}$ with
$N\varphi=p\varphi N$, and a decreasing filtration $\{\operatorname{Fil}^iB_{\mathrm{dR}}\}_{i\in\mathbb{Z}}$
on $B_{\mathrm{dR}}$.  For any finite dimensional $p$-adic representation $V$ of $G_{\mathbb{Q}_p}$, we associate the finite dimensional $\mathbb{Q}_p$-vector spaces
$$D_\dagger(V):=(V\otimes_{\mathbb{Q}_p}B_\dagger)^{G_{\mathbb{Q}_p}}\quad\text{for }\dagger\in\{\mathrm{dR},\mathrm{st},\mathrm{crys}\},$$
and the $\mathbb{Q}_p^{\mathrm{ur}}$-vector space$$D_{\mathrm{pst}}(V):=\varinjlim_{K}\,(V\otimes_{\mathbb{Q}_p}B_{\mathrm{st}})^{G_K},$$
where $K$ varies over finite extensions of $\mathbb{Q}_p$ contained in $\overline{\mathbb{Q}}_p$, and $\mathbb{Q}_p^{\mathrm{ur}}$ denotes the maximal unramified extension of $\mathbb{Q}_p$. One always has $\dim_{\mathbb{Q}_p}D_\dagger(V)\le\dim_{\mathbb{Q}_p}V$ for $\dagger\in\{\mathrm{dR},\mathrm{st},\mathrm{crys}\}$, and $\dim_{\mathbb{Q}_p^{\mathrm{ur}}}D_{\mathrm{pst}}(V)\le\dim_{\mathbb{Q}_p}V$. Further, for $\dagger\in\{\mathrm{dR},\mathrm{st},\mathrm{crys}\}$, $\dim_{\mathbb{Q}_p}D_\dagger(V)=\dim_{\mathbb{Q}_p}V$ if and only if the natural comparison map $B_\dagger\otimes_{\mathbb{Q}_p}D_\dagger(V)\to B_\dagger\otimes_{\mathbb{Q}_p}V$ is an isomorphism, and then the $p$-adic representation $V$ is called de Rham, semistable and crystalline depending on $\dagger=\mathrm{dR}$, $\dagger=\mathrm{st}$ and $\dagger=\mathrm{crys}$, respectively. The representation $V$ is called potentially semistable if $\dim_{\mathbb{Q}_p^{\mathrm{ur}}}D_{\mathrm{pst}}(V)=\dim_{\mathbb{Q}_p}V$, equivalently, if $V|_{G_K}$ is semistable for some finite extension $K$ of $\mathbb{Q}_p$. We have natural inclusions $D_{\mathrm{crys}}(V)\subseteq D_{\mathrm{st}}(V)\subseteq D_{\mathrm{dR}}(V)$, the equality $D_{\mathrm{st}}(V)=D_{\mathrm{pst}}(V)^{G_{\mathbb{Q}_p}}$, and a natural injection
\[
  \overline{\mathbb{Q}}_p\otimes_{\mathbb{Q}_p^{\mathrm{ur}}}D_{\mathrm{pst}}(V)\hookrightarrow\overline{\mathbb{Q}}_p\otimes_{\mathbb{Q}_p}D_{\mathrm{dR}}(V).
\]
Consequently, $V$ crystalline $\Longrightarrow$ $V$ semistable $\Longrightarrow$ $V$ potentially semistable $\Longrightarrow$ $V$ de Rham. In fact, if $V$ is a de Rham representation of $G_{\mathbb{Q}_p}$, then it is potentially semistable \cite[Theorem~0.7]{Ber}.

Moreover, $D_{\mathrm{dR}}(V)$ is equipped with the filtration of $\mathbb{Q}_p$-vector spaces
\[
  \operatorname{Fil}^iD_{\mathrm{dR}}(V):=(V\otimes_{\mathbb{Q}_p}\operatorname{Fil}^iB_{\mathrm{dR}})^{G_{\mathbb{Q}_p}},\qquad i\in\mathbb{Z},
\]
such that
\begin{itemize}
  \item $\operatorname{Fil}^{i+1}D_{\mathrm{dR}}(V)\subseteq\operatorname{Fil}^iD_{\mathrm{dR}}(V)$ for all $i\in\mathbb{Z}$.
  \item $\operatorname{Fil}^iD_{\mathrm{dR}}(V)=0$ for $i\gg0$, and $\operatorname{Fil}^iD_{\mathrm{dR}}(V)=D_{\mathrm{dR}}(V)$ for $i\ll0$.
\end{itemize}

We further define, for $i\in \Z$, 
\begin{equation}
\label{definition_of_d}
    d^i (V):= \dim_{\FF}\  (\mathrm{Fil}^i D_\dR(V)/\mathrm{Fil}^{i+1}D_\dR(V)), \quad d^-(V):= \sum_{i<0} i d^i(V) , \quad d(V):= \sum_{i\in \Z} i d^i(V). 
\end{equation}
\begin{definition}\label{panchishkin}\textbf{Panchishkin condition:}  Consider a finite dimensional $p$-adic representation $V$ of $G_{\Q_p}$, which is de Rham. Then  $V$ satisfies the Panchishkin condition if there exists an exact sequence of $p$-adic representations of $G_{\Q_p}$
\begin{equation}
0 \lra V^+ \lra V \lra V^-\lra 0
\end{equation}
so that $V^+\otimes \C_p \cong \oplus_{i>0} \C_p(i)^{\oplus k_i}$ and $V^-\otimes \C_p \cong  \oplus_{i\leq 0} \C_p(i)^{\oplus k_i}$ as $G_{\Q_p}$-modules, where $k_i \in \Z_{\geq 0}$. 
\end{definition}

\begin{definition}\label{p-ordinary}
    \textbf{Ordinary:} A finite dimensional $p$-adic representation $V$ of $G_{\Q_p}$ is called ordinary (in the sense of Greenberg \cite{Greenberg_Iwasawa_p_adic_representation}) if there exists a filtration $F^i V$ of $\Q_p$-subspaces of $V$ $(i\in \Z)$
with the properties : 
\begin{itemize}
    \item $F^{i+1}V\subseteq F^iV$, $F^iV=V$ for $i\ll 0$ and $F^iV=0$ for $i\gg 0$. 
    \item  $F^i V$ is invariant for the action of $G_{\Q_p}$. Further, the inertia subgroup $I_p$ of $G_{\Q_p}$ acts on $F^iV/F^{i+1}V$ by $\chi_p^i$. 
\end{itemize}
\end{definition}
Note that an ordinary representation of $G_{\Q_p}$ satisfies the Panchishkin condition.

\n \textbf{Hypothesis (Mov):}\label{Hyp:Mov}
A finite dimensional $p$-adic Galois representation $V$ of $G_\Q$ is said to satisfy \hyperref[Hyp:Mov]{(Mov)} if all of the following hold.
\begin{enumerate}
    \item  $V$ is unramified outside finitely many primes of $\Q$.
    \item $V$ is de Rham at $p$, i.e., $V|_{G_{\Q_p}}$ is  a de Rham representation.
    \item The coefficient field of $V$ is a finite extension $\mathcal F$ of $\Q_p$.
\end{enumerate}
 Throughout this article, we will assume that all the Galois  representations $V$ of $G_\Q$ that appear satisfy the hypothesis \hyperref[Hyp:Mov]{(Mov)}. Note that by the Fontaine–Mazur conjecture, it is expected from $(1)$ and $(2)$ that $V$ comes from a motive defined over a number field, and in that case $(1)$ and $(2)$ together implies $(3)$.

\begin{definition}\label{ssd}
    \textbf{Symplectic self-dual:} For $G\in\{G_\Q, G_{\Q_\ell}\}$, a finite dimensional $p$-adic representation $V$ of $G$ is called symplectic self-dual if there exists a non-degenerate, skew-symmetric and $G$-equivariant pairing $ \langle\ ,\ \rangle: V \times V \lra \FF(1)$ such that $V\cong V^*(1)= \Hom_{\FF} (V, \FF(1))$,  where $\FF(1):= \FF \otimes \chi_p$. 
    
    A $G$-stable lattice $T$ of $V$ is called symplectic self-dual if the $G$-equivariant, non-degenerate, skew-symmetric  pairing on $V$ induces a $G$-equivariant perfect pairing  $\langle\ ,\ \rangle: T\times T \lra \OO(1)$.
\end{definition}



\begin{remark}\label{self-dual_lattice}
Let $V$ be a symplectic self-dual finite-dimensional $p$-adic representation of $G_\Q$ and let $T$ be a $G_\Q$-stable lattice of $V$. If either  $\dim_\FF\ V =2$ or if the residual representation of $V$ is an irreducible $G_\Q$-module, then
with respect to a multiple of $\langle\ ,\ \rangle$ by a suitable scalar $a\in \FF ^\times$,  $T$  is also symplectic self-dual (see \cite[\S 2.5]{Ne}). In general, after replacing $\FF$ by  a finite extension $\FF'$ and replacing $V$ by an equivalent symplectic self-dual representation, one obtains a    $G_\Q$-stable symplectic self-dual lattice    (see \cite[Theorem 4.8(ii)]{bhkt}).
\end{remark}

\subsection{Weil-Deligne representation and the local root numbers}
\label{Weil-Deligne representation and local root numbers}  
For a finite dimensional $\Q_p$-vector space $V$,  let
$\rho:G_\Q\longrightarrow \mathrm{Aut}(V)$
 be a $p$-adic Galois representation of $G_\Q$ satisfying hypothesis \hyperref[Hyp:Mov]{(Mov)}.  
 Let $\ell \neq p$ be  primes of $\Q$ and consider $V$ as a representation of $G_{\Q_\ell}$. We fix an isomorphism of $\C_p $ with $\C$. 
 By Grothendieck's monodromy theorem, inertia subgroup $I_\ell$ acts quasi-unipotently, so after passing to an open subgroup of $I_\ell$, one can write
$\rho(\sigma)=\exp\bigl(t_p(\sigma)N\bigr)$
for a nilpotent operator $N$, where $t_p : I_\ell  \lra \Z_p$ is the tame character. The Weil group $W_{\Q_\ell}$  is the subgroup of $G_{\Q_\ell}$ consisting of elements whose image in $\operatorname{Gal}(\overline{\F}_\ell/\F_\ell)\simeq \widehat{\mathbb Z}$ is an integral power of Frobenius. It fits into the exact sequence
$$1\longrightarrow I_\ell\longrightarrow W_{\Q_\ell}\longrightarrow \mathbb Z\longrightarrow 1.$$
Pick a lift $\phi$  of $\Frob_\ell$ to $W_{\Q_\ell}$. The associated Weil-Deligne
representation is a pair $\operatorname{WD}(V)=(r,N)$,
where $r(w)=\rho(w)\exp\bigl(-t_p(\sigma)N\bigr)$, for $w=\phi^n\sigma \in W_{\Q_\ell}$. The representation $r$ 
has finite inertia image  and it satisfies $r(w)Nr(w)^{-1}=|w|\,N$
for all $w\in W_{\Q_\ell}$, where $|*|$ is a normalized valuation on $W_{\Q_\ell}$ via the reciprocity isomorphism (for more details, see \cite{Deligne_73}, \cite[\S 4]{Rohrlich_weil_deligne_group}, \cite{taylor}).

Now let  $\ell=p$, and $V$ be a de Rham   representation $G_{\Q_p}$  defined over a finite extension $\FF/\Q_p$.   Fix a field $E\supset\Q_p^{\ur}$ such that there is an embedding $\tau:\FF\lra E$.
 Then define the Weil-Deligne representation $\WD(V)$ associated to $V$ (see \cite[\S 3.1]{Ne3}), by 
 \begin{equation*}
     \WD(V):= D_\pst(V)\otimes_{\Q_p^\ur \otimes_{\Q_p}\FF,id\otimes\tau }E.
 \end{equation*}
 There exists an $\FF$-linear Frobenius operator $\varphi:D_\pst(V)\lra D_\pst(V)$ and a $\Q_p^\ur \otimes_{\Q_p}\FF$-linear nilpotent endomorphism $N:D_\pst(V)\lra D_\pst(V)$, called monodromy operator,
satisfying $N\varphi=p\varphi N$. The action of the monodromy operator  $N$ is extended to  $\WD(V)$  via $N=N\otimes id $ (see \cite[\S 3.2.1]{Ne3}).

 Deligne \cite{Deligne_73} defined the local epsilon factor at any finite prime $\ell$ (including $\ell=p$), $\epsilon_\ell(\WD(V),\psi,dx)$ associated to the Weil-Deligne representation $\WD(V)$, a non-trivial additive character $\psi : \Q_\ell \rightarrow \C^\times$ and a Haar measure $dx$ on $\Q_\ell$. The local root number of $V$ at $\ell$ is defined by 
\begin{equation}
\label{Rohrlich_root}
    W(V,\psi,\ell) = \mathrm{sign}(\epsilon_\ell(\WD(V),\psi,dx)):=\frac{\epsilon_\ell(\WD(V),\psi,dx)}{|\epsilon_\ell(\WD(V),\psi,dx)|}.
\end{equation}
 Note that, if $V|_{G_{\Q_\ell}}$ is symplectic self-dual representation then $W(V,\ \psi,\ \ell)$ is independent of $\psi$, $dx$ and its value is $\pm1$ (see \cite[\S12]{Rohrlich_weil_deligne_group}), and we will denote it by $W_\ell (V)$.  Since $V$ is de Rham at $p$, we can define the local root number at an Archimedean prime as in  \cite[\S 4.1.3]{Ne3}.
 If $V$ comes from a motive, then this definition of local root number at an Archimedean prime  matches with the corresponding definition obtained using the Hodge-Tate decomposition of the de Rham realisation of $V$ (see
\cite[\S 4.4]{Tate}, \cite[\S 5.3]{Deligne_79}).
Now  define the global root number  of  $V$ by 
\begin{equation}\label{global_root}
W(V):=\prod_{\text{all places } \ell} W_\ell (V).
\end{equation}

 For a finite prime $\ell\neq p$,  Deligne also defines \cite[(5.1)]{Deligne_73} the modified local epsilon factor 
\begin{equation}
    \epsilon_0(\WD(V),\psi, dx, \ell):= \epsilon_\ell(\WD(V),\psi, dx)\det\left(-\mathrm{Frob}_\ell \vert V^{I_\ell}\right).
\end{equation}


Let $V$ be as in \S\ref{Weil-Deligne representation and local root numbers} above and
let the coefficient field of  $V$ be $\FF$.  Suppose $T$ is a $G_\Q$-stable $\OO$-lattice of $V$. Let $T$ be a $G_\Q$-stable symplectic self-dual lattice of $V$ (see Remark \ref{self-dual_lattice}) and
 $\bar{T}:=T/\pi T$ be the  residual representation associated to $T$. 
The semisimplification of $\bar{T}$ is independent of the choice of the lattice $T$. Further, if $\bar{T}$ is an irreducible $G_\Q$-module, then $\bar{T}$ is also independent of the choice of $T$. 
 
It follows from \cite[Theorem 6.5(c)]{Deligne_73}, that for $\ell\neq p$, $\epsilon_0(\WD(V), \psi, dx, \ell)\pmod{\pi}$ depends only on the semisimplification $\bar{T}^\mathrm{ss}$ of the residual representation $\bar{T}$. Then, using \eqref{Rohrlich_root}, we obtain: 

\begin{lemma}\label{root_number_V_1_V_2} 
Let $V_1,\ V_2$ be two finite-dimensional $p$-adic Galois representations of $G_{\Q_\ell}$ with coefficient field $\FF$. Also assume that both $V_1$ and $V_2$ are symplectic self-dual Galois representations and their residual representations $\bar{T}_1\cong \bar{T}_2$ as $G_{\Q_\ell}$-modules. Then for a finite prime $\ell \neq p$, we have 
    \begin{equation}
    \label{root_number_relation}
        \frac{W_{\ell}(V_1)}{W_\ell(V_2)} \equiv  \frac{\det \left(-\mathrm{Frob}_\ell \vert (V_2)^{I_\ell}\right)}{\det \left(-\mathrm{Frob}_\ell \vert (V_1)^{I_\ell}\right)} \pmod{\pi}.\qed
    \end{equation} 
\end{lemma}

We end this section with a definition  to be used later.
\begin{definition}
   \textbf{Purity:}\label{purity}  A finite dimensional $p$-adic representation $V$ of $G_{\Q_\ell}$ is called pure of weight $m\in \Z$ if the associated Weil-Deligne representation $\WD(V)$ is pure of weight $m$; see \cite[\S1.4]{Ne3}. If $V$ is self-dual, i.e., $V\cong V^*(1)$ and is pure of weight $m\in \Z$, then $m=-1$. 
\end{definition}

\section{A generalization of a
theorem of Coates–Greenberg}\label{Kummer_theory}

Let $V$ a finite dimensional $\FF$-vector space endowed with a continuous $G_{\Q_p}$ action.
Assume that $V$ satisfies the  Panchishkin condition at $p$ (Definition \ref{panchishkin}). 
 Let $T$ be a $G_{\Q_p}$-stable  $\OO$-lattice of $V$ and put $A:= V/T$ and let $A^+$ be the image of $V^+$ under the natural map $V\lra A$ and set $A^-:= A/A^+$.

Let $L$ be a number field,  $v\mid p$ be a finite prime of $L$ and  $I_v$ denote the inertia subgroup of $G_{L_v}$ and  $L_v^{\ur}$ denote the maximal unramified extension of $L_v$.
For a number field $L$, $L_\infty$ will denote the cyclotomic $\mathbb{Z}_p$-extension of $L$ with $\Gamma_L:= \Gal(L_\infty /L)\cong \mathbb{Z}_p$ and let $L\subset L_n\subset L_\infty$ with $\Gal(L_n/L)\cong \Z/p^n \Z$. 
 Let $L_{\infty,v}$ be the cyclotomic $\Z_p$-extension of $L_v$ and $L_{n,v}$ be its $n$-th layer. 


We define the Bloch-Kato vector space local condition at $v$  by 
\begin{equation*}
H^\BK (L_v,V):= \Ker\left(H^1(L_v,V)\rightarrow H^1(L_v,V \otimes_{\Q_p}B_{\crys} )\right).
\end{equation*}
 We have a natural map $H^1(L_v,V) \xrightarrow{p_*} H^1(L_v,A)$ induced by the natural projection map $V \xrightarrow{p} A$. Then define  the 
Bloch-Kato local condition at $v$, 
$H ^{\BK}(L_v, A): = p_{*}(H_v^\BK (L,V))$.
 Also define the strict local condition of $A$ over $L$ at $v$  by 
\begin{equation*}
     H ^\str (L_v, A):= \Ker\left(H^1(L_v,A)\rightarrow H^1(L_v,A^-)\right)
\end{equation*}

Let $B$ be an abelian variety defined over $L$ with ordinary reduction at the  prime $v\mid p$. Then using the Kummer map, Coates-Greenberg \cite[(4.9)]{CG} showed that $H^\BK(\mathcal{L}_w,B)\subset H^\str(\mathcal{L}_w,B)$ for any Galois extension $\mathcal{L}/L$ and any prime $w\mid v \mid p$ of $\mathcal{L}$. Moreover, as a key result in the article (\cite[Prop. 4.3]{CG}), it is  proved that for any {\it deeply  ramified extension} $\mathcal{L}/L$, the above two subgroups are equal. In particular, if we consider the cyclotomic $\Z_p$-extension $L_\infty$ of $L$, then it follows from \cite[Prop. 4.3]{CG} that for any prime $v\mid p$ of $L$,
\begin{equation*}
H^\BK(L_{\infty,v},B)=H^\str(L_{\infty,v},B).
\end{equation*}

Nekov\'{a}\v{r} has also  compared  the Bloch-Kato and strict local conditions for $p$-ordinary abelian varieties  in \cite[Prop. 9.6.7.5(i)]{Ne2}.
We now extend the results of Coates-Greenberg  for  more general $p$-adic Galois representations which satisfy the Panchishkin condition at $p$. In particular,  our result includes representations  which  exhibit exceptional zero phenomenon (cf. \cite{Gr3}). 

Let $V$, $T$, $A$ be as   above. 
Our goal is to establish the equality between $H^\str(L_{\infty,v},A):=\varinjlim_n H^\str(L_{n,v},A)$ and  $H^\BK(L_{\infty,v},A):=\varinjlim_n H^\BK(L_{n,v},A)$.

\begin{lemma}\label{Sen_theory}
 For a number field $L$ and for a prime $v|p$ of $L$,    Let $\rho: G_{L_v} \lra \mathrm{Aut}_{\Q_p} (V)$ be a $p$-adic Galois representation which is semistable.  Then $H^0(G_{L_v^\ur (\mu_{p^\infty})},V) \sim \bigoplus_{k=1}^{r} \chi_p ^{i_k}$ as a $\Gal(L_v^\ur (\mu_{p^\infty})/L_v^\ur )$-module.
\end{lemma}

\n \textbf{Proof:} Put $\tilde{L}_v=L_v^\ur (\mu_{p^\infty})$. Set $X=V^{G_{\tilde{L}_v}}$ and  $\widetilde{X}= X\otimes_{\Q_p} \C_p$. Then by Sen's theory \cite[Theorem 3.31]{FO}, there exists an open subgroup $\Gamma_0=\langle\gamma_0 \rangle$ of $\Gamma = \Gal(L_{v}^\mathrm{ur} (\mu_{p^\infty})/ L_{v}^\mathrm{ur})$ of finite index such that after fixing a basis of $\widetilde{X}$, $\rho : \Gamma_0 \lra \mathrm{Aut}(\widetilde{X})$ is given by $\rho(\gamma_0)= \exp (\Phi \log (\chi_p (\gamma_0)))$, where $\Phi$ is the matrix representation of the Sen's operator. This $\Phi$
 has coefficients in $L_v^{\ur}$. On the other hand, by \cite[Prop. 5.5]{FO}, $\Phi$ is semisimple and its eigenvalues are in $\Z$, say $i_1,\dots, i_r$ (need not be distinct). Hence we have
\begin{small}\begin{equation*}
     \rho(\gamma_0) \sim \begin{pmatrix}
\chi_p(\gamma_0)^{i_1} & 0 & \cdots & 0 \\
0 & \chi_p(\gamma_0)^{i_2} & \cdots & 0 \\
\vdots & \vdots & \ddots & \vdots \\
0 & 0 & \cdots & \chi_p(\gamma_0)^{i_r}
\end{pmatrix}.
\end{equation*}\end{small}
Since $\chi_p (\gamma_0)\in \Z_p^\times$, there exists a basis over $L_{v}$ such that $\rho|_{\Gamma_0}$ is diagonal and 
\begin{equation*}
    \rho|_{\Gamma_0} \sim \bigoplus_{k=1}^{r} \chi_p^{i_k}.
\end{equation*}
Consider an eigenvalue $i_k$ of $\Phi$, and corresponding to this, there exists a $\Gamma_0$-stable subspace of $X$, say $X_k$, where $\rho|_{\Gamma_0}$ acts via $\chi_p^{i_k}$.
Then as a $\Gamma_0$-module, $X$ can be written as $X\cong \bigoplus_{k=1}^{r} X_k$. 
As $\Gal(\tilde{L}_v/L_v)$ is abelian, $X_k$ is a $\Gal(\tilde{L}_v/L_v)$-subrepresentation of $X$ and hence semistable. Therefore $X_k\otimes  \chi_p^{-i_k}$ is also semistable.
 Note that $\Gamma_0$ acts trivially on $X_k\otimes  \chi_p^{-i_k}$ and hence action of $I_{L_v}$  on $X_k\otimes  \chi_p^{-i_k}$ factors through a finite quotient of $\Gamma$.  Further, it follows from   \cite[Prop. 7.17]{FO} that $I_{L_{v}}$ acts trivially on $X_k\otimes  \chi_p^{-i_k}$ i.e., $I_{L_{v}}$ (and $\Gamma$) acts on $X_k$ via $\chi_p^{i_k}$. The result follows from this.  \qed

\medskip

Recall $v$ be a prime of a number field $L$ dividing $p$. Let $V$ be a $p$-adic semistable Galois representation of $G_{L_v}$ defined over $\FF$. Fix a field $E\supset\Q_p^{\ur}$ such that there is an embedding $\tau:\FF\lra E$.
 Then as in \S \ref{Weil-Deligne representation and local root numbers},  the Weil-Deligne representation associated to $V$ can be defined by (see also \cite[\S 3.1]{Ne3}), $\WD(V):= D_\st(V)\otimes_{L_0\otimes_{\Q_p}\FF,id\otimes\tau }E$,
 where $L_0= L_v\cap \Q_p^\ur$.

\begin{proposition}\label{trivial_action_gamm_v}
Let $L$ be a number field and $v\mid p$ be a prime such that $L_v= L_v^\ur \cap L_v (\mu_{p^\infty})$.
Let $V$ be a p-adic Galois representation of $G_{L_v}$,   which is self-dual symplectic, semistable, and pure (of weight $-1$) satisfying the Panchishkin condition.  Then $\Gamma_v:=\Gal(L_{\infty,v}/L_v)$ acts on $H^0(L_{\infty,v},V^-)$ trivially.
\end{proposition}
\begin{proof}
 As $V$ is semistable, $H^0(L_v(\mu_{p^\infty}),V^-)$  is also semistable. Assume that $H^0(L_{\infty,v},V^-)\neq 0$.
Since, $ \Gal(L_v^{\ur}(\mu_{p^\infty})/L_v(\mu_{p^\infty}))$ acts trivially on $H^0(L_v(\mu_{p^\infty}),V^-)$ and $\Gal(L_v^{\ur}(\mu_{p^\infty})/L_v)\cong \Gal(L_v^\ur/L_v)\times \Gal(L_v(\mu_{p^\infty})/L_v)$, it
follows from Lemma \ref{Sen_theory}, that $H^0(L_v(\mu_{p^\infty}),V^-)$  is isomorphic to $\oplus_{r=1}^s \chi_p^{i_r}$ as a $\Gal(L_v(\mu_{p^\infty})/L_v)$-module for some $s\geq 1$ and $i_r\in \Z$. Using the Panchishkin condition, we get that $i_r\leq 0$ for all $r$. Then $H^0(L_{\infty,v},V^-)\simeq \oplus_{d|i_r}\chi_p^{i_r}$ where $d:= [L_v(\mu_p):L_v]$. It suffices to show that $i_r=0$ whenever $d\mid i_r$. On the contrary, assume that   $t:=i_r <0$ for some such $r$. Then we have an exact sequence 
\[ 0 \lra V^+ \lra U \lra \mathcal{F}(t) \lra 0,\]
where $\FF(t):= \Q_p(t)\otimes \FF$ and $U\subset V$ is the preimage of $\FF(t)\subset V^-$ via $V\twoheadrightarrow  V^-$. Note that $V^+$, $U$ and $\FF(t)$ are again semistable representations. This gives us the following exact sequence of $(\varphi,N)$ modules:
\[ 0 \lra D_\st(V^+) \lra D_\st (U) \stackrel{f}{\lra} D_\st (\mathcal{F}(t)) \lra 0. \]
Applying the monodromy operator $N$, we get the following exact sequence 
\[ 0 \lra D_\st (V^+)^{N=0} \lra D_\st (U)^{N=0} \stackrel{f}{\lra} D_\st(\mathcal{F}(t))^{N=0} \stackrel{\delta}{\lra} D_\st(V^+)/ND_\st(V^+)  .\]
Since $\mathcal{F}(t)$ is crystalline, $D_\st(\FF(t))=D_\st(\FF(t))^{N=0}$.
Let $0\neq x\in D_\st(\mathcal{F}(t))$ and $\bar{x}\in D_\st(U)$ such that $f(\bar x)=x$. Since $Nx=0$, $N\bar x\in D_\st(V^+)$ and the connecting homomorphism is given by $\delta(x)=N\bar x \pmod{ND_\st(V^+)}$. 
We have $\varphi(x)=p^{-t }x$ and $N\varphi=p\varphi N$. Therefore $\delta(\varphi(x))=N\varphi(\bar x)=p\varphi N\bar x$ and hence $\varphi\delta(x)= p^{{-t }-1}\delta (x)$.

Consider the slope decomposition $D_\st(V^+)/ND_\st(V^+)\cong \oplus_{\alpha_i} D(\alpha_i)$ as in \cite[Lemma 8.1.11]{BC}, where $D(\alpha_i)$ are subobjects (called isoclinic) with Frobenius slope $\alpha_i$ such that $\alpha_1< \alpha_2<\cdots < \alpha_n$. 
Assuming $\delta$ is a non-zero map, from the equality $\varphi\delta(x)= p^{{-t }-1}\delta (x)$ we get that $\alpha_i={-t}-1$ 
for some $i$. Therefore, we get a surjection $D_\st(V^+) \twoheadrightarrow D_\st(V^+)/ND_\st(V^+) \twoheadrightarrow D({-t}-1)$. Let $t_N(D({-t}-1))=({-t}-1)\dim D({-t}-1))$ be the Newton number of $D({-t}-1)$ as defined in \cite[page 106]{BC}.   

Since $V$ is semistable and $V^+$ is a $G_{L_v}$-stable
subrepresentation, $V^+$ is semistable. Consequently,
$D_{\mathrm{st}}(V^+)$, equipped with its Hodge filtration,
is a weakly admissible filtered $(\varphi,N)$-module in the sense of \cite[Definition 8.2.9]{BC} (see also \cite[Theorem 9.3.4]{BC}). Since it  is a quotient in the category of filtered $(\varphi,N)$-modules, weak admissibility, together with $t_H\bigl(D_{\mathrm{st}}(V^+)\bigr)
=
t_N\bigl(D_{\mathrm{st}}(V^+)\bigr)$
implies that
$t_N\bigl(D(\alpha_i)\bigr)
\leq
t_H\bigl(D(\alpha_i)\bigr)$ for every $i$.
In particular, we have $t_N(D({-t}-1))\leq t_H(D({-t}-1))$,
where $t_H(D({-t}-1))$ is the Hodge number of $D({-t}-1))$ (see \cite[Definition 8.1.1]{BC}). Note that $t_H(-)$ is the same as $d(-)$ as defined in \eqref{definition_of_d}. Since $V$ satisfy the Panchishkin condition, $ \mathrm{Fil}^0 D_{\dR}(V^+)=0$, where $\mathrm{Fil}^i D_{\dR}(V^+)$ is as defined in \eqref{definition_of_d}. 
Recall $L_0=L_v \cap \Q_p^\ur $.  Using the isomorphism $D_{\st}(V^+)\otimes_{L_0} L_v \cong D_{\dR}(V^+)$ (see \cite[\S 5.6.3]{Fo2}) and the surjection  $D_\st(V^+)  \twoheadrightarrow D({-t}-1)$, we see that the induced filtration on the quotient $D({-t}-1)$ satisfy $\mathrm {Fil}^0 D({-t}-1)\otimes_{L_0} L_v =0$. Therefore $(-t -1)\dim D({-t}-1)) =t_N(D({-t}-1))\leq t_H(D({-t}-1))<0$, which is a contradiction as $t<0$. Therefore, $\delta$ must be zero, and we have an exact sequence 
 \begin{equation}\label{sur1}
 0 \lra D_\st(V^+)^{N=0} \lra D_\st(U)^{N=0} \stackrel{f}{\lra} D_\st(\mathcal{F}(t) )\lra 0 .  \end{equation} 
Recall the increasing monodromy filtration $M_{\bullet}$ of the Weil-Deligne representation  $\WD(V)$ (\cite[1.3.1]{Ne3}) 
\[ M_n\WD(V) := \sum_{i-j=n+1} \Ker(N^i)\cap \operatorname{Im}(N^j) \]  satisfies 
$N(M_n \WD(V))\subseteq M_{n-2}\WD(V) $  and, for every $n\geq 0$, \[ N^n:\operatorname{Gr}^M_n \WD(V) \xrightarrow{\sim} \operatorname{Gr}^M_{-n}\WD(V). \]
From the definition, we see that  $\Ker(N)=\Ker(N)\cap\operatorname{Im}(N^0)\subseteq M_0\WD(V)$. Since $D_\st(E)^{N=0}\subset \Ker(N)$ and $V$ is pure, we deduce from the exact sequence \eqref{sur1} that there exists an $i\leq 0$ such that, Weil weight of $\operatorname{Gr}^M_i\WD(V)$ is $-2t>0$. 
However, using $V$ is pure of weight $-1$, Weil weight of $\operatorname{Gr}^M_i\WD(V) $ is $-1+i$ (see \cite[\S 1.4.3]{Ne3}), which is negative. Thus we arrive at a  contradiction.  Therefore $i_r=0$ for all $d\mid i_r$ and hence $\Gamma_v$ acts trivially on $H^0(L_{\infty,v},V^-)$. 
   \end{proof}

\begin{corollary}\label{mod_2_Qcyc_Q}
Let $V$ be a finite dimensional  $\mathcal{F}$ vector space endowed with a continuous 
 $G_{L_v}$ action. Assume that $V$ is symplectic self-dual representation which is potentially semistable,  is pure (of weight $-1$) and satisfy the Panchishkin condition.  Then 
\begin{equation*}
    \dim_\mathcal{F} H^0(L_{\infty,v},V^-)\equiv \dim_\mathcal{F}  H^0(L_v,V^-) \pmod 2.
\end{equation*}
\end{corollary}
\begin{proof}
Let $K_v$ be a finite Galois extension of $L_v$ such that $V$ is semistable over $K_v$. Extending to a finite extension, if needed, we can assume that  $K_v= K_v^\ur \cap K_v(\mu_{p^\infty})$. Consider the slope decomposition of $D_{\mathrm{st}}(V^+)= \bigoplus_{i=1}^r
D(\alpha_i)$, with isoclinics $D(\alpha_i)$ and $\alpha_i<\alpha_{i+1}$ for all $i$. Since multiplication by $N$ reduces slope, $\bigoplus_{i<r} D(\alpha_i)$ is a $(\phi,N)$ submodule. In particular, $D_{\mathrm{st}}(V^+)/(\bigoplus_{i<r} D(\alpha_i)) =D(\alpha_r)$ is a $(\varphi,N)$-module quotient of a weakly admissible $(\varphi,N)$-module $D_{\mathrm{st}}(V^+)$, and hence $\alpha_r \dim D(\alpha_r) = t_N (D(\alpha_r))\leq  t_H(D(\alpha_r))$.
 Let $K_0$ be the maximal unramified extension of $\Q_p$ in $K_v$.
Using the isomorphism $D_\mathrm{st}(V^+)\otimes_{K_0}K_v\cong D_\mathrm{dR}(V^+)$ and the Panchishkin condition $\mathrm{Fil}^0 D_{\mathrm{dR}} (V^+)=0$,  we have $t_H(D(\alpha_r))<0$ and consequently, $\alpha_r <0$.  In other words, slopes of eigenvalues of $\varphi$ are always negative. 
Hence $D_{\mathrm{st}}(V^+)^{\varphi^h -1}=D_{\mathrm{st}}(V^+)/{(\varphi^h -1}) =0$. Therefore 
\begin{equation}\label{frob1}
D_{\mathrm{st}}(V)^{\varphi^h -1}=D_{\mathrm{st}}(V^-)^{\varphi^h -1},   \end{equation}
where $h:=[K_0:\Q_p]$.

Recall $E$ is an extension of $\Q_p^{ur}$ and we have an embedding $\tau : \mathcal{F}_v\hookrightarrow E$. Then the Frobenius  element $\Frob_v$ in the Weil-Deligne group $W_{K_v}$ acts on $WD_\tau(V)=D_{\mathrm{st}}(V)
\otimes_{K_0\otimes_{\mathbf Q_p}L,\;\mathrm{id}\otimes\tau}E$ via $\rho(\Frob_v)=\varphi^h\otimes\mathrm{id}$ (see \cite[Example 3.2.2 (2)]{Ne3}).  We also assume that all the representations of $\Gal(K_v/L_v)$ are defined over $E$.  The relation between the action of $\Frob_v$ and $\varphi$ gives us
\begin{equation}\label{frob2}
WD_\tau(V)^{\Frob_v -1}=D_{\mathrm{st}}(V)^{\varphi^h -1}\otimes_{\mathcal{F}_v,\tau} E.    \end{equation}
From  \cite[Proposition 3.3.2 (1)]{Ne3}, equation \eqref{frob1} and equation \eqref{frob2} we get
\begin{equation}\label{WD-frob}
WD_\tau(V)^{\Frob_v -1}=H^0(K_v,V^-)\otimes_{\mathcal F_v,\tau} E. \end{equation}

Recall that we have an increasing monodromy filtration on $WD_\tau(V)$ with graded components  $\operatorname{Gr}_r^M \WD(V)$, which are strictly pure of weight $r-1$. In particular, $\operatorname{Gr}_r^M \WD(V)^{\Frob_v-1}=0$ unless $r=1$. Then from equation \eqref{WD-frob}, we deduce $WD_\tau(V)^{\Frob_v -1}=\operatorname{Gr}_1^M \WD(V)^{\Frob_v -1}=H^0(K_v,WD_\tau(V))$.
As $V$ is semistable, the inertia subgroup acts trivially on $\operatorname{Gr}_1^M \WD(V)$ and therefore,
\begin{equation}\label{WD-frob1}
H^0(K_v,\operatorname{Gr}_1^M \WD(V))=\operatorname{Gr}_1^M \WD(V)^{\Frob_v -1}=H^0(K_v,V^-)\otimes_{\mathcal F_v,\tau} E. \end{equation}
 
Let $\eta$ be a non-trivial sub-representation of $\Gal(K_v/L_v)$ in $H^0(K_v,V^-)$.  Since  there exist a non-degenerate orthogonal $G_{L_v}$-equivariant pairing on $\operatorname{Gr}^M_1 \WD(V)$ (\cite[1.3.5]{Ne3}), the multiplicities of $\eta$ and its dual $\eta^\vee$ are the same in $\operatorname{Gr}^M_1 \WD(V)$. Moreover,  $G_{K_v}$ is contained in the kernel of $\eta$, and thus  $\eta^\vee \subset \operatorname{Gr}^M_1 \WD(V)^{\Frob_v -1}=H^0(K_v,\operatorname{Gr}^M_1 \WD(V))=H^0(K_v,V^-)\otimes_{\mathcal F_v,\tau} E $. 

Now consider a special case where $\eta$ is a character of $\Gal(L_{n,v}/L_v)$ with $L_{n,v}:=K_v\cap L_{\infty,v}$. Then both $\eta$ and $\eta^{-1}$ occur in $H^0(L_{n,v},V^-)\otimes_{\mathcal F_v,\tau} E$ with the same multiplicity. Note that $V$ is semistable over $K_v$. Then by Proposition  \ref{trivial_action_gamm_v}, $\Gal(K_{\infty,v}/K_v)$ acts on $H^0(K_{\infty,v},V^-)$ trivially and hence $\Gal(L_{\infty,v}/L_{n,v})$ acts on $H^0(L_{\infty,v},V^-)$ trivially as well. We deduce that, for every character $\eta$ of $\Gal(L_{\infty,v}/L_v)$, the multiplicity of $\eta$ and $\eta^{-1}$ are the same in $H^0(L_{n,v},V^-)\otimes_{\mathcal F_v,\tau} E$. 
Now the corollary follows from the assumption that $p$ is an odd prime. 
\end{proof}

\subsection{Almost Cartesian condition at $p$} Panchishkin condition for $V$ at $p$ induces a Panchishkin condition at $p$ for $T$, i.e., there exists an 
 exact sequence of $G_{L_v}$-modules $0\to T^+\to T\to T^-\to0$, where $T^-$ is the image of $T$ in $V^-$ under the natural map $V\twoheadrightarrow V^-$ and  $T^+:= \Ker (T\lra T^-)$. Define
$H^{\str}(L_{n,v},T)
:=\mathrm{Im}\bigl(H^1(L_{n,v},T^+)\to H^1(L_{n,v},T)\bigr)$,
and put
\begin{equation*}
    Q_n:=H^1(L_{n,v},T)/H^{\str}(L_{n,v},T).
\end{equation*}

\begin{definition}[Bounded Cartesian defect at $p$]
We say that the strict local condition has
\emph{uniformly bounded Cartesian defect at $p$ along $L_{\infty,v}/L_v$}
if there exists a constant $C$, independent of $n$, such that for all $n$,
\[
\operatorname{length}_{\mathcal O}\bigl((Q_n)_{\mathrm{tor}}\bigr)
\leq C.
\]

Equivalently, the groups $(Q_n)_{\mathrm{tor}}$ are finite of uniformly bounded order independent of $n$. \qed
\end{definition}
The notion of a Cartesian local condition was introduced by Mazur-Rubin (\cite[Definition 1.1.4]{MR1}). If $(Q_n)_{\mathrm{tor}}=0$, then  the strict local condition is indeed Cartesian (see \cite[Lemma 3.7.1(i)]{MR1}). Our condition can be considered as a variant of the Mazur-Rubin condition.  Instead of vanishing of $(Q_n)_{\mathrm{tor}}$, we are seeking it to be uniformly bounded in the cyclotomic tower.

The exact sequence $0\to A^+\to A\to A^-\to0$ gives a connecting map
$\partial_n:H^0(L_{n,v},A^-)\to H^1(L_{n,v},A^+)$, hence induces a map
$\partial_n^{\mathrm{red}}:
\frac{H^0(L_{n,v},A^-)}{H^0(L_{n,v},A^-)_\mathrm{div}}
\to \frac{H^1(L_{n,v},A^+)}{H^1(L_{n,v},A^+)_\mathrm{div}}$.
Define
\begin{equation}
D_n:=
\ker\left(
\frac{H^0(L_{n,v},A^-)}{H^0(L_{n,v},A^-)_\mathrm{div}}
\xrightarrow{\;\partial_n^{\mathrm{red}}\;}
\frac{H^1(L_{n,v},A^+)}{H^1(L_{n,v},A^+)_\mathrm{div}}
\right).
\end{equation}
The exact sequence $0\to T^+\to T\to T^-\to0$ gives a connecting homomorphism $d_n:H^1(L_{n,v},T^-)\to H^2(L_{n,v},T^+)$ and  we obtain
\begin{equation}
Q_n
\simeq
\ker\!\left(
H^1(L_{n,v},T^-)\xrightarrow{d_n}H^2(L_{n,v},T^+)
\right).
\end{equation}
There are
canonical identifications
\begin{equation}\label{identification}
\frac{H^0(L_{n,v},A^-)}{H^0(L_{n,v},A^-)_\mathrm{div}}
\simeq H^1(L_{n,v},T^-)_{\tor},
\qquad
\frac{H^1(L_{n,v},A^+)}{H^1(L_{n,v},A^+)_\mathrm{div}}
\simeq H^2(L_{n,v},T^+)_{\tor}.
\end{equation}
Under these identifications, $\partial_n^{\mathrm{red}}$ is the restriction to torsion
of the connecting map $d_n$. More precisely, 
from \cite[Proposition 4.1]{CE}, we get the following anticommutative diagram 
\[
\begin{array}{ccc}
\frac{H^0(L_{n,v},A^-)}{H^0(L_{n,v},A^-)_\mathrm{div}}
&
\xrightarrow{\ {\partial_n^{\mathrm{red}}}}
&
\frac{H^1(L_{n,v},A^+)}{H^1(L_{n,v},A^+)_\mathrm{div}}
\\[2mm]
\Big\downarrow\vcenter{\rlap{$\beta^- \scriptstyle \cong$}}
&&
\Big\downarrow\vcenter{\rlap{$\beta^+ \scriptstyle \cong$}}
\\[2mm]
H^1(L_{n,v},T^-)_{\mathrm{tor}}
&
\xrightarrow{\ d_n\ }
&
H^2(L_{n,v},T^+)_{\mathrm{tor}} 
\end{array}
\]
satisfying $\beta^{+}\circ \partial_n^{\mathrm{red}}
=-d_n\circ \beta^{-},$ and therefore
$D_n\simeq (Q_n)_{\tor}.$

Consequently, we have the following lemma.
\begin{lemma}\label{bdd}
The local strict Selmer condition at $p$ has uniformly bounded Cartesian defect  if and only if there exists $C>0$ (independent of $n$) such that $
\operatorname{length}_{\OO}\ (D_n)<C$ for all $n$. \qed
\end{lemma}

Recall $V^*= \Hom_\FF (V,\FF)$ and further define  $T^*:= \Hom(T, \OO)$ and $A^*:= T^*\otimes_\OO \FF/\OO$. We also define $(T^\pm )^*$ and  $(A^\pm)^*$  analogously.
From \eqref{identification}, we have $\frac{H^1(L_{n,v},A^+)}{H^1(L_{n,v},A^+)_\mathrm{div}}
\cong H^2(L_{n,v}, T^+)_\tor $ and applying Tate duality we obtain $H^2(L_{n,v}, T^+)_\tor \cong \frac{H^0(L_{n,v}, (T^+)^*(1))^\vee}{H^0(L_{n,v}, (T^+)^*(1))^\vee_\mathrm{div}}=  \frac{H^0(L_{n,v}, (A^+)^*(1))^\vee}{H^0(L_{n,v}, (A^+)^*(1))^\vee_\mathrm{div}}$. Further, applying the Panchishkin condition, we get that $(A^+)^*(1)\cong (A^*(1))^-$, and hence $\frac{H^1(L_{n,v},A^+)}{H^1(L_{n,v},A^+)_\mathrm{div}}\cong \frac{H^0(L_{n,v}, (A^+)^*(1))^\vee}{H^0(L_{n,v}, (A^+)^*(1))^\vee_\mathrm{div}}\cong  \frac{H^0(L_{n,v},A^*(1)^-)^\vee}{H^0(L_{n,v},A^*(1)^-)^\vee_\mathrm{div}}$.
Hence we get an exact sequence
\begin{equation}\label{coker}
0\lra D_n \lra \frac{H^0(L_{n,v},A^-)}{H^0(L_{n,v},A^-)_\mathrm{div}}
\longrightarrow
\Big(\frac{H^0(L_{n,v},A^*(1)^-)}{H^0(L_{n,v},A^*(1)^-)_\mathrm{div}}\Big)^\vee \lra \frac{H^{\str}(L_{n,v},A)}{H^{\str}(L_{n,v},A)_\mathrm{div}} \lra 0.
\end{equation}

\begin{proposition}\label{trvial-action_implies_defect} Let $V$ be a Panchishkin representation of $G_{L_v}$.
Suppose that $\Gamma_v=\Gal(L_{\infty,v}/L_v)$ acts trivially on $H^0(L_{\infty,v},V^-)$. Then  the  local strict Selmer condition at $p$ has uniformly bounded Cartesian defect. 
\end{proposition}
\begin{proof}
The exact sequence $0 \lra T^- \lra V^- \lra A^- \lra 0$
gives us
\[
H^0(L_{\infty,v},A^-)_{\mathrm{div}}
=
\operatorname{Im}\left(
H^0(L_{\infty,v},V^-)
\longrightarrow
H^0(L_{\infty,v},A^-)
\right).
\]
Since $\Gamma_v$ acts trivially on $H^0(L_{\infty,v},V^-)$, then the same is true for all 
$\Gamma_{v,n}:=\Gal(L_{\infty,v}/L_{n,v})$, and hence $H^0(L_{n,v},V^-)
=
H^0(L_{\infty,v},V^-)^{\Gamma_{v,n}}
=
H^0(L_{\infty,v},V^-).$
Therefore $H^0(L_{n,v},A^-)_{\mathrm{div}}
=
H^0(L_{\infty,v},A^-)_{\mathrm{div}}.$
Moreover, $H^0(L_{n,v},A^-)
=
H^0(L_{\infty,v},A^-)^{\Gamma_{v,n}}
\subseteq
H^0(L_{\infty,v},A^-),$
and consequently
\[
\frac{H^0(L_{n,v},A^-)}
{H^0(L_{n,v},A^-)_{\mathrm{div}}}
\hookrightarrow
\frac{H^0(L_{\infty,v},A^-)}
{H^0(L_{\infty,v},A^-)_{\mathrm{div}}}.
\]
It follows that 
$\frac{H^0(L_{n,v},A^-)}{H^0(L_{n,v},A^-)_\mathrm{div}}$ is finite and  uniformly bounded by the finiteness of $\frac{H^0(L_{\infty,v},A^-)}{H^0(L_{\infty,v},A^-)_\mathrm{div}}$. In particular,  $D_n$ is finite and uniformly bounded by \eqref{coker}, and the proposition follows from Lemma \ref{bdd}. 
\end{proof}

\begin{theorem}\label{str_bk_at_p}
   Let $V$ be a self-dual symplectic $p$-adic Galois representation of $G_{\Q_p}$ which is de Rham. Also assume that $V$ satisfies the Panchishkin condition and is pure (of weight $-1$). 
    Then the following is true.
    \begin{enumerate}[label=(\alph*)]
        \item The  local strict Selmer condition at $p$ has uniformly bounded Cartesian defect.
        \item $H^{\str}(\Q_{\infty,p},A)=H^\BK(\Q_{\infty,p},A).$
    \end{enumerate}
 \end{theorem}
\begin{proof} As $V$ is de Rham representation of $G_{\Q_p}$, it is potentially semistable by  \cite[Theorem 0.7]{Ber}. So there exists a finite extension $J$ of $\Q_p$ such that $V$ is semistable as a $G_J$-representation. Extending to a finite extension, if needed, we can assume that  $J= J^\ur \cap J(\mu_{p^\infty})$. Let $\Q_{m,p}= J\cap \Q_{\infty,p}$, where $\Gal(\Q_{m,p}/\Q_p)\cong \Z/p^m\Z$ for a natural number $m$. Now $\Gal(J_\infty/J)$ acts trivially on $H^0({J_\infty},V^-)$ by Proposition  \ref{trivial_action_gamm_v}. We have $\Gal({J}_\infty /\Q_{m,p})\cong \Gal(\Q_{\infty,p}/\Q_{m,p})\times \Gal(J/\Q_{m,p})$,  and $\Gal(\Q_{\infty,p}/\Q_{m,p})$ acts trivially on $ H^0(\Q_{\infty,p},V^-)=H^0(J_{\infty},V^-)^{\Gal(J_\infty/\Q_{\infty,p})}$. 
Now the assertion (a) follows from Prop.  \ref{trvial-action_implies_defect}.

For the proof of (b), observe that the self-duality of $V$, (i.e., $V\cong V^*(1))$ and the uniqueness of the Panchishkin filtration imply that $V^*(1)^-\cong V^-$ and hence $A^*(1)^-\cong A^-$. Now from the exact sequence \eqref{coker} and Lemma \ref{bdd}, it follows that $\frac{H^{\str}(\Q_{n,p}, A)}{H^{\str}(\Q_{n,p}, A)_\mathrm{div}}$  is uniformly bounded if and only if the local strict Selmer condition at $p$ has uniformly bounded Cartesian defect. Then it follows from (a) that $\frac{H^{\str}(\Q_{n,p},A)}{H^{\str}(\Q_{n,p},A)_\mathrm{div}}$
 is bounded independent of $n$. Also recall the fact that $H^{\str}(\Q_{n,p},A)_{\mathrm{div}}=H^{\BK}(\Q_{n,p},A)$ for all $n$ (see \cite[Prop. 4.2.2(3)]{Ne3}) and hence 
$H^{\str}(\Q_{\infty,p},A)/H^\BK(\Q_{\infty,p},A)$ is finite. Since $\Q_\infty/\Q$ is a deeply ramified extension, $G_{\Q_{\infty,p}}$ has $p$-cohomological dimension one and   therefore  $H^1(\Q_{\infty,p},A^+)$ is divisible.  Consequently, $H^{\str}(\Q_{\infty,p},A)$, being the image of  $H^1(\Q_{\infty,p},A^+)$, is also divisible. On the other hand, by definition, $H^\BK(\Q_{\infty,p},A)$ is divisible. Now the theorem follows from the finiteness of $H^{\str}(\Q_{\infty,p},A)/H^\BK(\Q_{\infty,p},A)$.
\end{proof}
\begin{remark}
Although for simplicity, Theorem \ref{str_bk_at_p} is stated over $\Q_{\infty,p}$, the equality of the Bloch-Kato and strict Selmer local condition at a prime over $p$ can be extended to any deeply ramified extension of $\Q_p$ which contains  $\Q_{ \infty,p}.$
\end{remark}
\section{Comparison of root numbers} \label{cmparison_of_root_number}
In this section, we  consider two $p$-adic Galois representations  with isomorphic residual representations and   compare their local root numbers via certain Iwasawa invariants.

 Let $V_1$ and $V_2$ be two finite dimensional symplectic self-dual $p$-adic representation of $G_{\Q_\ell}$ (including $\ell=p$). Let $\FF$, $\OO$, $\pi$, $\mathfrak{f}$ be as before (see \S\ref{preliminaries}) such that $\FF$ contains the coefficient fields of $V_1$ and $V_2$. For $m\in \{1,2\}$, let $T_m$ be a $G_{\Q_\ell}$-stable symplectic self-dual $\OO$-lattice (see Remark \ref{self-dual_lattice}) of $V_m$ and put $A_m:=V_m/T_m$. Notice that we can identify $\bar{T}_m$ with $A_m[\pi]$. We assume that $A_1[\pi]\cong A_2[\pi]$ as $G_{\Q_\ell}$-modules.  
\subsection{Comparison of root numbers at a finite prime $\ell\neq p$} Let $\ell\neq p$ be a finite prime of $\Q$. 
From \cite[Prop. 2.12]{Ne}, there exists a non-degenerate symmetric bilinear pairing 
 \begin{equation}
 \label{pairing_on_H1}
     \langle \ , \ \rangle : H^1(I_\ell , T_m)_{\tor} \times H^1(I_\ell , T_m )_{\tor} \lra \mathcal{F} / \OO
 \end{equation}
for which $\Frob_\ell$ is an isometry. Further, by \cite[Corollary 1.5]{Ne}, we have 
\begin{equation}
    (-1)^{\dim_\mathfrak{f}\ \mathrm{Im}\left( H^1(I_\ell, T_m)_{\tor}^{\Frob_\ell =1} \lra H^1(I_\ell, T_m)_{\tor}/ \pi H^1(I_\ell, T_m)_{\tor} \right)}= \det (-\Frob_\ell | \ H^1(I_\ell, T_m)[\pi]).
\end{equation}
Now from the proof of \cite[Prop. 1.4(i)]{Ne}, we obtain that 
\begin{equation}
    \dim_\mathfrak{f}\ \mathrm{Im}\left( H^1(I_\ell, T_m)_{\tor}^{\Frob_\ell =1} \lra H^1(I_\ell, T_m)_{\tor}/ \pi H^1(I_\ell, T_m)_{\tor} \right) \equiv  \dim_\mathfrak{f}\ (H^1(I_\ell, T_m)_{\tor})_1[\pi] \pmod 2,
\end{equation}
where $(H^1(I_\ell, T_m)_{\tor})_1 : = \cup_{n\geq 1} \Ker\ (\Frob_\ell -1)^n$. Combining these two, we get that 
\begin{equation}
\label{det_and_generalised_1-space}
    (-1)^{\dim_\mathfrak{f}\ (H^1(I_\ell, T_m)_{\tor})_1[\pi]}= \det (-\Frob_\ell | \ H^1(I_\ell, T_m)[\pi]).
\end{equation}


From the exact sequence of $G_{\Q_\ell}$-modules
\begin{equation}
\label{exact_sequence}
    0 \lra T_m^{I_\ell}/{\pi T_m ^{I_\ell}} \lra (\bar T_m )^{I_\ell} \lra H^1(I_\ell , T_m )[\pi] \lra 0,
\end{equation}
we have that $\det (\Frob_\ell | \bar T_m ^{I_\ell})= \det (\Frob_\ell | \frac{T_m ^{I_\ell}}{\pi T_m ^{I_\ell}})\det (\Frob_\ell |  H^1 (I_\ell , T_m )[\pi])$.
The isomorphism $\bar T_1 \cong  \bar T_2$ implies that $\det (\Frob_\ell | \bar{T}_1^{I_\ell})=\det (\Frob_\ell | \bar{T}_2^{I_\ell})$. Hence, the ratio of the local root numbers satisfies the following: 
\begin{equation}
\label{root_number_comparison_equation}
\begin{aligned}
    \frac{W_\ell(V_1)}{W_\ell(V_2)} &= \frac{\det (-\Frob_\ell | V_2^{I_\ell})}{\det (-\Frob_\ell | V_1^{I_\ell})} \pmod{\pi} \\
&= \frac{\det (-\Frob_\ell |  \frac{T_2 ^{I_\ell}}{\pi T_2 ^{I_\ell}})}{\det (-\Frob_\ell |  \frac{T_1 ^{I_\ell}}{\pi T_1 ^{I_\ell}})}  \\
&\overset{\eqref{exact_sequence}}{=} \frac{ \det (-\Frob_\ell | H^1(I_\ell , T_1 )[\pi] )}{ \det (-\Frob_\ell |  H^1(I_\ell , T_2 )[\pi])} \\
&\overset{\eqref{det_and_generalised_1-space}}{=} (-1)^{\dim_\mathfrak{f}\ (H^1(I_\ell, T_1)_{\tor})_1[\pi]-\dim_\mathfrak{f}\ (H^1(I_\ell, T_2)_{\tor})_1[\pi] }.
\end{aligned}
 \end{equation}


 Fix a finite rational prime  $\ell \neq p$. For a prime $\eta$ of $\Q_\infty$ dividing $\ell$, let $\Q_{\infty,\eta}$ be a completion of $\Q_\infty$ at $\eta$.  Define 
\begin{equation}\label{delta_at_l}
\delta(A_m,\ell):=\rk_\mathcal{O}\    \Big( \prod_{\eta\mid \ell} H^1(\Q_{\infty,\eta} , A_m) ^\vee \Big) \quad \text{and} \quad \tau(A_m,\eta):= \rk_\mathcal{O}\  H^1(\Q_{\infty,\eta}, A_m)^\vee.
\end{equation}
Note that $\tau(A_m,\eta_1)=\tau(A_m,\eta_2)$ for any two primes $\eta_1,\eta_2$ of $\Q_\infty$ dividing $\ell$ (cf. \cite[\S5]{Greenberg_projectivemodule_modularrepresentation}). In addition, the number of primes of $\Q_\infty$ above $\ell$ is odd, so  
$\delta(A_m,\ell) \equiv \tau(A_m,\eta) \pmod2$ for any $\eta|\ell$. Recall that $\mathcal{O}/\pi=\mathfrak{f}$. We will use the following result of Greenberg-Vatsal:

\begin{proposition}\cite[Prop. 2.4]{Greenberg_Vatsal}
\label{multiplicity}
    Put $P_\ell^m(X)= \det\left(1-\mathrm{Frob}_\ell  X|_{(V_m)_{I_\ell}}\right) \in \mathcal{O}[X]$. Let \; $\widetilde{}$\; denote the reduction modulo $\pi$ map. Then $\tau(A_m,\eta)$ is equal to the multiplicity of $X=\widetilde{l}^{-1}$ as a root of $\widetilde{P_\ell^m}(X) \in \mathfrak{f} [X]$. \qed
\end{proposition}

Now observe that
\begin{align*}
    \det\left(1- X \Frob_\ell | (V_m)_{I_\ell}\right)&=    \det\left(1- X \Frob_\ell |( V_m^*(1))_{I_\ell}\right)\\
    &= \det(1- X \ell  \Frob_\ell | (V_m^{I_\ell})^*)\\
    &= \det(1- X\ell \Frob_\ell^{-1} | V_m^{I_\ell})\\
    &=\det( \Frob_\ell^{-1} | V_m^{I_\ell}) \det(\Frob_\ell - X\ell  | V_m^{I_\ell}).
\end{align*}

 Therefore the  multiplicity of $\widetilde \ell ^{-1}$ in $ \det(1- X \Frob_\ell | (V_m)_{I_\ell}) \pmod{\pi}$ is the same as the multiplicity of $\bar 1$ in $ \det(\Frob_\ell - X  | V_m^{I_\ell}) \pmod{\pi}$ and that is the same as the multiplicity of $\widetilde 1$ in $ \det(\Frob_\ell - X  |\frac{ T_m^{I_\ell}}{\pi T_m
^{I_\ell}} ) \pmod{\pi}$. 
Combined with Proposition \ref{multiplicity} and the above arguments, we obtain
\begin{equation}
\label{lambda_invariant_comaprison_equation}
\begin{aligned}
    {\delta(A_1,\ell) -\delta(A_2,\ell)} &\equiv  \text{multiplicity of $\widetilde 1$ in $ \det(\Frob_\ell - X  |\frac{ T_1^{I_\ell}}{\pi T_1^{I_\ell}} ) \pmod{\pi}$}\\
    & \quad \quad - \text{multiplicity of $\widetilde 1$ in $ \det(\Frob_\ell - X  |\frac{ T_2^{I_\ell}}{\pi T_2^{I_\ell}} ) \pmod{\pi}$}\pmod2 \\
& \overset{\eqref{exact_sequence}}{\equiv } \text{multiplicity of $\widetilde 1$ in $ \det(\Frob_\ell - X  |H^1(I_\ell, T_1)[\pi] ) \pmod{\pi}$} \\ & \quad \quad - \text{multiplicity of $\widetilde 1$ in $ \det(\Frob_\ell - X  |H^1(I_\ell, T_2)[\pi] ) \pmod{\pi}$} \pmod2 \\
&\equiv  {\dim_{\mathfrak{f}}\ (H^1(I_\ell , T_1)[\pi])_1 - \dim_{\mathfrak{f}}\ (H^1(I_\ell , T_2)[\pi])_1} \pmod2 .
\end{aligned}
 \end{equation}  

Now we are ready to relate local root numbers at $\ell \neq p$ with the local Iwasawa invariants defined above.             
\begin{proposition}\label{root_lambda_for_lneqp}
   Let $\ell\neq p$ be an integer prime and let $V_1$ and $V_2$ be two finite dimensional $\FF$-vector spaces endowed with a continuous $G_{\Q_\ell}$ action. Also, assume that they are symplectic self-dual and satisfy $A_1[\pi] \cong A_2[\pi]$ as $G_{\Q_\ell}$-modules. Then,   
\begin{equation*}
     \frac{W_\ell(V_1)}{W_\ell(V_2)} = (-1)^{\delta(A_1,\ell) -\delta(A_2,\ell)}. \qed
\end{equation*} 
\end{proposition}
\n \textbf{Proof:} Using equations \eqref{root_number_comparison_equation} and \eqref{lambda_invariant_comaprison_equation}, it suffices to show  that for $m\in \{ 1,2\}$, 
\begin{equation}\label{claim_on_1_space_in_out}
     (H^1(I_\ell , T_m)[\pi])_1 =  (H^1(I_\ell, T_m)_{\tor})_1[\pi].
\end{equation}
We establish this in the rest of the proof. 
  Clearly, $(H^1(I_\ell, T_m)_{\tor})_1[\pi] \subseteq (H^1(I_\ell , T_m)[\pi])_1$. For the other inclusion, let $x\in (H^1(I_\ell , T_m)[\pi])_1\subset H^1(I_\ell , T_m)_{\tor}$.
  There is an orthogonal decomposition $H^1(I_\ell, T_m)_{\tor} = (H^1(I_\ell, T_m)_{\tor})_1 \oplus (H^1(I_\ell, T_m)_{\tor})_1^\perp$, where $(H^1(I_\ell, T_m)_{\tor})_1^\perp = \cap_{n\geq 1} \mathrm{Im}\ (\Frob_\ell -1)^n$. Moreover,
  $\Frob_\ell -1$ is nilpotent (respectively invertible) on $(H^1(I_\ell, T_m)_{\tor})_1$ (respectively on $(H^1(I_\ell, T_m)_{\tor})_1^\perp$) (see \cite[Prop. 1.4(i)]{Ne}).
  Then $x$ can be written uniquely as $x=x_1 + x_2$, where $x_1\in (H^1(I_\ell, T_m)_{\tor})_1$ and $x_2 \in (H^1(I_\ell, T_m)_{\tor})_1^\perp$.    From the definition of $(\cdot )_1$, there exists $r_1, r_2 \in \N$ such that $(\Frob_\ell -1)^{r_1}(x)=0$ and $(\Frob_\ell -1)^{r_2}(x_1)=0$. So for $r= \max\{ r_1, r_2\}$, $(\Frob_\ell -1)^{r}(x_2)=0$ and hence $x_2\in (H^1(I_\ell, T_m)_{\tor})_1$ which concludes that $x_2=0$. So $x=x_1\in (H^1(I_\ell, T_m)_{\tor})_1[\pi]$ and  $(H^1(I_\ell , T_m)[\pi])_1  \subseteq (H^1(I_\ell, T_m)_{\tor})_1[\pi]$.  \qed 

 \subsection{Comparison of local root numbers at $p$} Next, we compare the local root numbers of two residually isomorphic  $p$-adic Galois representations  of $G_{\Q_p}$. This case is more delicate and we  make some  additional hypotheses on the $G_{\Q_p}$ representations, as discussed below. 

Let $V_1$ and $V_2$ be two finite dimensional $\FF$-vector spaces endowed with a continuous $G_{\Q_p}$ action such that $V_1$ and $V_2$ are symplectic self-dual representations and $A_1[\pi] \cong A_2[\pi]$ as $G_{\Q_p}$-modules. We assume that $V_1$ and $V_2$ are  de Rham representation, pure of weight $-1$ and  satisfy the Panchishkin condition. 
 For $m\in\{1,2\}$, let $A_m^+$ be the image of $V_m^+$ (defined in \eqref{panchishkin}) under the natural map $V_m\lra A_m$ and set $A_m^-:= A_m/A_m^+$. We further assume that the isomorphism $A_1[\pi] \cong A_2[\pi]$ induces   isomorphism $A_1^-[\pi] \cong A_2^-[\pi]$ of $G_{\Q_p}$ modules. 
   As $V_m$ is a symplectic self-dual representation, the epsilon factor $\varepsilon_p(\WD (V_m))$ (see \S \ref{Weil-Deligne representation and local root numbers}) takes the value $\pm 1$. 
Further, by (cf. \cite[Prop. 3.3.3(4)]{Ne3})
\begin{equation}
\label{epsilon_at_p}
    \varepsilon_p(\WD (V_m)) = (-1)^{\dim_\FF\  H^0(\Q_p, V_m^-)} (-1)^{d^-(V_m)} (\mathrm{det}\  V_m^+)(-1).
\end{equation}
Define the completed epsilon constant $\hat{\varepsilon}_p(\WD (V_m))$ by (cf. \cite[Definition 2.10, Lemma 2.8]{bkno})
\begin{equation}\label{completed_epsilon}
    \hat{\varepsilon}_p(\WD (V_m)):= \varepsilon_p(\WD (V_m)) (-1)^{\sum_{i>0} i d^i(V_m)}.
\end{equation}
Recall from \eqref{definition_of_d} that $(-1)^{d(V_m)}= (-1)^{d^-(V_m)+  \sum_{i>0} i d^i(V_m)}$. 
From \cite[Prop. 12.1.6(ii)]{Be},
 we have $d (V_m)= d(\det (V_m))$ and consequently, we obtain 
 \begin{equation*}
     \hat{\varepsilon}_p(\WD (V_m))= (-1)^{\dim_\FF\  H^0(\Q_p, V_m^-)} (-1)^{d(\det (V_m))} (\mathrm{det}\  V_m^+)(-1).
 \end{equation*}
Since $V_m$ is symplectic self-dual, from \cite[Prop. 7.17]{FO}, we get $\det(V_m)= \chi_p^{\frac{\dim_\FF V_m}{2}}$. 
Now the isomorphism of residual representations $A_1[\pi] \cong A_2[\pi]$ forces $\dim_\FF V_1= \dim_\FF V_2$ and hence, we obtain ${d(\det (V_1))}={d(\det (V_2))}$.

Next, using the isomorphisms of residual representations $A_1[\pi] \cong A_2[\pi]$ and $A_1^-[\pi] \cong A_2^-[\pi]$, we obtain $\det(V_1)\equiv \det(V_2)\pmod{\pi}$ and $\det(V_1^-)\equiv \det(V_2^-)\pmod{\pi}$. Since determinant splits over exact sequence, i.e., we have $\det(V_m)= \det(V_m^+)\det(V_m^-)$, and hence $\det(V_1^+)\equiv  \det(V_2^+)\pmod{\pi}$. Thus we obtain $(\det(V_1^+))(-1)\equiv (\det(V_2^+))(-1)\pmod{\pi}$ and consequently $(\det(V_1^+))(-1)= (\det(V_2^+))(-1)$ as $p$ is odd.


 Combining everything, we have the following result:
\begin{proposition}\label{root_lambda_at_p}
    Let $V_1$ and $V_2$ be two finite dimensional symplectic self-dual $p$-adic Galois  representations of $G_{\Q_p}$ such that $A_1[\pi] \cong A_2[\pi]$ as $G_{\Q_p}$-modules.  We make the following additional assumptions.
\begin{enumerate}
    \item $V_1$ and $V_2$ are de Rham representations, pure (of weight $-1$) and satisfy the Panchishkin condition.
    \item The isomorphism $A_1[\pi] \cong A_2[\pi]$  in turn induces an isomorphism $A_1^-[\pi] \cong A_2^-[\pi]$ of $G_{\Q_p}$-modules.
\end{enumerate} 
 \begin{equation*}
 \text{Then} \quad    \frac{\hat{\varepsilon}_p(\WD (V_1))}{\hat{\varepsilon}_p(\WD (V_2))}= (-1)^{\dim_\FF\  H^0(\Q_p, V_1^-)-\dim_\FF\  H^0(\Q_p, V_2^-)}. \tag*{\qed}
\end{equation*}
\end{proposition}

 Now we define the local Iwasawa invariants at $p$ as follows. 
 Define
\begin{equation}
    H^\str (\Q_{\infty,p}, A_m[\pi]):= \Ker\left( H^1(\Q_{\infty,{p}}, A_m[\pi]) \lra H^1(\Q_{\infty,{p}}, A_m^-) \right).
\end{equation}
We keep the hypothesis that the isomorphism $A_1[\pi] \cong A_2[\pi]$  of $G_{\Q_p}$-modules further induces a $G_{\Q_p}$ isomorphism $A_1^-[\pi] \cong A_2^-[\pi]$. Then it follows that  
 $\mathcal{L}_m:= \frac{ H^\str (\Q_{\infty,p}, A_m[\pi])}{ H^\str (\Q_{\infty,p}, A_1[\pi]) \cap  H^\str (\Q_{\infty,p}, A_2[\pi])}$ is finite (see the proof of Proposition \ref{str_root_at_p})  and put 
\begin{small}\begin{equation}\label{delta_at_p}
    \delta(A_m,p):= \dim_\mathfrak{f}\  \mathcal{L}_m = \dim_\mathfrak{f}\ \frac{ H^\str (\Q_{\infty,p}, A_m[\pi])}{ H^\str (\Q_{\infty,p}, A_1[\pi]) \cap  H^\str (\Q_{\infty,p}, A_2[\pi])}, 
\end{equation}\end{small}
where the intersection is taken via the identification $H^1(\Q_{\infty,{p}}, A_1[\pi])=H^1(\Q_{\infty,{p}}, A_2[\pi])$.

Now we are ready  to relate the local root number at $p$ with the local Iwasawa invariants at $p$ defined as above.
\begin{proposition}\label{str_root_at_p}
   We keep the setting and  the hypotheses of Proposition \ref{root_lambda_at_p}. 
   Then 
\begin{small}\begin{equation*}
\frac{\hat{\varepsilon}_p(\WD (V_1))}{\hat{\varepsilon}_p(\WD (V_2))}= (-1)^{\delta(A_1,\ p) -\delta(A_2,\ p)}.
\end{equation*}\end{small}
\end{proposition}
\begin{proof}
Define $X_m:= \Ker(H^1(\Q_{\infty,p}, A_m[\pi]) \lra H^1(\Q_{\infty,p}, A_m^-[\pi]) )$ and
consider the commutative diagram 
\begin{equation}\label{cd_proposition_2.7}
\begin{tikzcd}
	0 & X_m & {H^1(\Q_{\infty,p}, A_m[\pi])} & {H^1(\Q_{\infty,p}, A_m^-[\pi])} & 0 \\
	0 & {H^\str (\Q_{\infty,p}, A_m[\pi])} & {H^1(\Q_{\infty,p}, A_m[\pi])} & {H^1(\Q_{\infty,p}, A_m^-)[\pi],}
	\arrow[from=1-1, to=1-2]
	\arrow[from=1-2, to=1-3]
	\arrow["\varphi_m", from=1-2, to=2-2]
	\arrow[from=1-3, to=1-4]
	\arrow["{= }", from=1-3, to=2-3]
	\arrow[from=1-4, to=1-5]
	\arrow["\psi_m"', from=1-4, to=2-4]
	\arrow[from=2-1, to=2-2]
	\arrow[from=2-2, to=2-3]
	\arrow[from=2-3, to=2-4]
\end{tikzcd}
    \end{equation}
 Using assumption $(2)$, we can identify $X_1$ with $X_2$ and denote it by $X$.
 Further using the injectivity of $\varphi_m$, we can view $X$ as a submodule of $  H^\str (\Q_{\infty,p}, A_1[\pi]) \cap  H^\str (\Q_{\infty,p}, A_2[\pi])$. By applying a Snake lemma in \eqref{cd_proposition_2.7}, we obtain that $\Coker (\varphi_m) \cong \Ker(\psi_m) =  H^0(\Q_{\infty,p}, A_m^-)/ \pi$, which is finite. This in turn  implies the finiteness of $\mathcal{L}_m$ in \eqref{delta_at_p}. Now
\begin{equation}\label{delta_rank_at_p}
\begin{aligned}
\delta(A_1,\ p) -\delta(A_2,\ p) = &\dim_\mathfrak{f}\   \frac{ H^\str (\Q_{\infty,p}, A_1[\pi])}{X}   - \dim_\mathfrak{f}\   \frac{ H^\str (\Q_{\infty,p}, A_1[\pi]) \cap  H^\str (\Q_{\infty,p}, A_2[\pi])}{X} \\ - &
\dim_\mathfrak{f}\   \frac{ H^\str (\Q_{\infty,p}, A_2[\pi])}{X}  + \dim_\mathfrak{f}\   \frac{H^\str (\Q_{\infty,p}, A_1[\pi]) \cap  H^\str (\Q_{\infty,p}, A_2[\pi])}{X}\\
=\ & \rank_{\mathfrak{f}}\ \Coker(\varphi_1)-  \rank_{\mathfrak{f}}\ \Coker(\varphi_2) \\
=\ &\rank_{\mathfrak{f}}\ H^0(\Q_{\infty,p}, A_1^-)/ \pi - \rank_{\mathfrak{f}}\ H^0(\Q_{\infty,p}, A_2^-)/ \pi\\
=\ & \rank_\OO\ H^0(\Q_{\infty,p}, A_2^-)^\vee - \rank_\OO\ H^0(\Q_{\infty,p}, A_1^-)^\vee.
\end{aligned}
\end{equation}
Further, $\rank_\OO\ H^0(\Q_{\infty,p}, A_m^-)^\vee = \dim_\FF\  H^0(\Q_{\infty,p}, V_m^-)$ and we obtain 
\begin{equation}
    \delta(A_1,\ p) -\delta(A_2,\ p) = \dim_\FF\  H^0(\Q_{\infty,p}, V_2^-)-  \dim_\FF\  H^0(\Q_{\infty,p}, V_1^-).
\end{equation}
Using Proposition \ref{root_lambda_at_p}, to complete the proof of this proposition, it suffices to show  
\begin{equation}\label{Q_cyc_Q_invariant}
    \dim_\FF\  H^0(\Q_p, V_m^-)\equiv \dim_\FF\  H^0(\Q_{\infty,p}, V_m^-) \pmod 2,
\end{equation}
which is proved in Corollary \ref{mod_2_Qcyc_Q} and the proposition follows.
\end{proof}

 \begin{proposition}\label{prop_epsilon_delta_p_without_hat}
     We keep the setting and  the hypotheses of Proposition \ref{str_root_at_p}. In addition we assume that both $V_1$ and $V_2$ are semistable  $G_{\Q_p}$-representations. Then 
     \begin{equation}\label{epsilon_delta_p_without_hat}
         \frac{{\varepsilon}_p(\WD (V_1))}{{\varepsilon}_p(\WD (V_2))}= (-1)^{\delta(A_1,\ p) -\delta(A_2,\ p)}.
     \end{equation}
 \end{proposition}
\begin{proof}
  In the proof of Proposition \ref{root_lambda_at_p}, we have shown that $\mathrm{det}\  (V_1^+)(-1)=\mathrm{det}\  (V_2^+)(-1)$. Further, it follows from  the proof of Proposition \ref{str_root_at_p}  that 
  \begin{equation*}
      \delta(A_1,p)- \delta(A_2,p) \equiv \dim_\FF\  H^0(\Q_p, V_1^-)-\dim_\FF\  H^0(\Q_p, V_2^-) \pmod2.
  \end{equation*}
 Then using equation \eqref{epsilon_at_p}, it suffices to establish that $(-1)^{d^-(V_1)}=(-1)^{d^-(V_2)}$. From \cite[Prop. 3.3.3(3)]{Ne3},  $(-1)^{d^-(V_m)}=(-1)^{d (V_m^+)} $  and moreover,  by \cite[Prop. 12.1.6(ii)]{Be}, we have $d (V_m^+)= d(\det (V_m^+))$.  Now
   $\det V_m = \bigwedge ^{\dim V_m} V_m$   and $\det V_m^-$ are (one dimensional) semistable representations of $G_{\Q_p}$ and  we have $\det (V_m) =  \chi_p ^{\frac{\dim_\FF V_m}{2}}$  and $\det(V_m^-)= \phi_m \chi_p^{t_m} $,
for some unramified character  $\phi_m$ and for $t_m\in \Z$ (see  \cite[Prop. 7.17]{FO}). 
Then using the isomorphisms of residual representations $A_1[\pi]\cong A_2[\pi]$ and  $A_1^-[\pi]\cong A_2^-[\pi]$, we deduce  $\dim_\FF V_1= \dim_\FF V_2$,  $t_1\equiv t_2 \pmod{p-1}$ and $\phi_1 \equiv \phi_2 \pmod{\pi}$. Consequently, $\det(V_m^+)= \phi_m^{-1} \chi_p^{s_m}$, where $s_m= \frac{\dim_\FF V_m}{2}- t_m$. As twisting by unramified character does not change the Hodge-Tate weights,
we conclude that $d(\det(V_1^+))= d(\chi_p^{s_1})= -s_1 \equiv -s_2 \pmod 2$, as required. 
\end{proof}

\begin{remark}\label{semistable_str-root_at_p_different_HT_weights}
In particular,  our result in Proposition \ref{prop_epsilon_delta_p_without_hat}  covers the cases where the semistable representations $V_1$ and $V_2$ have different Hodge-Tate weights. 

On the other hand, suppose that $V_1$ and $V_2$ are de Rham (but not necessarily semistable) representations of $G_{\Q_p}$ with the same Hodge-Tate weights. In this case,  $d(V_1^+)= d(V_2^+)$ and from \cite[Prop. 3.3.3(3)]{Ne3}, we get $(-1)^{d^-(V_m)}= (-1)^{d (V_m^+)}$ and equation \eqref{epsilon_delta_p_without_hat} continues to hold.
\end{remark}

\subsection{Comparison of local root numbers at an Archimedean prime}
Let $V_1$ and $V_2$ be two finite dimensional $\FF$-vector spaces endowed with a continuous $G_{\Q}$ action satisfying \hyperref[Hyp:Mov]{(Mov)}.  Now we compare the local root numbers at the Archimedean prime.

\begin{proposition}\label{root_lambda_at_infty}
    Let $V_1$ and $V_2$ be two symplectic self-dual, $p$-adic Galois representations of $G_{\Q}$ satisfying \hyperref[Hyp:Mov]{(Mov)} such that $A_1[\pi] \cong A_2[\pi]$  as $G_{\Q_p}$-modules. Assume $V_1$ and $V_2$ are pure (of weight $-1$) and satisfy the Panchishkin condition at $p$.   We also assume that
the isomorphism $A_1[\pi] \cong A_2[\pi]$ in turn induces an isomorphism $A_1^-[\pi] \cong A_2^-[\pi]$ of $G_{\Q_p}$-modules. {Then }  
 \begin{equation*}
\frac{W_\infty(V_1)}{W_\infty(V_2)}= (-1)^{d^-(V_1)-d^-(V_2)}.
\end{equation*}
\end{proposition}
\n\textbf{Proof:} Let $\infty$ denote the Archimedean place of $\Q$.  Then in this case, we have the formula, $W_\infty (V_m)= (-1)^{d^-(V_m)}$ (cf. \cite[\S 4.1.3]{Ne3}). 
 \qed

    From the definition of the $\hat{\varepsilon}(\WD(V_m))$ (equation \eqref{completed_epsilon}) and  Proposition \ref{root_lambda_at_infty}, we observe that 
    \begin{equation*}
\frac{W_p(V_1)W_\infty(V_1)}{W_p(V_2)W_\infty(V_2)}= \frac{\hat{\varepsilon}_p(\WD(V_1))(-1)^{d^-(V_1)+\sum_{i>0} i d^i(V_1)}}{\hat{\varepsilon}_p(\WD(V_2))(-1)^{d^-(V_2)+\sum_{i>0} i d^i(V_2)}}.
    \end{equation*}
On the other hand,  $(-1)^{d^-(V_m)+\sum_{i>0} i d^i(V_m)}=(-1)^{d(V_m)}$ (see Definition \eqref{definition_of_d}). It follows from the proof of Proposition \ref{root_lambda_at_p} that $d(V_1)=d(\det\ V_1)= d(\det\ V_2)=d(V_2)$. 
Hence we arrive at the following result:
\begin{lemma}\label{root_p_infty_epsilonhat}
We keep the setting and hypotheses of Proposition \ref{root_lambda_at_infty}. Then 
\begin{equation*}
\frac{W_p(V_1)W_\infty(V_1)}{W_p(V_2)W_\infty(V_2)}= \frac{\hat{\varepsilon}_p(\WD(V_1))}{\hat{\varepsilon}_p(\WD(V_2))}. \tag*{\qed}
\end{equation*}
\end{lemma}

\subsection{Comparison of global root number}
Combining the  results on the comparison of the local root numbers of two congruent $p$-adic representations, we arrive at the global root number comparison in terms of our local invariants $\delta(-,\ell)$ stated in  Theorem \ref{global_root_lambda}.

\textbf{Proof of Theorem \ref{global_root_lambda}:}
From the results obtained for local root numbers at $\ell\neq p$ (Proposition \ref{root_lambda_for_lneqp}), at $\ell=p$ (Proposition \ref{str_root_at_p}) and for the Archimedean prime (Proposition \ref{root_lambda_at_infty} and Lemma \ref{root_p_infty_epsilonhat}), it is immediate that 
\begin{equation*}
   \frac{W(V_1)}{W(V_2)}=(-1)^{\sum_{{\ell } }\  (\delta(A_1,\ell)- \delta(A_2,\ell))}.  
\end{equation*}
Further, if  $V_1$ and $V_2$ both are unramified at a finite prime $\ell\neq p$, then  by Lemma \ref{root_number_V_1_V_2}, we have $W_\ell(V_1)=W_\ell(V_2)$ and  $\delta(A_1,\ell)=\delta(A_2,\ell)$ by Proposition \ref{multiplicity}. So the summation in Theorem \ref{global_root_lambda} can be taken over finitely many finite primes; more specifically the primes $\ell\ (\neq p)$ where either $V_1$ or $V_2$ is ramified and the prime $p$. \qed


\subsection{Arithmetic local constants} Mazur-Rubin developed the theory of arithmetic local constants in \cite{MR}, which is used by Nekov\'{a}\v{r} to prove many cases of $p$-parity conjecture \cite{Ne, Ne4}. Nekov\'{a}\v{r} related the  ratio of local root numbers with arithmetic  local constants. On the other hand, we have expressed the ratio of local root numbers in terms of certain Iwasawa invariants at $\ell\neq p$ (see \eqref{delta_at_l}) and also at $\ell=p$ (see \eqref{delta_at_p}). In Proposition \ref{lnp}, we  relate   the arithmetic local constants with the same Iwasawa invariants, which shows the compatibility of our results with the existing literature.

Let $V_1$ and $V_2$ be two finite dimensional symplectic self-dual  $\FF$-vector spaces equipped with a continuous $G_\Q$ action satisfying \hyperref[Hyp:Mov]{(Mov)}. We recall   the  definition of the local (vector space) Bloch-Kato Selmer condition. 

Let $v$ be a finite prime of a number field $L$ and $I_v$ be the inertia subgroup of $G_{L_v}$. Define 
\begin{equation}\label{local_BK_selmer}
    H^\BK(L_{v},V_i):=\left\{ \begin{array}{lc}
       \Ker\left(H^1(L_v,V_i)\rightarrow H^1(I_v,V_i)\right)  & \text{if}\;\; v \nmid p  \\
       \Ker\left(H^1(L_v,V_i)\rightarrow H^1(L_v,V_i \otimes_{\Q_p}B_{crys} )\right) & \text{if}\;\; v \mid p
    \end{array}\right.
\end{equation}

Consider any prime $\ell$ (including $p$) of $\Q$. Define 
\begin{align*}
& H_f^1(\Q_\ell, T_i):=\Ker\left(H^1(\Q_\ell,T_i)\lra \frac{H^1(\Q_\ell,V_i)}{H^\BK(\Q_\ell,V_i)}\right) \quad \text{and}\\
&\mathfrak{F}_\ell^i:= \mathrm{Image}(H_f^1(\Q_\ell, T_i) \lra H^1(\Q_\ell,A_i[\pi]).
\end{align*}
 In addition, also assume that $A_1[\pi]\cong A_2[\pi]$ as symplectic $G_\Q$-modules, i.e. the Galois isomorphism preserves the symplectic pairings on $A_1[\pi]$ and $A_2[\pi]$. The subspace $\mathfrak{F}_\ell^i$ is a finite dimensional $\mathfrak{f}$-vector space and is a Lagrangian subspace (cf. \cite[\S4]{MR}, \cite{Ne}) in $H^1(\Q_\ell, A_1[\pi])\cong H^1(\Q_\ell, A_2[\pi])$. For any two Lagrangian subspaces $\mathfrak{G}$, $\mathfrak{H}$ of $H^1(\Q_\ell, A_1[\pi])\cong H^1(\Q_\ell, A_2[\pi])$, define the arithmetic local constant $d(\mathfrak{G}, \mathfrak{H})$ by 
\begin{equation}\label{arithmetic_local-constant}
    d(\mathfrak{G}, \mathfrak{H}):= \dim_{\mathfrak{f}}\left( \frac{\mathfrak{G}}{\mathfrak{G}\cap \mathfrak{H}} \right) \pmod2 = \dim_{\mathfrak{f}}\left( \frac{\mathfrak{H}}{\mathfrak{G}\cap \mathfrak{H}} \right) \pmod2 \in \Z/2\Z.
\end{equation}

Let $\Sigma$ be a finite set of primes of $\Q$ containing the prime $p$, the primes where either $V_1$ or $V_2$ is ramified and the infinite prime, and for $m =1,2$, let $S_\BK(A_m/\Q)$ be  the Bloch-Kato Selmer group of $A_m$ over $\Q$ defined in equation \eqref{global_BK_selmer_group}. 
It is easy to see that if $\ell \notin \Sigma$, then $d(\mathfrak{F}_\ell^1, \mathfrak{F}_\ell^2)=0$.
Now under the assumption that $A_1[\pi]\cong A_2[\pi]$ as symplectic $G_\Q$-modules and $V$ is pure of weight $-1$, it follows from \cite[Theorem 1.4]{MR} and \cite[Eq. 4.6.2]{Ne} that 
\begin{equation}\label{MRtheorem}
    \rk_\OO\ S_\BK(A_1/\Q)^\vee -\rk_\OO\ S_\BK(A_2/\Q)^\vee \equiv \sum_{\ell\in \Sigma} d(\mathfrak{F}_\ell^1, \mathfrak{F}_\ell^2) \pmod2.
\end{equation}
 
Recall that two $\mathrm{GSp}_{2d}$ representations are said to be equivalent by a similitude if the associated bilinear pairing are equal up to multiplication by a non-zero scalar. 
Next we illustrate an example of a pair of Galois representations that are residually equivalent but not equivalent by symplectic similitude. 
\begin{example}\label{example_isomorphic_residual_but_not_symplectic_similitude} 
Let $A$, $B$ and $C$ be three elliptic curves defined over $\Q$. Suppose that $A[3]$, $B[3]$ and $C[3]$ are all irreducible $G_\Q$-module such that $A[3]$ and $C[3]$ are not isomorphic over $\Q$.  Further suppose that there is a degree $n$ isogeny $\phi : A \longrightarrow B$ with $n \equiv 2 \pmod3$.  As $A[3]$ is irreducible, $\phi : A[3] \longrightarrow B[3]$ is an isomorphism of Galois module. Suppose that endomorphism rings of all the three $G_\Q$-modules $A[3]$, $B[3]$ and $C[3]$ are isomorphic to $\Z/3\Z$. Then any morphism $g :A[3]\longrightarrow B[3]$ is given by a scalar multiple of $\phi$. Let $e_X :X[3] \times X[3]\longrightarrow \mu_3$ be the Weil-pairing for $X\in \{A, B, C\}$. Then $e_B(\phi(x),\phi(y))=e_A(x,y)^n=e_A(x,y)^2$ (see \cite[Proposition 8.2, Chapter III]{S}). In particular, $\phi$ is a symplectic similitude isomorphism but not a symplectic isomorphism from $A[3]$ to $B[3]$. Define $\bar{\rho}_1:=A[3]\oplus C[3]$ and  $\bar{\rho}_2:=B[3]\oplus C[3]$. 
We get a pair of representations
\[ \bar\rho_i : G_\Q \longrightarrow \mathrm{GSp}_4(\F_3)\]
for $i=1,2$. Also define 
\[  e_1 : \bar\rho_1 \times \bar\rho_1 \longrightarrow \mu_3,  \  \  \  \  e_2 : \bar\rho_2 \times \bar\rho_2 \longrightarrow \mu_3 \] 
by the formula $e_1((x,y),(z,w)):=e_A(x,z)e_C(y,w)$  and $e_2((x,y),(z,w)):=e_B(x,z)e_C(y,w)$. 

We note that even though the isomorphism $\phi : A[3] \longrightarrow B[3]$ is symplectic similitude, the induced isomorphism $\tilde{\phi} :  \bar\rho_1 \longrightarrow \bar\rho_2$ is not symplectic similitude. 
In fact, there does not exist an isomorphism $f : \bar{\rho}_1 \lra \bar{\rho}_2  $ which is symplectic similitude. To see this, suppose such an isomorphism $f$ exists. As $A[3]$ is not isomorphic to $C[3]$ any such isomorphism $f$ is of the form $(g,h)$ where $g : A[3] \lra B[3]$ and $h : C[3] \lra C[3]$ are isomorphisms. By assumption $g=a\phi$ and $h=b I_2$ for $a,b\in (\Z/3\Z)^{\times}$.  Since $f$ is assumed to be symplectic similitude there exist a $\lambda \in (\Z/3\Z)^{\times}$ such that $e_2(f(x,y),f(z,w))=e_1((x,y),(z,w))^\lambda$ for all $x,z\in A[3]$ and $y,w\in C[3]$. This gives $e_B(g(x),g(z))e_C(h(y),h(w))=e_A(x,z)^\lambda e_C(y,w)^\lambda $ for all $x,z\in A[3]$ and $y,w\in C[3]$. Therefore 
\[ e_B(g(x),g(z))=e_A(x,z)^\lambda   \  \  \text{and} \  \ e_C(h(y),h(w))= e_C(y,w)^\lambda \] 
for all $x,z\in A[3]$ and $y,w\in C[3]$. Substituting $g=a\phi$ and $h=b I_2$ we get 
\[ e_B(\phi(x),\phi(z))^{a^2}=e_A(x,z)^\lambda   \  \  \text{and} \  \ e_C(y,w)^{b^2}= e_C(y,w)^\lambda .\]
We have $a^2=1=b^2 \pmod3$. Now $e_C(y,w)^{b^2}= e_C(y,w)^\lambda$ this implies that $\lambda=1$. Therefore $e_B(\phi(x),\phi(z))^{a^2}=e_A(x,z)$. But this contradicts the fact that $\phi$ is not a symplectic isomorphism. Hence, no such $f$ exists. 

Notice that $p=3$ implies $\mu_3\cong \F_3(1)\cong \F_3(-1)$. Now from \cite[Theorem 10.2.1]{BCGP} there exists an infinite family of principally polarised abelian varieties $A_i$ over $\Q$ up to twist, such that $\mathrm{End}(A_i)=\Z$ and    $A_i[3] \cong \bar\rho_i$ symplectically for $i=1,2$. Moreover,  if $A$, $B$, $C$ have good and ordinary reduction at $3$, then the $A_i$ can be chosen to be an infinite family of good ordinary abelian varieties.   Note that given an element in this family of abelian variety there can be only finitely many that are isogenous to it and furthermore $\mathrm{End}(A_i)=\Z$ implies that there can be finitely many which are isomorphic to it up to a twist.

Any pair of non-isogenous abelian varieties $A_1$ and $A_2$ over $\Q$ constructed in this manner have isomorphic residual representation at $p=3$ but not by a  symplectic similitude. \qed   
\end{example}
The above example of pair of Galois representations are residually reducible. If the residual representations are absolutely irreducible, then in that case it can be shown that any two such residually equivalent Galois representation  are in fact residually equivalent by symplectic similitude (see \cite{Ne4}). More generally, when residual representation are not  necessarily absolutely irreducible the following is true :

\begin{proposition}\label{sym}
Let $\bar\F$ be a separably closed field of characteristic $p\neq 2$ and $\tau_i :G_\Q \lra \mathrm{GSp}_{2d}(\bar{\F})$  for $i=1,2$ be two Galois representation with a similitude character $\chi$. Suppose $\tau_1$ and $\tau_2$ are equivalent as $G_\Q$ modules. Then $\tau_1$ and $\tau_2$ are $\mathrm{GSp}_{2d}$-equivalent by similitude.    
\end{proposition}
\begin{proof}
This is a well-known algebraic fact (see, for instance, \cite[Theorem~2.1]{K}). Although \cite[Theorem~2.1]{K} treats complex semisimple representations, the statement remains true over any separably closed field of odd characteristic (see \cite[Lemma 2.1.3]{BCGP}). The key ingredient is that every invertible self-adjoint endomorphism of a symplectic space admits a self-adjoint square root, which may be chosen as a polynomial in the endomorphism.
\end{proof}

  Consider any prime $\ell\neq p$. Then the unramified cohomology $\mathfrak{F}_\ell^\ur := H^1(G_{\Q_\ell}/I_\ell,A_i[\pi]^{I_\ell})$ is a Lagrangian subspace of $H^1(\Q_\ell, A_i[\pi])$ (see \cite[\S 2.6]{Ne}). By a result of Klagsbrun-Mazur-Rubin (\cite[Prop. 2.9]{Ne}), we have 
\begin{equation}
    \label{KMR}d(\mathfrak{F}_\ell^1, \mathfrak{F}_\ell^2)= d(\mathfrak{F}_\ell^1, \mathfrak{F}_\ell^\ur)- d(\mathfrak{F}_\ell^\ur, \mathfrak{F}_\ell^2).
\end{equation}

\begin{proposition}\label{lnp} Let $\ell\neq p$ be a prime and
 $V_1$ and $V_2$ be two finite dimensional symplectic self-dual  $\FF$-vector spaces equipped with a continuous $G_{\Q_\ell}$ action. In addition, also assume that $A_1[\pi]\cong A_2[\pi]$ as symplectic $G_{\Q_\ell}$-modules. Then
     we have \begin{equation*}
    d(\mathfrak{F}_\ell^1, \mathfrak{F}_\ell^2) \equiv \delta(A_1,\ell) -\delta(A_2,\ell) \pmod2.
\end{equation*}
\end{proposition}
\n \textbf{Proof:} 
By \cite[Prop. 2.13]{Ne}, we have $(-1)^{d(\mathfrak{F}_\ell^i, \mathfrak{F}_\ell^\ur)}= \det(-\Frob_\ell \ |\ H^1(I_\ell, T_i)[\pi])$. Now using equation \eqref{det_and_generalised_1-space}, we get that 
\begin{equation*}
(-1)^{d(\mathfrak{F}_\ell^i, \mathfrak{F}_\ell^\ur)}= (-1)^{\dim_\mathfrak{f}\ (H^1(I_\ell, T_i)_{\tor})_1[\pi]}.
\end{equation*}
Further, it follows from equations \eqref{lambda_invariant_comaprison_equation}  and \eqref{claim_on_1_space_in_out} that 
\begin{equation*}
    \dim_\mathfrak{f}\ (H^1(I_\ell, T_1)_{\tor})_1[\pi] - \dim_\mathfrak{f}\ (H^1(I_\ell, T_2)_{\tor})_1[\pi]\equiv \delta(A_1,\ell) -\delta(A_2,\ell) \pmod2.
\end{equation*}
Now the result follows from \eqref{KMR}. \qed


\begin{remark}\label{remark 4.15}
Note that Nekov\'a\v r's comparison of the local root numbers for $A_1[\pi]$ and $A_2[\pi]$ requires a symplectic isomorphism between the modules.
 We stress that we do not assume the $G_{\Q_\ell}$ isomorphism (including $\ell=p$) between $A_1[\pi]$ and $A_2[\pi]$ preserves the symplectic pairing, i.e., in our setting $A_1[\pi]$ and $A_2[\pi]$ need not be symplectically isomorphic. 
 
Let us first consider the case where $\ell\neq p$. In view of Proposition \ref{sym}, we can always extend scalars to a finite extension of $\mathcal{F}$ so that the residual representation in Proposition  \ref{lnp} become symplectically isomorphic. 
Note that the root number and $\delta(-,\ell)$ does not change under the extension of scalars. Thus the relation between the arithmetic local constants and the local Iwasawa invariants for $\ell\neq p$, as stated in Proposition \ref{lnp}, would always hold. This does not seem to be pointed out in the literature.

   The case of $\ell = p$ is more subtle. The relation between the arithmetic local constant $d(\mathfrak{F}_p^1, \mathfrak{F}_p^2)$ and root numbers $W(-,p)$ are known only in special cases (\cite{Ne4}, \cite{Ne} \cite{bkno}), under the assumption that $A_1[\pi]$ and $A_2[\pi]$ are symplectically isomorphic.

  In this article, we have introduced the local invariants $\delta(-,p)$ for  general $p$-adic de Rham representation satisfying Panchishkin condition at $p$ (see \eqref{delta_at_p}). Further, we have expressed  the ratio of the local root numbers at $p$ in terms of $\delta(-,p)$ and  it is a key ingredient in establishing the relative $p$-parity conjecture (Theorem \ref{relative_p_parity_over_Q}).
    In fact, as a novel application of   our main theorem (Theorem \ref{relative_p_parity_over_Q}), for $\ell=p$, we have shown that $d(\mathfrak{F}_p^1, \mathfrak{F}_p^2) \equiv \delta(A_1,p) -\delta(A_2,p) \pmod2$ in Proposition \ref{leqp}. 

Note that for $\ell=p$, with appropriate extension of scalar we may assume that the residual representations are symplectically equivalent. However, it is not  clear that this isomorphism maps $A_1^-[\pi]$ to $A_2^-[\pi]$.
\qed 
\end{remark}

\begin{remark}
In \cite[\S 5.1.2]{Ne3}, Nekov\'a\v r considers a complete local Noetherian domain of dimension $2$ and studies the Galois representations arising from such a family. He refers to these as a parametric family of Galois representations and compares their root numbers and local Selmer conditions. The Galois representations attached to modular forms lying on a fixed branch of a Hida family (cf. \cite{epw}) provide a natural example of such a parametric family. However, if the modular forms lie on different branches of a Hida family, the universal deformation ring parametrising the resulting family need not be a domain  (cf. \cite{epw}). Our results also  apply to such families of Galois representations, parametrised by Noetherian local rings that need not be domains and need not be of Krull dimension equal to $2$. \qed
\end{remark}

\section{Selmer Groups}\label{Selmer Groups}  Let $V$ be a $p$-adic Galois representation of $G_\Q$ satisfying the hypothesis \hyperref[Hyp:Mov]{(Mov)}. Also assume that $V$ satisfies the  Panchishkin condition at $p$ (Definition \ref{panchishkin}). 
Recall $\FF$, $\OO$, $\pi$ be as before. Let $T$ be a $G_\Q$-stable  $\OO$-lattice  of $V$ and put $A:= V/T$ and let $A^+$ be the image of $V^+$ under the natural map $V\lra A$ and set $A^-:= A/A^+$. 
Let $\Sigma$ be a finite set of places of $\Q$ containing the prime $p$, the infinite prime and the primes where $V$ is ramified. Suppose $\Q_{\Sigma}$ is the maximal algebraic extension of $\Q$ unramified outside $\Sigma$.  For a number field $L\subset \Q_\Sigma$, let $v$ be a finite prime of $L$ and  $I_v$ denote the inertia subgroup of $G_{L_v}$.

 We have a natural map $H^1(L_v,V) \xrightarrow{p_*} H^1(L_v,A)$ induced by the natural projection map $V \xrightarrow{p} A$. Then define  the 
Bloch-Kato local condition for $A$ at $v$,
\begin{equation}\label{global_BK_selmer_group}
    H ^{\BK}(L_v,\ A): = p_{*}(H^\BK (L_v,V)),
\end{equation}
where $H^\BK (L_v,V)$ is defined in equation \eqref{local_BK_selmer}. 
Next, for a finite rational prime $\ell$ and a number field $L$ and $X \in \{V,A\}$, define 
{\small
\begin{equation}
\label{local_condition_strict}
\begin{aligned}
      H ^\str (L_v, X)&:=\left\{ \begin{array}{lc}
        \Ker\left(H^1(L_v,X)\rightarrow H^1(I_v,X)\right)  & \text{if}\;\; v \nmid p,  \\
       \Ker\left(H^1(L_v,X)\rightarrow H^1(L_v, X^-)\right) & \text{if}\;\; v \mid p.
     \end{array}\right.
\end{aligned}
\end{equation}
}

For  $\Sigma_0\subset \Sigma \backslash \{p,\infty\}$ and $\dag \in \{ \str, \BK\}$, define the $\Sigma_0$-imprimitive $\dag$ Selmer group $S_{\dag}^{\Sigma_0}(X/{L})$ of $X$ over $L$ as 
\begin{equation}
\label{Imprimitive_dag_selmer_group}
\begin{aligned}
    S_{\dag}^{\Sigma_0}(X/{L}) :=\Ker\Big( H^1(\Q_\Sigma/{L},X) \rightarrow \prod_{v|\ell,\ \ell\in \Sigma\backslash\Sigma_0} \frac{H^1(L_v,X)}{H^\dag(L_{v},X)}\Big).
\end{aligned}
\end{equation}

Note that the $\dag$ Selmer group $S_{\dag}^{\Sigma_0} (X/L)$ is independent of the choice of $\Sigma$ as long as $\Sigma$ contains $p,\ \infty$ and the primes where $V$ is ramified. In fact, as $p$ is odd, the local conditions at the infinite primes do not play a role in the definition of Selmer group. Also note that $\rank_\OO\ S_{\dag}^{\Sigma_0} (A/L)^\vee = \rank_\FF\ S_{\dag}^{\Sigma_0} (V/L)$.



If $\Sigma_0$ is empty, denote $ S_{\dag}^{\Sigma_0}(-/{L})$ as $S_{\dag}(-/{L})$.
Next,  the $\Sigma_0$-imprimitive strict selmer group $S_{\str}^{\Sigma_0}(A[\pi]/{L})$ of $A[\pi]$ over ${L}$ is defined by
\begin{equation*}
    S_{\str}^{\Sigma_0}(A[\pi]/{L}) :=\Ker\Big( H^1(\Q_\Sigma/{L},A[\pi]) \rightarrow \prod_{v|\ell,\ \ell\in \Sigma\backslash\Sigma_0} \frac{H^1(L_v,A)}{H^\str(L_{v},A)}\Big).
\end{equation*}  
For each of the Selmer groups defined above, the definition extends to the cyclotomic $\Z_p$-extension $L_\infty$ of $L$ by taking the direct limit of the corresponding Selmer group over $L_n$, as $n$ varies. 


\subsection{Comparison between Bloch-Kato and strict Selmer groups} Let $V$ be a $p$-adic Galois representation of $G_\Q$ and $L$ be a number field, as defined in \S \ref{Selmer Groups}.

Let $v$ be a prime of $L$ not dividing $p$ and $w$ be a prime of $L_\infty$ dividing $v$.
Then  $L_{\infty,w}$ is an 
unramified $\mathbb{Z}_p$-extension of $L_v$, so $\mathrm{Gal}(L_v^{\mathrm{ur}}/L_{\infty,w})$ has pro-order prime to $p$. Thus,
\[
  H^{\mathrm{str}}(L_{\infty,w},A)\cong H^1\bigl(\mathrm{Gal}(L_v^{\mathrm{ur}}/L_{\infty,w}),A^{I_w}\bigr)=0,
\]
and hence $H^{\mathrm{BK}}(L_{\infty,w},A)$ being a subgroup of  $H^{\mathrm{str}}(L_{\infty,w},A)$ vanishes  as well. Combining this with Theorem~\ref{str_bk_at_p} for the primes above $p$, we  we arrive at the following theorem:

\begin{theorem}\label{BK_str_equality_over_cyclotomic}
     Let $V$ be a symplectic self-dual $p$-adic representation of $G_\Q$ satisfying \hyperref[Hyp:Mov]{(Mov)}. Assume that the representation $V$ satisfies the Panchishkin condition at $p$ and is pure (of weight $-1$).
     Then for any number field $L$,
     \begin{equation*}
S_\BK(A/L_\infty)= S_\str(A/L_\infty). \tag*{\qed}
\end{equation*}
   \end{theorem}

\begin{remark}
For a prime $v \nmid p$ of $L$, the two conditions $H^{\mathrm{BK}}(L_v,A)$ and $H^{\mathrm{str}}(L_v,A)$ may differ. Indeed, for a prime  $v\nmid p$, $H^{\mathrm{BK}}(L_v,A)=0$ by purity, whereas
$H^{\mathrm{str}}(L_v,A)\cong\bigl(A^{I_v}/(A^{I_v})_{\mathrm{div}}\bigr)_{\mathrm{Frob}_v}$ is a finite
group of order prime to $p$, but need not be zero.
\end{remark}
We end this subsection with the following lemma which will be used later.
\begin{lemma}\label{vanishing_V_m}
   Let $V$ be a $p$-adic representation of $G_\Q$ satisfying \hyperref[Hyp:Mov]{(Mov)}. Assume that $V$ is pure of weight $-1$, then $H^0(\Q_\infty, V) = 0$.  
\end{lemma}
\begin{proof}
Suppose, to the contrary, that $H^0(\Q_\infty, V)\neq 0$. Let $L$ be a number field such that  $V$ is a semistable representation of $G_{L_v}$ for a prime $v\mid p$ and $L_v= L_v^\ur \cap L_v(\mu_{p^\infty})$. As $H^0(\Q_\infty, V)\neq 0$, we have $H^0(L_{\infty,v}, V)\neq 0$ and using Sen's theory (see Lemma \ref{Sen_theory}), we have $H^0(L_{\infty,v}, V)\cong \bigoplus_{k=1}^r\chi_p^{i_k}$ as $\Gal(L_{\infty,v}/L_v)$-module with $i_k\in \Z$. Consider a one dimensional subspace $U$ of $H^0(L_{\infty,v}, V)$ such that  $\Gal(L_{\infty,v}/L_v)$ acts on $U$ via $\chi_p^{i_k}$. As $L_\infty/L$ is totally ramified extension at the prime $v$, the decomposition subgroup $\Gal(L_{\infty,v}/L_v)$ is the same as $\Gal(L_\infty/L)$ and consequently, $\Gal(L_\infty/L)$ acts on $U$ via $\chi_p^{i_k}$. Choose a prime
$w \nmid p$ of $L$ (with  $V$  is unramified at $w$), then $G_{L_w}$ acts on $U$ via $\chi_p^{i_k}$ and the  complex absolute value $|\chi_p(\Frob_\ell)^{i_k}|=\ell^{i_k}$. In particular, the action of $\Frob_\ell$ on $V$ has one of the eigenvalues with complex absolute value $\ell^{i_k}$. On the other hand,  as $V$ is pure of weight $-1$, all the   eigenvalues of $\Frob_\ell$ on $V$ should be of complex absolute value $\ell^{\frac{1}{2}}$, which is a contradiction as $i_k\in \Z$. 
Hence $H^0(L_{\infty,w}, V)=0$ and consequently, $H^0(\Q_\infty, V)= 0$.
\end{proof}

\subsection{Comparison of Selmer coranks} Consider two $p$-adic Galois representations $V_1$ and $V_2$ of $G_\Q$ satisfying \hyperref[Hyp:Mov]{(Mov)}  with coefficient field $\FF$. We assume the residual representations $A_1[\pi]\cong A_2[\pi]$ as $G_\Q$-modules and this isomorphism induces an isomorphism $A_1^-[\pi] \cong A_2^-[\pi]$ of $G_{\Q_p}$ modules. Then under suitable conditions, we relate the difference of the $\OO$-ranks of $S_\str(A_1/\Q_\infty)^\vee$ and $S_\str(A_2/\Q_\infty)^\vee$
in terms of the same local invariants $\delta(A_1,\ell)$ and $\delta(A_2,\ell)$,  $\ell \in \Sigma $, that appeared in comparison of the root numbers of $V_1$ and $V_2$ (see Theorem \ref{global_root_lambda}).

Recall $L_\infty$ denotes the cyclotomic $\mathbb{Z}_p$-extension of $L$ with $\Gamma_L= \Gal(L_\infty /L)\cong \mathbb{Z}_p$. The Iwasawa algebra $\Lambda(\Gamma_L)$ is defined by  $\Lambda(\Gamma_L):=\varprojlim_{n} \OO[\Gal(L_n/L)]$. By choosing a topological generator $\gamma$ of $\Gamma_L$, we can identify $\Lambda(\Gamma_L)$ with $\OO [[X]]$, via the map $\gamma \longmapsto 1+X$. For brevity, we write $\Gamma=\Gamma_\Q$ and $\Lambda=\Lambda(\Gamma_\Q)$. Let $\Sigma$ be a finite set of places of $\Q$ containing the prime $p$, the infinite prime and the primes where $V$ is ramified. It is well known that $H^1(\Q_\Sigma/\Q_\infty, A)^\vee$ (and hence $S_{\str}(A/\Q_\infty)^\vee$)  is a finitely generated $\Lambda$-module. By the  structure theorem of finitely generated modules over $\Lambda\cong \OO [[X]]$, if $S_{\str}(A/\Q_\infty)^\vee$ is $\Lambda$-torsion and the Iwasawa $\mu$-invariant, $\mu(S_{\str}(A/\Q_\infty)^\vee)$ vanishes, then $S_{\str}(A/\Q_\infty)^\vee$ is a finitely generated module over the DVR $\OO$.

Recall,  from equations \eqref{delta_at_l} and \eqref{delta_at_p}, for any finite prime $\ell$ of $\Q$,
\begin{equation*}
    \delta(A_m,\ell)= \left\{ \begin{array}{ccc}
       \rk_\mathcal{O}\  \Big( \prod_{\eta\mid \ell} H^1(\Q_{\infty,\eta} , A_m) ^\vee \Big)  & \text{if} & \ell\neq p, \\
        \dim_\mathfrak{f}\ \frac{H^\str(\Q_{\infty,p}, A_m[\pi])}{H^\str(\Q_{\infty,p}, A_1[\pi])\cap H^\str(\Q_{\infty,p}, A_2[\pi])} & \text{if} & \ell =p.
    \end{array}\right.
\end{equation*}


\begin{proposition}
\label{pi_inside_outside} Let $V_1$ and $V_2$ be two finite dimensional $\FF$-vector spaces endowed with a continuous $G_\Q$ action satisfying \hyperref[Hyp:Mov]{(Mov)} such that $A_1[\pi] \cong A_2[\pi]$ as $G_\Q$-modules. We assume that $V_1$ and $V_2$ are symplectic self-dual representations (i.e., $V_i\cong V_i^*(1)$ for $i=1,2$) and they satisfy the Panchishkin condition at $p$ and pure (of weight $-1$). Let $\Sigma_0$ be a subset of $\Sigma\backslash\{p,\infty\}$ containing all the primes (excluding $p$), where either $V_1$ or $V_2$ is ramified. We make the following additional hypotheses.
\begin{enumerate}
    \item \label{hyp:H1}  $S_{\BK}(A_1/\Q_\infty)^\vee$ is $\Lambda$-torsion and  $\mu(S_{\BK}(A_1/\Q_\infty)^\vee)$ is equal to $0$. 
    \item \label{hyp:H4} The given isomorphism $A_1[\pi] \cong A_2[\pi]$ when restricted to $G_{\Q_p}$, induces an isomorphism $A_1^-[\pi] \cong A_2^-[\pi]$.
\end{enumerate}
\begin{equation*}\label{thm1eq}
        \text{Then }
     \ \mathrm{rk}_\mathcal{O}\ S_{\BK}^{\Sigma_0}(A_1/\Q_\infty)^\vee - \delta(A_1, p) = \mathrm{rk}_\mathcal{O}\ S_{\BK}^{\Sigma_0}(A_2/\Q_\infty)^\vee - \delta(A_2,p).
    \end{equation*}
\end{proposition}


\n \textbf{Proof:} 
By Theorem \ref{BK_str_equality_over_cyclotomic}, we have $S_\BK(A_1/\Q_\infty) = S_\str (A_1/\Q_\infty)$ and so $S_{\str}(A_1/\Q_\infty)^\vee$ is $\Lambda$-torsion and $\mu(S_{\str}(A_1/\Q_\infty)^\vee)=0$. For $m \in \{1,2\}$ and  $X\in \{A_m, A_m[\pi]\}$, set
\begin{equation*}
\G_\ell ^\str (\Q_\infty,X):=\left\{  \begin{array}{lc}
       \prod_{\eta \mid \ell} H^1(\Q_{\infty,\eta}, X)  & \text{if $\ell\neq p$,}  \\
       H^1(\Q_{\infty,p}, X^-)  & \text{if $\ell=p$.}
    \end{array} \right.
\end{equation*}
Then multiplication by $\pi$ on the $\pi$-divisible module $A_m$ 
induces the  commutative diagram 
{\small 
\begin{equation}
\label{cd}
\begin{tikzcd}[ampersand replacement=\&]
	0 \& {S_{\str}^{\Sigma_0}(A_m[\pi]/\Q_\infty)} \& {H^1(\Q_\Sigma/\Q_\infty,A_m[\pi])} \& {\prod_{\ell\in \Sigma\backslash\Sigma_0}\G_\ell^\str(\Q_\infty, A_m[\pi])} \\
	0 \& {S_{\str}^{\Sigma_0}(A_m/\Q_\infty)[\pi]} \& {H^1(\Q_\Sigma/\Q_\infty,A_m)[\pi]} \& {\prod_{\ell\in \Sigma\backslash\Sigma_0}\G_\ell^\str(\Q_\infty, A_m)[\pi].}
	\arrow[from=1-1, to=1-2]
	\arrow[from=1-2, to=1-3]
	\arrow["\psi_m", from=1-2, to=2-2]
	\arrow["\varrho_m", from=1-3, to=1-4]
	\arrow["\varphi_m",from=1-3, to=2-3]
	\arrow["h_m"', from=1-4, to=2-4]
	\arrow[from=2-1, to=2-2]
	\arrow[from=2-2, to=2-3]
	\arrow["\theta_m", from=2-3, to=2-4]
\end{tikzcd}
\end{equation}
}
By hypothesis \hyperref[hyp:H4]{(2)}, we have $S_{\str}^{\Sigma_0}(A_1[\pi]/\Q_\infty)\cong S_{\str}^{\Sigma_0}(A_2[\pi]/\Q_\infty)$ and by construction, the map $\varphi_m$ is  surjective with finite kernel.  Note that $\Ker(h_m)$ is finite. Thus it is easy to see that condition \hyperref[hyp:H1]{(1)} implies $S_{\str}(A_2/\Q_\infty)^\vee$ is also $\Lambda$-torsion and $\mu(S_{\str}(A_2/\Q_\infty)^\vee)=0$.

Next, we  claim that  $S_\str ^*(A_m[\pi]/\Q_\infty):= \varprojlim_n S_\str ^{\Sigma_0}(A_m[\pi]/\Q_n)=0$.   Indeed, under assumption \hyperref[hyp:H1]{(1)}, we have the following exact sequence (cf. \cite[1.3.a]{Greenberg_projectivemodule_modularrepresentation})
\begin{equation}
\label{strict_selmer_for_Am}
      0 \rightarrow S_{\str}(A_m/\Q_\infty) \rightarrow H^1(\Q_\Sigma/\Q_\infty,A_m) \rightarrow \prod_{\ell\in \Sigma}\G_\ell^{\str}(\Q_\infty,A_m) \rightarrow 0.
\end{equation}
Thus we obtain,
\begin{equation}\label{sigma_sigma_0}
    S_\str^{\Sigma_0}(A_m/\Q_\infty)/ S_\str(A_m/\Q_\infty) \simeq \prod_{\ell\in \Sigma_0} \G_\ell^\str (\Q_\infty,A_m).
\end{equation}
 Now, for each $\ell\in \Sigma_0$, $\G_\ell^\str(\Q_\infty, A_m)$ is a cofinitely generated $\mathcal{O}$-module, therefore from \hyperref[hyp:H1]{(1)},
 $S_\str^{\Sigma_0}(A_m/\Q_\infty)^\vee$ is also $\Lambda$-torsion, and   $\mu(S_\str^{\Sigma_0}(A_m/\Q_\infty)^\vee)=0$. In particular,  $S_\str^{\Sigma_0}(A_m/\Q_\infty)[\pi]$ is finite and applying a Snake lemma in  \eqref{cd}, we deduce that $S_\str^{\Sigma_0}(A_m[\pi]/\Q_\infty)$ is finite.
 Then it follows  from  \cite[Lemma 5.5]{Lim-Sujatha_JNT_2018}  
 that $S_\str^*(A_m[\pi]/\Q_\infty)=0$, as claimed.  By the Poitou-Tate duality, $\Coker(\varrho_m)$ injects into $S_\str ^*((A_m[\pi])^*(1)/\Q_\infty)$ and by the self-duality of $A_m$, we conclude that $\varrho_m$ is surjective. 

Now we will show that $S_\str^{\Sigma_0}(A_m/\Q_\infty)^\vee$ has no non-zero pseudo-null $\Lambda$-submodule. Since $\varrho_m$ and $h_m$ are both surjective, it is  immediate that $\theta_m$ is also surjective.
Next  we again have  a short exact sequence, similar to \eqref{strict_selmer_for_Am}, for the defining global to local map of   $S_\str^{\Sigma_0}(A_m/\Q_\infty)^\vee$, as it is a $\Lambda$-torsion module. 
 Thus, we get $0 \rightarrow   \frac{S_\str^{\Sigma_0}(A_m/\Q_\infty)}{\pi S_\str^{\Sigma_0}(A_m/\Q_\infty)} \hookrightarrow   \frac{H^1(\Q_\Sigma/\Q_\infty,A_m)}{\pi H^1(\Q_\Sigma/\Q_\infty,A_m)}.$ Under the assumption \hyperref[hyp:H1]{(1)}, it is known that (see for example, \cite[Proposition 4.4]{Sujatha_elliptic_curve_iwasawa_mu0_conjecture}) $H^2(\Q_\Sigma/\Q_\infty, A_m[\pi]) =0$. Hence we deduce that $\frac{H^1(\Q_\Sigma/\Q_\infty,A_m)}{\pi H^1(\Q_\Sigma/\Q_\infty,A_m)} =0=\frac{S_\str^{\Sigma_0}(A_m/\Q_\infty)}{\pi S_\str^{\Sigma_0}(A_m/\Q_\infty)}$. Thus $S_\str^{\Sigma_0}(A_m/\Q_\infty)$ is divisible and consequently, $S_\str^{\Sigma_0}(A_m/\Q_\infty)^\vee$  has no non-zero pseudo-null $\Lambda$-submodule. Hence
 {\small
 \begin{equation}
\label{pi_in_out_selmer_rank_compute}
 \rk_\mathcal{O}\ S_\str^{\Sigma_0}(A_m/\Q_\infty)^\vee = \rk_\mathfrak{f}\ S_\str^{\Sigma_0}(A_m/\Q_\infty)[\pi] = \rk_\mathfrak{f} \ S_\str^{\Sigma_0}(A_m[\pi]/\Q_\infty) + \rk_\mathfrak{f}\ \Coker(\psi_m) - \rk_\mathfrak{f}\ \Ker(\psi_m).
\end{equation}}

Moreover, by the surjectivity of $\varrho_m$ and using the Snake lemma in \eqref{cd}, we get the following exact sequence:
\begin{equation}\label{4_term}
    0 \lra \Ker (\psi_m) \lra \Ker(\varphi_m) \lra \Ker(h_m) \lra \Coker(\psi_m) \lra 0.
\end{equation}
As all terms in this exact sequence are of finite order, we have 
\begin{equation*}
    \rk_\mathfrak{f}\ \Coker(\psi_m) - \rk_\mathfrak{f}\ \Ker(\psi_m)= \rk_\mathfrak{f}\ \Ker(h_m) - \rk_\mathfrak{f}\ \Ker(\varphi_m).
\end{equation*}
Observe that $\Ker(\varphi_m)= H^0(\Q_{\Sigma}/\Q_\infty, A_m)/\pi$.  By Lemma \ref{vanishing_V_m},  $H^0(\Q_{\Sigma}/\Q_\infty, A_m)$ is finite and we obtain  $\rk_\mathfrak{f}\  H^0(\Q_{\Sigma}/\Q_\infty, A_m)/\pi= \rk_\mathfrak{f}\ H^0(\Q_{\Sigma}/\Q_\infty, A_m)[\pi]=\rk_\mathfrak{f}\   H^0(\Q_{\Sigma}/\Q_\infty, A_m[\pi])$, which is the same for both $A_1$ and $A_2$ and hence 
\begin{equation}\label{phi_equality}
    \rk_\mathfrak{f}\ \Ker(\varphi_1)= \rk_\mathfrak{f}\ \Ker(\varphi_2).
\end{equation}
 Let $\eta_\ell$ be a prime of $\Q_\infty$ dividing $\ell$, where $\ell\in \Sigma\backslash\Sigma_0$, $\ell\neq p$. As $G_{\Q_{\infty,\eta_\ell}}/I_{\eta_\ell}$ is a  profinite group of order prime to $p$ and $A_m$ is unramified at $\eta_\ell$,    $H^0(G_{\Q_{\infty,\eta_\ell}}, A_m)$ is $\pi$-divisible and we compute,  $\Ker(h_m) \simeq \frac{H^0(\Q_{\infty,p},\ A_m^-)}{\pi}$. 
Now, from \eqref{pi_in_out_selmer_rank_compute}, we have the following  identity
{\small
\begin{equation}\label{prop_equn_4}
\begin{aligned}
\mathrm{rk}_{\mathcal{O}}\ S_\str^{\Sigma_0}(A_m/\Q_\infty)^\vee + \mathrm{rk}_\mathcal{O}\ (A_m^-)_{G_{\Q_{\infty,p}}}^\vee & = \rk_\mathfrak{f}\ S_\str^{\Sigma_0}(A_m[\pi]/\Q_\infty) + \rk_\mathfrak{f}\ \Ker(h_m)  +\mathrm{rk}_\mathcal{O}\ (A_m^-)_{G_{\Q_{\infty,p}}}^\vee - \rk_\mathfrak{f}\ \Ker(\varphi_m)
    \\
     &= \rk_\mathfrak{f}\ S_\str^{\Sigma_0}(A_m[\pi]/\Q_\infty) + \rk_\mathfrak{f}\ \left((A_m^-)^{G_{\Q_{\infty,p}}}[\pi]\right) - \rk_\mathfrak{f}\ \Ker(\varphi_m).  
\end{aligned}
\end{equation}}
Recall from Theorem \ref{BK_str_equality_over_cyclotomic}, we have 
\begin{equation}\label{prop_equn_6}
    S_\str^{\Sigma_0}(A_m/\Q_\infty)= S_\BK^{\Sigma_0}(A_m/\Q_\infty).
\end{equation}
On the other hand, from equation \eqref{delta_rank_at_p}, we know that 
\begin{equation}\label{prop_equn_5}
\mathrm{rk}_\mathcal{O}\ (A_1^-)_{G_{\Q_{\infty,p}}}^\vee - \mathrm{rk}_\mathcal{O}\ (A_2^-)_{G_{\Q_{\infty,p}}}^\vee = \delta(A_2,p)- \delta(A_1, p).
\end{equation}
Now the proposition follows from equations \eqref{phi_equality}, \eqref{prop_equn_4}, \eqref{prop_equn_6} and  \eqref{prop_equn_5}. \qed


\begin{corollary}
\label{str_imprimitive_str_compare}
    We keep the setting and hypotheses of Proposition \ref{pi_inside_outside}. Then using \eqref{sigma_sigma_0}, we further deduce: 
\begin{align*}
 \mathrm{rk}_\mathcal{O}\ S_{\BK}(A_1/\Q_\infty)^\vee & + \sum_{\ell\in \Sigma_0} \delta(A_1,\ell) - \delta(A_1,p)= \mathrm{rk}_\mathcal{O}\ S_{\BK}(A_2/\Q_\infty)^\vee + \sum_{\ell\in \Sigma_0} \delta(A_2,\ell) - \delta(A_2,p).\tag*{\qed}
\end{align*}
\end{corollary}

\section{Relative $p$-parity conjecture}\label{relative_p_parity_conjecture}
In this section, we establish  the relative $p$-parity conjecture under suitable hypotheses. At first, we relate the ratio of global root numbers with the difference of coranks of  Bloch-Kato Selmer groups over the cyclotomic $\Z_p$ extension.
\begin{theorem}\label{relative_p_parity_over_cycyclotomic}
         Let $V_1$, $V_2$ be two $p$-adic representations of $G_\Q$ satisfying \hyperref[Hyp:Mov]{(Mov)}, that are  symplectic self-dual, pure (of weight $-1$) and Panchishkin at $p$, and satisfy $A_1[\pi]\cong A_2[\pi]$ as $G_\Q$-modules. We make the following additional assumptions.
     \begin{enumerate}
         \item $S_\BK(A_1/\Q_\infty)^\vee$ is $\Lambda$-torsion and $\mu(S_\BK(A_1/\Q_\infty)^\vee)$ is equal to $0$.
    \item  The given isomorphism $A_1[\pi] \cong A_2[\pi]$, when restricted to $G_{\Q_p}$,  induces an isomorphism $A_1^-[\pi] \cong A_2^-[\pi]$.
     \end{enumerate}
\begin{equation*}
  \text{ Then  }  \ \quad   \frac{W(V_1)}{W(V_2)}= (-1)^{\mathrm{rk}_\mathcal{O}\ S_{\BK}(A_1/\Q_\infty)^\vee- \mathrm{rk}_\mathcal{O}\ S_{\BK}(A_2/\Q_\infty)^\vee}.
     \end{equation*}
\end{theorem}
\begin{proof}
The result follows from applying Theorem \ref{global_root_lambda} in Corollary \ref{str_imprimitive_str_compare}.
\end{proof}

\subsection{Comparison of Selmer corank over $\Q$ and $\Q_\infty$}
\label{Existence_of_pairing}
Now we will complete the proof of relative $p$-parity conjecture by comparing the $\mathcal O$-ranks of the dual Bloch-Kato Selmer groups $S_{\BK}(A_m/\Q)^\vee$ and $S_{\BK}(A_m/\Q_\infty)^\vee$ for $m\in \{1,2\}$, under the assumptions of Theorem \ref{relative_p_parity_over_cycyclotomic}. 

As $V_m$ is symplectic self-dual,  we have a $G_\Q$-invariant, non-degenerate  skew-symmetric bilinear pairing $V_m\times V_m\rightarrow \mathcal{F}(1)$. For any number field ${L}$, set $\Sha(A_m/{L}):= \frac{S_{\BK}(A_m/{L})}{S_{\BK}(A_m/{L})_{\mathrm{div}}}$, where $(-)_{\mathrm{div}}$ denotes the maximal $\pi$-divisible submodule of $(-)$. Flach \cite{Flach90}  has constructed a  non-degenerate skew-symmetric  pairing 
\begin{equation}
\label{pairing}
    \langle\; ,\; \rangle : \Sha(A_m/{L}) \times \Sha(A_m/{L}) \rightarrow\mathcal{F}/\mathcal{O}.
\end{equation}
Using this pairing, we can now complete the proof of relative $p$-parity conjecture stated in Theorem \ref{relative_p_parity_over_Q}.

\textbf{Proof of Theorem \ref{relative_p_parity_over_Q}:}
 It is clear from Theorem \ref{relative_p_parity_over_cycyclotomic}, that we only have to prove for $m\in \{1,2\}$,
\begin{equation*}
\mathrm{rk}_{\mathcal{O}}\ S_{\BK}(A_m/\Q_\infty)^\vee \equiv \mathrm{rk}_{\mathcal{O}}\ S_{\BK}(A_m/\Q)^\vee \pmod2.
\end{equation*}
 Since $S_\BK(A_m/\Q_\infty)^\vee$ is $\Lambda$-torsion and  $\mu(S_\BK(A_m/\Q_\infty)^\vee)$ is equal to $0$, $S_{\BK}(A_m/\Q_\infty)^\vee $ is a  finitely generated $\mathcal{O}$-module. 
Observe that the kernel of  the natural restriction map $S_{\BK}(A_m/\Q_n)\longrightarrow S_{\BK}(A_m/\Q_\infty)$ is finite and bounded independent of $n$ by Lemma \ref{vanishing_V_m} and hence $\mathrm{rk}_\mathcal{O}\ (S_{\BK}(A_m/\Q_n))^\vee$ is bounded by $\lambda:= \mathrm{rk}_\mathcal{O}\ (S_{\BK}(A_m/\Q_\infty))^\vee$. Then there exists some $n_0\in \mathbb{N}$ such that $\mathrm{rk}_\mathcal{O}\ (S_{\BK}(A_m/\Q_n))^\vee =\lambda'$ for all $n\geq n_0$. 
Put $T_n:= S_{\BK}(A_m/\Q_n)_{\mathrm{div}}$  and set $T_\infty:=\varinjlim_n  T_n$. 
Clearly, the restriction map $T_0\rightarrow T_{n_0}^{\Gal(\Q_{n_0}/\Q)}$ has  finite kernel and finite cokernel.  We can write $T_{n_0}^\vee\otimes\FF $ as direct sum of $\Gal(\Q_{n_0}/\Q)$-modules,  $T_{n_0}^\vee \otimes\FF= (T_{n_0}^\vee \otimes\FF) ^{\Gal(\Q_{n_0}/\Q)} \bigoplus X$, where $\Gal(\Q_{n_0}/\Q)$ acts nontrivially on $X$. Now by the self-duality property of $S_{\BK}(A_m/\Q_{n_0})$ (see \cite[(3.3)]{Ne5}) and our assumption that $p$ is an odd prime, we conclude that $\mathrm{rk}_\mathcal{O}\ T_0^\vee \equiv \mathrm{rk}_{\mathcal{O}}\ T_{n_0}^\vee \pmod 2$.
Thus we obtain 
\begin{equation*}
     \mathrm{rk}_\mathcal{O}\ S_{\BK}(A_m/\Q)^\vee=\mathrm{rk}_\mathcal{O}\ T_0^\vee \equiv \mathrm{rk}_{\mathcal{O}}\ T_{n_0}^\vee = \lambda' \pmod2 .
\end{equation*}
 Put $\Sha(A_m/{\Q_\infty}): =\varinjlim_n \Sha(A_m/{\Q_n})$. Then we have a commutative diagram 
 \begin{equation}
 \label{selmer_sha_cd}
 \begin{tikzcd}[ampersand replacement=\&]
	0 \& {T_{n}} \& {S_{\BK}(A_m/\mathbb{Q}_{n})} \& \Sha(A_m/{\Q_n}) \& 0 \\
	0 \& {T_\infty} \& {S_{\BK}(A_m/\mathbb{Q}_\infty)} \& {\Sha(A_m/{\Q_\infty})} \& 0
	\arrow[from=1-1, to=1-2]
	\arrow[from=1-2, to=1-3]
	\arrow[from=1-2, to=2-2]
	\arrow[from=1-3, to=1-4]
	\arrow[from=1-3, to=2-3]
	\arrow[from=1-4, to=1-5]
	\arrow[from=1-4, to=2-4]
	\arrow[from=2-1, to=2-2]
	\arrow[from=2-2, to=2-3]
	\arrow[from=2-3, to=2-4]
	\arrow[from=2-4, to=2-5]
\end{tikzcd}
\end{equation}
Consider $n\geq n_0$. Then it is plain that the restriction map $T_{n}\rightarrow T_\infty$ is surjective and kernel of this map is finite and bounded independent of $n$.
 Also $\Ker\big(S_{\BK}(A_m/\mathbb{Q}_{n})\rightarrow S_{\BK}(A_m/\mathbb{Q}_\infty)\big)$ is finite and is bounded independent of $n$.
By using a Snake Lemma in \eqref{selmer_sha_cd}, we deduce that 
$\Ker\big(\Sha(A_m/{\Q_n})\rightarrow \Sha(A_m/{\Q_\infty})\big)$ is also finite and bounded independent of $n \geq n_0$.
 Since we have a non-degenerate, skew-symmetric pairing $\Sha(A_m/{\Q_n})\times \Sha(A_m/{\Q_n})\rightarrow  \mathcal{F}/\mathcal{O}$ (see \eqref{pairing}), it follows from \cite[Lemma 3 and Prop. 6]{Guo_Li} that $\mathrm{rk}_\mathcal{O}\ \Sha(A_m/{\Q_\infty})^\vee $ is even. Combining all these, we have 
\begin{equation*} \lambda= \mathrm{rk}_\mathcal{O}\ S_{\BK}(A_m/\Q_\infty)^\vee  \equiv \mathrm{rk}_\mathcal{O}\ T_\infty^\vee 
    \equiv\mathrm{rk}_\mathcal{O}\ T_{n_0}^\vee 
\equiv \mathrm{rk}_\mathcal{O}\ S_{\BK}(A_m/\Q)^\vee  \pmod2. \qed 
\end{equation*}

 \begin{remark}\label{relative_p_parity_over_L}
For simplicity, we have presented Theorem \ref{relative_p_parity_over_Q} over $\Q$. Using similar arguments, an analogue of Theorem \ref{relative_p_parity_over_Q}  remains valid over  arbitrary number field  $L$ instead of $\Q$, with parallel hypotheses over $L$ and cumbersome notation. 
\end{remark}

\subsection{Arithmetic local constant at $p$}
Under the setting of Theorem \ref{relative_p_parity_over_Q}, assume in addition that there is an isomorphism $A_1[\pi] \cong A_2[\pi]$  as symplectic $G_\Q$-modules. Then recall the definition of arithmetic local constant at $p$, $d(\mathfrak{F}_p^1, \mathfrak{F}_p^2)$ of  Mazur-Rubin-Nekov\'a\v r (see equation  \eqref{arithmetic_local-constant}). 
As an application of Theorem \ref{relative_p_parity_over_Q}, we can now relate $d(\mathfrak{F}_p^1, \mathfrak{F}_p^2)$ with the difference of the local Iwasawa invariants $\delta(-,p)$, defined in \eqref{delta_at_p}. 
\begin{proposition}\label{leqp}
We keep the setting and hypotheses of Theorem \ref{relative_p_parity_over_Q}. In addition assume that $A_1[\pi]$ is isomorphic to  $A_2[\pi]$  as symplectic $G_\Q$-modules.
\begin{equation*}
     \text{Then} \quad d(\mathfrak{F}_p^1, \mathfrak{F}_p^2) \equiv \delta(A_1,p) -\delta(A_2,p) \pmod2.
 \end{equation*}
\end{proposition}
\n \textbf{Proof:} Proof is immediate from Theorem \ref{relative_p_parity_over_Q}, Proposition \ref{lnp} and equation \eqref{MRtheorem}. \qed

\subsection{Arithmetic Examples}\label{arithmetic_examples} Now we will give some examples of classes of $p$-adic representations  for which we can apply our relative $p$-parity conjecture stated in Theorem \ref{relative_p_parity_over_Q}. 

\subsubsection{Modular forms:}\label{example_modular-forms}
 Let  $f=\sum_n a_n(f)q^n \in S_{k}(\Gamma_0(N))$ be a elliptic cuspidal newform of even weight $k\geq 2$, level $N$ and trivial nebentypus.
 Let $\Q_{f}: = \Q(\{a_n(f)\}_{n\geq 1})$ be the Hecke field of $f$. Let  $\mathcal{F}$ be the $p$-adic completion of  $\Q_{f}$ and $$\rho_{f} : G_\Q \rightarrow \mathrm{Aut}(V_{f})\cong \mathrm{GL}_2(\FF)$$ be  Deligne's  $p$-adic Galois representation associated with $f$, where $V_{f}$ is a two dimensional vector space over $\mathcal{F}$.  For every prime $\ell \nmid Np$,  $\rho_{f}$ satisfies 
$\mathrm{Tr}(\rho_{f}(\mathrm{Frob}_\ell))= a_\ell(f)$ and $\mathrm{det}(\rho_{f}(\mathrm{Frob}_\ell))=\chi_p(\mathrm{Frob}_\ell)^{k-1}=\ell^{k-1}$; here $\mathrm{Frob}_\ell$ denote the  arithmetic Frobenius at $\ell$.
Let us assume that  $f$ is ordinary at $p$, i.e. $\iota_p(a_p(f))$ is a $p$-adic unit. 
Let $T_{f}$ be a $G_\Q$ lattice of $V_{f}$ and set $A_{f}:= V_{f}/T_{f}$. Then 
 the representation $V_{f}$ of $G_\Q$ is semistable (hence de Rham) at $p$ and is pure of weight $1-k$. As $f$ is  $p$-ordinary, the representation $V_{f}$ satisfies the  Panchishkin condition at $p$ and let 
\begin{equation*}
    0\lra V_{f}^+ \lra V_{f} \lra V_{f}^- \lra 0
\end{equation*} 
 be the Panchishkin filtration of $V_{f}$. Let $A_f^+$ be the image of $V_f^+$ under the natural map $V_f \lra A_f$ and put $A_f^-:= A_f/ A_f^+$.


Let   $\sigma: G_\Q \longrightarrow \mathrm{GL}_n(\bar{\mathbb{Q}}_p)$ 
 be an irreducible orthogonal Artin representation (see \cite[\S9]{Greenberg_projectivemodule_modularrepresentation}). 
 Enlarging $\mathcal F$, if necessary, assume that $\mathcal{F}$  contains the coefficients of $\sigma$. 
 Set $V_f(\sigma):= V_{f}(1-\frac{k}{2}) \otimes_\FF W_\sigma$, where $W_\sigma$ is the representation space of $\sigma$.
 Since $f$ is even weight $k$ and of trivial nebentypus,  $V_f(1-\frac{k}{2})$ is a symplectic self-dual representation and together with  $W_\sigma$ is an orthogonal representation, we have $V_f(\sigma)$ is also a  symplectic self-dual representation of $G_\Q$.
 It is clear that $V_f(\sigma)$ is pure of weight $-1$ and is a de Rham representation of $G_{\Q_p}$.   The filtration
 \begin{equation*}
    0 \lra V_{f}^+(1-\frac{k}{2})\otimes_\FF W_\sigma \lra V_{f}(1-\frac{k}{2}) \otimes_\FF W_\sigma \lra V_{f}^-(1-\frac{k}{2}) \otimes_\FF W_\sigma \lra 0
\end{equation*}
gives  a Panchishkin filtration of $V_f(\sigma)$. Following the notation  in \eqref{panchishkin}, put $V_f(\sigma)^\pm := V_{f}^\pm(1-\frac{k}{2}) \otimes_\FF W_\sigma$.  For a $G_{\Q}$ stable $\OO$-lattice $L_\sigma$ of $W_\sigma$, put $A_f(\sigma):= A_{f}(1-\frac{k}{2}) \otimes_\OO L_\sigma$ and $A_f(\sigma)^\pm := A_{f}^\pm(1-\frac{k}{2}) \otimes_\OO L_\sigma$.

Let $f_1 \in S_{k_1}(\Gamma_0(N_1))$ and $f_2 \in S_{k_2}(\Gamma_0(N_2))$ be two $p$-ordinary elliptic cuspidal newforms with trivial nebentypus and  even weight $k_1$ and $k_2$, and level $N_1$ and $N_2$, respectively. 
Also assume that $\bar{\rho}_{f_1} \cong \bar{\rho}_{f_2}$ as $G_\Q$-modules. Note that $A_{f_m}^-[\pi]$ is the maximal unramified quotient of $A_{f_m}[\pi]$. This follows from the fact that $A_{f_m}^-[\pi]$ is unramified and inertia acts on $A_{f_m}^+[\pi]$ by $\omega_p^{k_m-1}$, which is nontrivial because $k_m$ is even and $p$ is odd. 
So the isomorphism $\bar{\rho}_{f_1} \cong \bar{\rho}_{f_2}$ when restricted to $G_{\Q_p}$, induces an isomorphism $A_{f_1}^-[\pi] \cong A_{f_2}^-[\pi]$. Moreover, if $k_1 \equiv k_2 \pmod{2(p-1)}$, then we immediately get the isomorphisms $A_{f_1}(\sigma)[\pi]\cong A_{f_2}(\sigma)[\pi]$ of $G_\Q$-modules and $A_{f_1}(\sigma)^-[\pi]\cong A_{f_2}(\sigma)^-[\pi]$ of $G_{\Q_p}$-modules.

Let $M$ be the Artin conductor of $\sigma$.
By Theorem \ref{global_root_lambda}, we deduce 
\begin{equation*}
   \frac{W(V_{f_1}(\sigma))}{W(V_{f_2}(\sigma))}=(-1)^{\sum_{\ell\mid  pMN_1N_2} \  (\delta(A_{f_1}(\sigma),\ \ell)- \delta(A_{f_2}(\sigma),\ \ell))}. 
\end{equation*}

Next, by a conjecture of Greenberg (see \cite[Conjecture 3.1]{Gr5}), for a  $p$-ordinary cuspidal newform $f_1$, $S_\BK(A_{f_1}(\sigma)/\Q_\infty)^\vee$ is $\Lambda$-torsion.
Thus the torsion condition in the hypothesis (1) of Theorem \ref{relative_p_parity_over_Q} is always expected to hold.
Note  that Kato \cite{Kato_p_adic_hodge_theory} proved that for a $p$-ordinary cuspidal newform $g$, 
$S_{\BK}(A_{f}/L_\infty)^\vee$ is $\Lambda(\Gamma_L)$-torsion, whenever  $L/\Q$ is abelian. 


It is conjectured by Greenberg (\cite[Conjecture 1.11]{Greenberg1997}) that if $E/\Q$ is an elliptic curve  with good ordinary reduction at $p$ then there exists a $\Q$-isogenous elliptic curve $E'$ with $\mu(S_{p^\infty}(E'/\Q_\infty)^\vee)$ $=0$. In particular, this implies  for a $p$-ordinary elliptic curve $E/\Q$, if $E[p]$ is an irreducible $G_\Q$-module, then $\mu(S_{p^\infty}(E/\Q_\infty)^\vee)=0$. For a general number field $L$,  in all known examples  where the $\mu$-invariant of a $p$-ordinary elliptic curve $E/L$ is positive,  the residual Galois representation of $E[p]$ is a reducible  $G_L$-module and  it is expected that for a  $p$-ordinary elliptic curve $E/L$ such that $E[p]$ is an irreducible $G_L$-module, the $\mu(S_{p^\infty}(E/\Q_\infty)^\vee)$  should vanish. Let $f_1, f_2, \sigma $  be as above and let $K$ be the fixed field of the kernel of the Artin representation $\sigma$. Then from the above discussion, it follows that if $\bar{\rho}_{f_1}$ is irreducible as a $G_K$-module, then we can expect $\mu(S_\BK(A_{f_1}/K_\infty)^\vee)=0$, which in turn implies that $\mu(S_\BK(A_{f_1}(\sigma)/\Q_\infty)^\vee)=0$.

In summary, under the assumption   that $\bar{\rho}_{f_1}$ is an irreducible  $G_K$-module, hypothesis (1) of Theorem \ref{relative_p_parity_over_Q} is conjectured to be satisfied and all other hypotheses of Theorem \ref{relative_p_parity_over_Q} are unconditionally verified to be  true and we see that the $p$-parity conjecture is true for $V_{f_1}(\sigma)$ if and only if it is true for $V_{f_2}(\sigma)$.




One of the most prominent examples of such congruence of modular forms comes from $p$-ordinary Hida family of modular forms.
For any newform $g$, let $g^\circ$ denote its $p$-stabilization
(see \cite[Lemma 3.3]{Hida_p_adic_measure_1}). If $p$
divides the level of $g$, then we use the convention $g^\circ=g$.

Let
$
f=\sum_{n\geq 1} a_n(f)q^n
\in S_{k}(\Gamma_0(N))
$
be a $p$-ordinary elliptic cuspidal newform of even weight $k\geq 2$ and trivial nebentypus, and suppose that its
residual Galois representation $\bar{\rho}_{f}$ is an irreducible $G_\Q$-module.
By Hida theory, $f^\circ$ belongs to a $p$-ordinary Hida family
\cite{Hida_Galois_representation_into_GL2_attached_ordinarycuspforms,
Wiles_ordinary_lambda_adic_representation}.
More precisely, there are infinitely many $p$-stabilized cuspidal newforms $f_n$ of weight $k_n$ such that $\bar{\rho}_{f_n}\cong \bar{\rho}_{f}$ as $G_\Q$-modules and moreover $k_n\equiv k \pmod{2(p-1)}$.

Note that we use the term \emph{Hida family} in the broad sense, allowing all
ordinary components associated with the fixed residual representation, as in \cite{epw}. In particular, the forms considered below
are not required to lie on the same branch of the Hida family and the universal deformation ring parameterizing the family need not be a domain.

Hence Theorem \ref{global_root_lambda} is applicable to infinitely many pairs $(V_f(\sigma), V_{f_m}(\sigma))$ such that $k_m\equiv k \pmod{2(p-1)}$  and
by the same Theorem we deduce 
\begin{equation*}
   \frac{W(V_{f}(\sigma))}{W(V_{f_m}(\sigma))}=(-1)^{\sum_{\ell\mid  pMN} \  (\delta(A_{f}(\sigma),\ \ell)- \delta(A_{f_m}(\sigma),\ \ell))}. 
\end{equation*}

Moreover, under the assumption   that $\bar{\rho}_{f}$ is an irreducible  $G_K$-module, hypothesis (1) of Theorem \ref{relative_p_parity_over_Q} is conjectured to be satisfied and all other hypothesis of Theorem \ref{relative_p_parity_over_Q} are unconditionally verified and we see that the relative $p$-parity conjecture is invariant for  $V_{f_m}(\sigma))$ for any $f_m$ in the Hida family  of $f$ with $k_m\equiv k \pmod{2(p-1)}$.

\begin{remark}
   There are several cases of elliptic curves $E/\Q$ whereas $E[p]$ is a reducible $G_\Q$-module but the Iwasawa $\mu$-invariant vanishes for the dual $p$-primary Selmer group over $\Q_\infty$. For some recent results, see \cite{ly}.
\end{remark}


\subsubsection{Twist by residually isomorphic Artin representations:}\label{twist_by_artin_representation}
Let $V$ be a $p$-adic representation of $G_\Q$ satisfying \hyperref[Hyp:Mov]{(Mov)} with a $G_\Q$-stable lattice $T$ and $A:=V/T$. Also assume that $V$ is a symplectic self-dual representation of $G_\Q$, pure of weight $-1$ and satisfies the Panchishkin condition at $p$.
For two Galois number fields $K_1$ and $K_2$, let $\sigma_m : \Gal(K_m/\Q) \lra \mathrm{GL}_n (\bar{\Q}_p)$ be two irreducible orthogonal Artin representations. Let $\FF$ be a $p$-adic completion of a  number field such that $V$, $\sigma_1$ and $\sigma_2$ are defined over $\FF$.
For $m\in \{1,2\}$, let $W_{\sigma_m}$ be the associated $\FF$ representation space and choose a $G_\Q$-stable lattice $L_{\sigma_m}$ of $W_{\sigma_m}$, which preserves the symmetric pairing on $W_{\sigma_m}$. 
Let the residual representations $\bar{\sigma}_1\cong \bar{\sigma}_2$  as $G_\Q$-modules. 
 Set $V(\sigma_m):= V\otimes_\FF W_{\sigma_m}$.  Then $V(\sigma_m)$ is also a symplectic self-dual representation of $G_\Q$ and pure of weight $-1$. As $V$ satisfies the Panchishkin condition at $p$, $V$ has a $G_{\Q_p}$ filtration given in \eqref{panchishkin}. Applying  $-{\otimes_\mathcal F}W_{\sigma_m}$, we obtain a Panchishkin filtration for $V(\sigma_m)$ as  $G_{\Q_p}$-modules:
\begin{equation*}
    0 \lra V^+\otimes_\FF W_{\sigma_m}\lra V\otimes_\FF W_{\sigma_m} \lra V^-\otimes_\FF W_{\sigma_m}\lra 0.
\end{equation*}
Let $A^+$ be the image of $V^+$ under the natural map $V \lra A$ and put $A^-:= A/ A^+$.
  Put $A(\sigma_m):= A \otimes_\OO L_{\sigma_m}$ and $A(\sigma_m)^\pm := A^\pm \otimes_\OO L_{\sigma_m}$. Also the isomorphism $\bar{\sigma}_1 \cong \bar{\sigma}_2$ of $G_\Q$-modules immediately gives the isomorphisms $A(\sigma_1)[\pi] \cong A(\sigma_2)[\pi]$ of $G_\Q$-modules and $A(\sigma_1)^-[\pi] \cong A(\sigma_2)^-[\pi]$ of $G_{\Q_p}$-modules.

Let $M_1$ and $M_2$ be the Artin conductors of $\sigma_1$ and $\sigma_2$, respectively and $N$ be the product of finitely many (finite) integer primes where $V$ is ramified.
By Theorem \ref{global_root_lambda}, we deduce 
\begin{equation*}
   \frac{W(V(\sigma_1))}{W(V(\sigma_2))}=(-1)^{\sum_{\ell \mid pNM_1M_2} \  (\delta(A(\sigma_1),\ \ell)- \delta(A(\sigma_2),\ \ell))}. 
\end{equation*}

In this setting (Panchishkin at $p$), Greenberg has conjectured  (see \cite[Conjecture 3.1]{Gr5})  that $S_\BK(A(\sigma_1)/\Q_\infty)^\vee$ is $\Lambda$-torsion. Further, if $\mu(S_\BK(A(\sigma_1)/\Q_\infty)^\vee)=0$, then all the hypotheses of Theorem \ref{relative_p_parity_over_Q} are satisfied and we see that the $p$-parity conjecture  for $V(\sigma_1)$ is equivalent to the $p$-parity conjecture  for $V(\sigma_2)$.


 In a special case, when $V=V_f(1-\frac{k}{2})$, where $f\in S_k(\Gamma_0(N))$ is a $p$-ordinary elliptic cuspidal newform with even weight $k$ and trivial nebentypus. If $\bar{\rho}_f$ is irreducible as a $G_{K_1}$-module, then $\mu(S_\BK(A(\sigma_1)/\Q_\infty)^\vee)$ is conjectured to be $0$ (as explained in Example \ref{example_modular-forms}). Thus the hypothesis (1) of Theorem \ref{relative_p_parity_over_Q} is expected to be true for $V(\sigma_1)$ and all other hypotheses of Theorem \ref{relative_p_parity_over_Q} are unconditionally verified and the relative $p$-parity conjecture is applicable for $(V(\sigma_1),V(\sigma_2))$.

\subsubsection{Symmetric powers of elliptic curves:} \label{symmetric_power_elliptic} Fix a positive integer $k$. Then for any prime $p>2k+1$, consider two $p$-ordinary elliptic curves $E_1$ and $E_2$ defined over $\Q$ with $E_1[p]\cong E_2[p]$ as $G_\Q$-modules. For $m=1,2$, let $V_pE_m$ denote the Galois representation associated with $E_m$ and $V_m:=(\Sym^{2k+1}V_pE_m)(-k)$ denote the symmetric $(2k+1)$-th power representation of $V_pE_m$. Let $T_pE_1$ and $T_pE_2$ be the Tate module associated to $E_1$ and $E_2$, respectively and set $T_m= (\Sym^{2k+1}T_pE_m)(-k)$. Then the map 
\begin{equation}\label{pairing_on_symmetric_cube}
(v_1\otimes \cdots \otimes  v_{2k+1}, w_1\otimes \cdots \otimes  w_{2k+1}) \mapsto \frac{1}{(2k+1)!}\sum_{\tau \in S_{2k+1}} \prod_{i=1}^{2k+1} \langle v_i, w_{\tau(i)}\rangle_m
\end{equation}
 induces a perfect skew-symmetric pairing on $T_m$, where $\langle\ ,\ \rangle_m$ is the Weil pairing on $T_pE_m$. 
In particular,  $T_m$ is a symplectic self-dual $G_\Q$-stable lattice of $V_m$ and put $A_m:= V_m/T_m$. Then 
\begin{enumerate}
    \item the representation $V_m$ of $G_\Q$ is semistable (hence de Rham) at $p$, symplectic self-dual, pure of weight $-1$ and satisfies the Panchishkin condition at $p$.
   \item Again $E_m^-[p]$ is the maximal unramified quotient of $E_m[p]$ and hence the isomorphism $E_1[p]\cong E_2[p]$ when restricted to $G_{\Q_p}$ induces the isomorphism $E_1^-[p]\cong E_2^-[p]$. 
  The $p$-ordinary filtration (see Definition \ref{p-ordinary}) on $V_pE_m$ induces a $p$-ordinary filtration on $V_m$. The isomorphism $E_1[p]\cong E_2[p]$ of $G_\Q$-modules induces an isomorphism $A_1[p]\cong A_2[p]$  and 
   the isomorphism $E_1^-[p]\cong E_2^-[p]$ induces an isomorphism $A_1^-[p]\cong A_2^-[p]$ as $G_{\Q_p}$-modules.
\end{enumerate}

Let $N_1$ and $N_2$ be the conductor of $E_1$ and $E_2$, respectively.
By Theorem \ref{global_root_lambda}, we deduce 
\begin{equation*}
   \frac{W(V_1)}{W(V_2)}=(-1)^{\sum_{\ell \mid pN_1N_2} \  (\delta(A_1,\ \ell)- \delta(A_2,\ \ell))}. 
\end{equation*}

Moreover, in this setting  $S_\BK(A_1/\Q_\infty)^\vee$ is expected to be $\Lambda$-torsion by  \cite[Conjecture 3.1]{Gr5} and further, if $\bar{\rho}_{E_1}$ satisfies certain big image hypothesis (see \emph{Hypothesis $N$} in \cite[Definition 1.1]{dr}), then it is conjectured by Drinen \cite{dr} that $\mu(S_\BK(A_1/\Q_\infty)^\vee)=0$.  

In particular, if $E_1$ is  a non CM elliptic curve defined over $\Q$, then by  Serre's open image theorem, the required big image hypothesis is satisfied in \cite{dr} for all but finitely many primes $p$ and for such a $p$,   it is conjectured (see \cite[Conjecture 1.1]{dr}) that if  $S_\BK(A_1/\Q_\infty)^\vee$ is $\Lambda$-torsion, then $\mu(S_\BK(A_1/\Q_\infty)^\vee)=0$. Thus hypothesis $(1)$ in Theorem \ref{relative_p_parity_over_Q} is expected  to be satisfied for $V_1$ and $V_2$ and all other hypotheses in Theorem \ref{relative_p_parity_over_Q} are verified to be true.
Thus  the $p$-parity conjecture   for  $(\Sym^{2k+1}V_pE_1)(-k)$ is equivalent to  the $p$-parity conjecture   for$(\Sym^{2k+1}V_pE_2)(-k)$.

\subsubsection{Tensor product of symmetric powers:}\label{tensor_product_symmetric_powers}
Fix any two natural numbers $m$ and $n$. Then for any prime $p>\max \{ 2m+1,\ 2n\}$, 
consider four $p$-ordinary elliptic curves $E_1$, $E_2$, $E_3$ and $E_4$ defined over $\Q$ with $E_1[p]\cong E_3[p]$ and $E_2[p]\cong E_4[p]$ as $G_\Q$-modules. Let $V_pE_k$ denote the Galois representation associated with $E_k$ for $k\in \{1,2,3,4\}$.  Consider the Galois representations of $G_\Q$
\begin{equation*}
    V_1:= (\Sym^{2m+1}V_pE_1)\otimes (\Sym^{2n}V_pE_2)(-m -n) \quad \text{and} \quad  V_2:= (\Sym^{2m+1}V_pE_3)\otimes (\Sym^{2n}V_pE_4)(-m -n).
\end{equation*}
Again for such $p$, equation \eqref{pairing_on_symmetric_cube} induces a perfect skew-symmetric pairing on $(\Sym^{2m+1}T_pE_k)(-m)$ and perfect symmetric pairing on $(\Sym^{2n}V_pE_k)(-n)$. Observe that
\begin{enumerate}
    \item the representations $V_1$ and $V_2$ of $G_\Q$ are semistable at $p$, symplectic self-dual, pure of weight $-1$ and satisfies the Panchishkin condition at $p$.
   \item The isomorphisms $E_1[p]\cong E_3[p]$ and $E_2[p]\cong E_4[p]$ of $G_\Q$-modules give the   isomorphisms $A_1[p]\cong A_2[p]$ of $G_\Q$-modules and when restricted to $G_{\Q_p}$, induces an isomorphism $A_1^-[p]\cong A_2^-[p]$.
\end{enumerate}

For $i\in \{1,2,3,4\}$, let $N_i$  be the conductor of $E_i$.
By Theorem \ref{global_root_lambda}, we deduce 
\begin{equation*}
   \frac{W(V_1)}{W(V_2)}=(-1)^{\sum_{\ell \mid pN_1N_2N_3N_4} \  (\delta(A_1,\ \ell)- \delta(A_2,\ \ell))}. 
\end{equation*}
Now under the assumption $S_\BK(A_1/\Q_\infty)^\vee$ is $\Lambda$-torsion (which is expected to be true by \cite[Conjecture 3.1]{Gr5}) and $\mu(S_\BK(A_1/\Q_\infty)^\vee)$ is equal to $0$, $V_1$ and $V_2$ satisfy all the conditions of Theorem \ref{relative_p_parity_over_Q}. Thus  the $p$-parity conjecture  for $V_1$ would be equivalent to the $p$-parity conjecture for $V_2$.

\begin{remark}
    Relative $p$-parity  results analogous to Examples \ref{symmetric_power_elliptic}, \ref{tensor_product_symmetric_powers} would hold if we consider symmetric powers of modular forms instead of symmetric powers of elliptic curves.  
\end{remark}

\section{New cases of the $p$-parity conjecture}\label{some_new_cases_p_parity} 

In this section, we will apply the relative  $p$-parity results of  \S \ref{relative_p_parity_conjecture} to deduce several cases of $p$-parity conjecture for 
certain families of $p$-adic Galois representations.

Combining the relative $p$-parity result in Theorem \ref{relative_p_parity_over_Q} with a result of Castella \cite[Theorem C]{Castella_Heegner_points}, we can now prove the  $p$-parity result  for certain Artin twists of cuspforms in a CM Hida family, stated in  Theorem \ref{castella_main_thm}.

\textbf{Proof of Theorem \ref{castella_main_thm}:}
  Under the hypotheses (Ca1)-(Ca4),  it is shown in   \cite[Theorem C]{Castella_Heegner_points} that $L(E\otimes \sigma/\Q,1)\neq 0 \Longrightarrow S_\BK(E\otimes\sigma/\Q)^\vee$ is finite. Put $V:= V_pE \otimes W_\sigma$. As $p\nmid NM$ and $V$ is pure of weight $-1$, we deduce $H^0(\Q_{p}, V^-)=0$ and also  from Lemma \ref{vanishing_V_m}, we have $H^0(\Q, V)$ is equal to zero. Hence, by a  proof similar to  Mazur's control theorem, we deduce that $(S_\BK(E\otimes\sigma/\Q_\infty)^\vee)^{\Gal(\Q_\infty/\Q)}$ is finite. Thus $S_\BK(E\otimes\sigma/\Q_\infty)^\vee$ is $\Lambda$-torsion. Further, by our assumption (5), $\mu(S_\BK(E\otimes\sigma/\Q_\infty)^\vee)=0$. Then it is easy to see that all the hypotheses of  Theorem \ref{relative_p_parity_over_Q} are satisfied  and the result is immediate from the same theorem. 

  Observe  that   $E$ being CM, $\sigma$ being two dimensional and assumptions (Ca1) to (Ca4) come from \cite[Theorem C]{Castella_Heegner_points}. \qed

\begin{remark}
Note that we have shown  $S_\BK(E\otimes\sigma/\Q_\infty)^\vee$ is $\Lambda$-torsion  in  Theorem \ref{castella_main_thm}. Further, if we assume $E[p]$ is an irreducible  $G_K$-module, then conjecturally $\mu(S_\BK(E\otimes\sigma/\Q_\infty)^\vee)=0$.
\end{remark}

Recall that Shou-Wu Zhang \cite{Shouwu_Zhang_Heights_Of_Heggner_points} has extended Gross-Zagier and Kolyvagin's result to totally real number fields. We  use his result    to establish   $p$-parity result for another infinite family of $p$-adic Galois representations stated in Theorem \ref{zhang_p_parity}.




\textbf{Proof of Theorem \ref{zhang_p_parity}:}
Let $B_f$ be the abelian variety defined over $K$ corresponding to the weight $2$ newform  $f\in S_2(\Gamma_0(M))$. 
For each prime $\lambda\mid p$ of  $\Q_f$, let $V_{f,\lambda}$ be the $\lambda$-adic Galois representation attached to $f$. Then the $p$-adic Galois representation $V_p(B_f)$ coming from the Tate module decomposes as $V_p(B_f) \cong \oplus_{\lambda \mid p} V_{f,\lambda}$.
Let $O_f$ be the ring of integers of the Hecke field $\Q_f$. 
By \cite[Theorem A]{Shouwu_Zhang_Heights_Of_Heggner_points}, 
\begin{equation}\label{zhang-result}
   \mathrm{ord}_{s=1} L(f/K,s)\leq 1 \Rightarrow  \rk_{\Q_p}\ S_{\BK}(V_p(B_f)/K) = \rk_{\Z}\ B_f(K)=[O_f:\Z]\ \mathrm{ord}_{s=1} L(f/K,s).
\end{equation}

By Artin formalism of $L$-function, we can write 
\begin{equation}
\label{L-function}
    L(f/K,s)=\prod_\sigma L(f\otimes \sigma/\Q, s)^{\dim \sigma},
\end{equation}
where $\sigma$ varies over irreducible Artin representation of $G=\Gal(K/\Q)$.

By base change results \cite[Theorem~1.1, with $n=2$]{CNT}, for every totally real number field $F$, the restriction of $V_f$ to $G_F$ corresponds to a Hilbert modular form of parallel weight $2$ over $F$.
Moreover, for every finite-order Hecke character $\chi$ of $F$,
the function $L(f_F\otimes\chi,s)$ is entire (see   \cite[Theorem 2]{MM}).
By Brauer induction,  we have the meromorphic continuation of
$L(f\otimes\sigma,s)$ for every $\sigma\in\operatorname{Irr}(G)$. Put $r=\operatorname{ord}_{s=1}L(f/K,s)\in\{0,1\}$, $r_\sigma=\operatorname{ord}_{s=1}L(f\otimes\sigma,s)\in \Z$ and $$\Theta=\sum_{\sigma\in\operatorname{Irr}(G)}r_\sigma\chi_\sigma,$$
where $\chi_\sigma$ is the character of $\sigma$.  
From~\eqref{L-function}, we get  $\Theta(1)=r$.

Let $C$ be a cyclic subgroup of $G$, set $F:=K^C$, and let
$\chi$ be a one-dimensional character of $C$, viewed also as
a finite-order Hecke character of $F$.
Frobenius reciprocity and Artin formalism give
\[
\begin{aligned}
\bigl\langle\operatorname{Res}^{G}_{C}\Theta,\chi\bigr\rangle_C
=\operatorname{ord}_{s=1}
L\bigl(f\otimes\operatorname{Ind}_{C}^{G}\chi,s\bigr)
=\operatorname{ord}_{s=1}L(f_F\otimes\chi,s)
\geq 0.
\end{aligned}
\]
Thus $\operatorname{Res}^{G}_{C}\Theta$ is a genuine character
of degree $r$.
Applying this to $C=\langle g\rangle$ shows that
$|\Theta(g)|\leq r$ for every $g\in G$.
Character orthogonality therefore yields
$$\sum_{\sigma\in\operatorname{Irr}(G)}r_\sigma^2
=\langle\Theta,\Theta\rangle_G
=\frac{1}{|G|}\sum_{g\in G}|\Theta(g)|^2
\leq r^2.$$
If $r=0$, it follows that $r_\sigma=0$ for every $\sigma$.
If $r=1$, there is a unique $\rho \in \mathrm{Irr} (G)$ such that $r_\rho= \pm 1$ and $r_\sigma=0$  for $\sigma\neq \rho$. As $\Theta(1)=r_\rho \chi_\rho(1)$, we get $r_\rho=1$ and $\dim \rho =1$. 
In particular, for  every $\sigma$, $L(f\otimes\rho,s)$ is holomorphic at $s=1$.

Also the isomorphism $V_p(B_f) \cong \oplus_{\lambda \mid p} V_{f,\lambda}$ gives an isomorphism of (vector space) Bloch-Kato Selmer groups 
\begin{equation*}
    S_{\BK}(V_p(B_f)/K) \cong \oplus_{\lambda \mid p} S_{\BK}(V_{f,\lambda}/K).
\end{equation*}

Case I: If $\mathrm{ord}_{s=1} L(f/K,s)=0$, then $S_{\BK}(V_{f,\lambda}/K)=0$  for any $\lambda\mid p$ of $\Q_f$ and hence  $S_{\BK}(V_{f,\lambda} \otimes W_\sigma/K)=0$, for any  Artin representation $\sigma$ of $G$.
Now using Hochschild-Serre, we have an isomorphism $ S_{\BK}(V_{f_\lambda}\otimes W_\sigma/K)^G \cong S_{\BK}(V_{f_\lambda} \otimes W_\sigma/\Q)$ and we obtain 
\begin{equation*}
    \rk_{O_f}\ S_{\BK} (A_{f_\lambda}\otimes\sigma/\Q)^\vee =0.
\end{equation*}

Thus 
\begin{equation*}
    \mathrm{ord}_{s=1}\ L(f\otimes \sigma/\Q, s)=0=\rk_{O_f}\ S_{\BK} (A_{f_\lambda}\otimes\sigma/\Q)^\vee.
\end{equation*}


\n Case II: If $\mathrm{ord}_{s=1} L(f/K,s)=1$, then 
\begin{equation*}
 [O_f:\Z]  = \rk_{\Q_p}\ S_{\BK}(V_p(B_f)/K)= \sum_{\lambda\mid p} [\Q_{f,\lambda}:\Q_p]\  \rk_{\Q_{f,\lambda}}\  S_{\BK}(V_{f,\lambda}/K). 
\end{equation*}
  Now using the fact that $\sum_{\lambda\mid p} [\Q_{f,\lambda}:\Q_p] = [O_f:\Z]$ and the known result of the $p$-parity conjecture for $V_{f,\lambda}$ over $K$ (see \cite[Theorem C]{Ne6}), we obtain that for all $\lambda\mid p$,
\begin{equation*}
    \rk_{\Q_{f,\lambda}}\  S_{\BK}(V_{f,\lambda}/K)=1=\mathrm{ord}_{s=1} L(f/K,s).
\end{equation*}

 Since $L(f/K,s)$ has simple zero at $s=1$, from \eqref{L-function}, it is clear that only one $1$-dimensional Artin representation of $\Gal(K/\Q)$ should contribute to the vanishing of $L(f/K,s)$, say $\sigma$. Using the functional equation, the conjugate of $\sigma$ would also contribute to the vanishing of $L(f/K,s)$ at $s=1$. This forces $\sigma$ to be a self-dual $1$-dimensional Artin representation, and hence $\sigma$ is either trivial or a quadratic character of $\Gal(K/\Q)$.
 And in both cases, the $p$-parity conjecture for $f\otimes \sigma /\Q$ is known due to \cite[Theorem B]{Ne4}.
 
 On the other hand, given $\mathrm{ord}_{s=1}L(f/K,s)=1$,
 if the dimension of an orthogonal Artin representation $\sigma$ is $\geq 2$, then $L(1,f\otimes \sigma/\Q)\neq 0$. Since $\rk_{O_f}\ S_{\BK}(V_{f,\lambda}/K)^\vee =1$ and degree of $\sigma$ is $\geq 2$ implies that the multiplicity of $\sigma$ in $ S_{\BK}(V_{f,\lambda}/K)$ is $0$. Again as an application of Hochschild-Serre, we obtain $ S_{\BK}(V_{f,\lambda} (\sigma)/\Q)=0$.

 Combining Case-I and II,  for any self-dual, orthogonal Artin representation $\sigma$ of $\Gal(K/\Q)$, the $p$-parity conjecture holds  for $f\otimes \sigma$. Now for every member $f_n$ of the  Hida family associated $f$ with $ k_n\equiv 2 \pmod{2(p-1)}$, the $p$-parity conjecture  for $f_n\otimes \sigma$ follows directly from the relative $p$-parity result in Theorem \ref{relative_p_parity_over_Q}.    \qed

Next we discuss a significant  generalization of $p$-parity results of Dokchitser-Dokchitser \cite{Dokchister-Dokchister_parity-conjecture} in our setting. One of the key ingredients in the $p$-parity results obtained in \cite{Dokchister-Dokchister_parity-conjecture} is their {\it Hypothesis $4.1$}. At first, using our relative $p$-parity result Theorem \ref{relative_p_parity_over_Q}, we  prove that the {\it Hypothesis $4.1$} as stated in \cite{Dokchister-Dokchister_parity-conjecture}, is satisfied in our setting.

\begin{lemma}\label{hypothesis_4.1}
 Let $V$ be a symplectic self-dual $p$-adic representation of $G_\Q$ satisfying \hyperref[Hyp:Mov]{(Mov)}. Also assume that $V$  satisfies the Panchishkin condition at $p$ and pure (of weight $-1$). Let $K/\Q$ be a
Galois extension of number fields. Assume that $S_\BK(A/K_\infty)^\vee$ is $\Lambda$-torsion and $\mu(S_\BK(A/K_\infty)^\vee)=0$. Then, for any subgroups \(N \triangleleft U\) of \(\operatorname{Gal}(K/\mathbb{Q})\) with
\(U/N \cong D_{2p^n}\), and for any representation
\(\tau=\sigma \oplus \mathbbm{1} \oplus \det\sigma\), where \(\sigma\) is a two-dimensional
representation of \(U/N\), the \(p\)-parity conjecture holds for the twists of \(V\) by \(\tau\) over $K^U$. 
\end{lemma}

\begin{proof}
Let $\sigma$ be a two dimensional representation of $U/N$. If $\sigma$ is reducible, then $\sigma \cong \mathbbm{1} \bigoplus \det \sigma$ or $\sigma \cong \varepsilon \bigoplus \varepsilon$, for a quadratic character $\varepsilon$ and the $p$-parity is immediate. On the other hand, if $\sigma$ is irreducible, then $\sigma\cong \mathrm{Ind} \chi$, where $\chi$ is a character of $p$-power order. Note that $\mathrm{Ind} \chi$ is a representation defined over the ring of integers of $\Q_p(\zeta_{p^n})$. Moreover, the corresponding $\Z_p[\zeta_{p^n}]$ lattice is self-dual  under the symmetric pairing given by the matrix $\left(\begin{array}{cc}
   0  & 1 \\
   1  & 0
\end{array}\right)$, and the residual representation $\bar\sigma $ is isomorphic to \(\bar{\mathbbm{1}} \bigoplus \bar{\det \sigma}\).

Then under the assumption that $S_\BK(A/K_\infty)^\vee$ is $\Lambda$-torsion and $\mu(S_\BK(A/K_\infty)^\vee)=0$, we deduce via Theorem \ref{relative_p_parity_over_Q} (and Remark \ref{relative_p_parity_over_L}), \begin{equation}\label{dd_ypothesis_4.1}
    \frac{W(V\otimes \sigma)}{W(V\otimes(\mathbbm{1} \oplus \det \sigma)) } = (-1)^{\rk_\OO\ S_\BK(A\otimes (\mathbbm{1} \oplus \det \sigma)/K^U)^\vee - \rk_\OO\ S_\BK(A\otimes \sigma/K^U)^\vee}.
\end{equation}
Now observe that $\frac{W(V\otimes \sigma)}{W(V\otimes(\mathbbm{1} \oplus \det \sigma)) }= W(V\otimes \sigma)W(V\otimes(\mathbbm{1} \oplus \det \sigma)) = W(V\otimes \tau)$ and $\rk_\OO\ S_\BK(A\otimes (\mathbbm{1} \oplus \det \sigma)/K^U)^\vee + \rk_\OO\ S_\BK(A\otimes \sigma/K^U)^\vee= \rk_\OO\ S_\BK(A\otimes \tau /K^U)^\vee$. Thus from equation \eqref{dd_ypothesis_4.1}, we obtain that 
\begin{equation*}
    W(V\otimes \tau)=(-1)^{ \rk_\OO\ S_\BK(A\otimes \tau /K^U)^\vee}.
\end{equation*}
\end{proof}

\begin{remark} It is shown in  \cite[Theorem 4.2]{Dokchister-Dokchister_parity-conjecture} that {\it Hypothesis $4.1$} is satisfied for abelian varieties $A$ satisfying suitable hypotheses.
  Under the assumption that $V=V_pA$ is $p$-ordinary, $S_\BK(A/K_\infty)^\vee$ is $\Lambda$-torsion and $\mu(S_\BK(A/K_\infty)^\vee)=0$,  Lemma \ref{hypothesis_4.1} gives a different proof  \cite[Theorem 4.2]{Dokchister-Dokchister_parity-conjecture} for $p>2$.\qed 
\end{remark}

Now using Lemma \ref{hypothesis_4.1},  we establish the following result which generalizes \cite[Theorem 4.5]{Dokchister-Dokchister_parity-conjecture} from $V=V_pE$ coming from an elliptic curve to a general $p$-adic  de Rham representation $V$ of $G_\Q$ that is  Panchishkin at $p$, under suitable hypotheses.
\begin{maintheorem}\label{dokchister_dokchister_theorem4.5_general}
     Let $V$ be a symplectic self-dual $p$-adic representation of $G_\Q$ satisfying \hyperref[Hyp:Mov]{(Mov)}. Also assume that $V$  satisfies the Panchishkin condition at $p$ and is pure (of weight $-1$). Let $K/\Q$ be a Galois extension of number fields such that the Sylow $p$-subgroup $P$ of $G=\Gal(K/\Q)$ is normal and $G/P$ is abelian. We further  assume that $S_\BK(A/K_\infty)^\vee$ is $\Lambda$-torsion and $\mu(S_\BK(A/K_\infty)^\vee)=0$. If the $p$-parity conjecture holds for $V$ over $\Q$ and its quadratic extensions in $K$, then it holds for all twists of $V$ by orthogonal representations of $G$.
\end{maintheorem}
\begin{proof}
From the proof of \cite[Theorem 4.3 $\&$ 4.5]{Dokchister-Dokchister_parity-conjecture} and by Lemma \ref{hypothesis_4.1}, it is clear that to prove our theorem it is enough to show the following. 
\begin{enumerate}
    \item [(a)] For an Artin representation $\tau :\Gal(K/\Q) \lra \mathrm{GL}_n(\C)$, 
\begin{equation*}
        \langle \tau , S_\BK(A/K)^\vee\rangle =  \langle \tau^* , S_\BK(A/K)^\vee\rangle \; \; \text{and} \;\; W(V\otimes\tau)= \overline{W(V\otimes\tau^*)},
\end{equation*}
where $\tau^*$ denotes the contragredient representation of $\tau$.
\item[(b)] If $\Q\subset L \subset K$ and for an Artin representation $\rho$ of $\Gal(K/L)$, let $\tau = \mathrm{Ind}_{\Gal(K/L)}^{\Gal(K/\Q)}\rho$. 
\begin{equation*}
 \text{Then }\;\;     \langle \tau , S_\BK(A/K)^\vee\rangle=  \langle \rho , S_\BK(A/K)^\vee\rangle \;\; \text{and} \;\; W(V/\Q \otimes \tau) = W(V/L \otimes \rho).
\end{equation*}
\end{enumerate}

For the proof of (a), $\langle \tau , S_\BK(A/K)^\vee\rangle =  \langle \tau^* , S_\BK(A/K)^\vee\rangle $ holds true because of \cite[(3.3)]{Ne5} and $W(V\otimes\tau)= \overline{W(V\otimes\tau^*)}$ is immediate from \cite[(3.4.7)]{Tate}.

For  the proof of (b), $\langle \tau , S_\BK(A/K)^\vee\rangle=  \langle \rho , S_\BK(A/K)^\vee\rangle$ is Frobenius reciprocity and $W(V/\Q \otimes \tau) = W(V/L \otimes \rho)$ holds by the inductivity property of the global root number $W$ (see \cite[Proposition 1]{Ro3}).

The rest of the proof is similar to \cite[Theorem 4.5]{Dokchister-Dokchister_parity-conjecture}. Note that the assumption $S_\BK(A/K_\infty)^\vee$ is $\Lambda$-torsion and $\mu(S_\BK(A/K_\infty)^\vee)=0$ are used in Lemma \ref{hypothesis_4.1}.
\end{proof}

As an  application of Theorem \ref{dokchister_dokchister_theorem4.5_general}, we derive the $p$-parity conjecture for  orthogonal Artin twists of modular forms stated in Theorem \ref{P_parity_as_DD}.

\textbf{Proof of Theorem \ref{P_parity_as_DD}:}
     By a result of Nekov\'a\v r (\cite[Theorem B]{Ne4}),  the $p$-parity conjecture holds for $V_f(1-k)$ and quadratic twists of $V_f(1-k)$.  Consequently, $p$-parity conjecture holds for $V_f(1-k)$ over any quadratic extensions of $\Q$ inside $K$. Now the result follows from Theorem \ref{dokchister_dokchister_theorem4.5_general}. \qed

Finally, combining the relative $p$-parity result in Theorem \ref{relative_p_parity_over_Q} with the result of H. Wang \cite[Theorem 5.8]{wa}, we can complete the proof of the $p$-parity result for symmetric cubes of modular forms stated in Theorem \ref{wang_result_p_parity}. 

\textbf{Proof of Theorem \ref{wang_result_p_parity}:}
    Under the assumptions $(\mathrm{Wa1})- (\mathrm{Wa5})$, H. Wang \cite[Theorem 5.8]{wa} has shown that $S_\BK(A/\Q)^\vee$ is finite. The stated $p$-parity result for the symmetric cube family $\{V_n\}$  holds by Theorem \ref{relative_p_parity_over_Q} and Example \ref{symmetric_power_elliptic}. \qed

\begin{remark}
 If we further assume $p\nmid N$,   then $S_\BK(A/\Q_\infty)^\vee$ is $\Lambda$-torsion. Indeed, considering  the Frobenius eigenvalue at $p$, we observe that $H^0(\Q_{p}, V^-)=0=H^0(\Q, V)$. By applying a control theorem, we deduce $S_\BK(A/\Q_\infty)^\vee$ is $\Lambda$-torsion. \qed 
\end{remark}

We discuss a numerical example where all the hypotheses are simultaneously met and we unconditionally deduce the $p$-parity conjecture for an Artin twist of a Hida family of cuspforms.

\subsection{Numerical Examples}\label{numerical_examples}
\begin{itemize}
    \item[(Ex-1)] Let $p=7$. Consider the elliptic curve $E/\Q$ with Cremona label $11a2$  and its quadratic twist $E_{-19}/\Q$ with Cremona label $3971b3$. Then $E$ and  $E_{-19}$ are ordinary at $p=7$. We  verify from \cite{lmfdb} that $\mu_{an}$ invariant of both $E$ and $E_{-19}$ over $\Q$ is $0$. Therefore the $\mu_{an}$ invariant is $0$ for $E/L$ with $L=\Q(\sqrt{-19})$. By Kato's work on   Iwasawa main conjecture (\cite{Kato_p_adic_hodge_theory}), we know $S_\str (E/L_\infty)^\vee$ is $\Lambda_L$-torsion and we also deduce   the algebraic  $\mu$ invariant  of  $S_\str (E/L_\infty)^\vee$ is $0$. Consider a Galois extension $F/L$ so that $\Gal(F/\Q)\cong D_{2\times7^n}$. (As an example, the number field $F$ with LMFDB label $14.0.20825499863802119464459.1$ is Galois over $\Q(\sqrt{-19})$ with $\Gal(F/\Q)\cong D_{14}$). Using Kida's formula \cite[Corollary 3.4]{Hachimori_Matsuno}, we obtain $S_\str (E/F_\infty)^\vee$ is a $\Lambda$-torsion module and the $\mu$ invariant of $S_\str (E/F_\infty)^\vee$ is $0$.  
 
 Let $\sigma$ be an irreducible, orthogonal representation of $\Gal(F/\Q)$, then $7$-parity conjecture is true unconditionally for $E\otimes \sigma$ by Theorem \ref{P_parity_as_DD} (see also \cite[Theorem 1.3]{Dokchister-Dokchister_parity-conjecture}, \cite[Theorem 1.15]{Roc}).  Then from our Theorem \ref{relative_p_parity_over_Q} and Example \ref{example_modular-forms}, we  get an infinite family of modular forms $\{f_n\}_{n}$  with weight $k_n \equiv 2 \pmod{12}$ in the Hida family associated to  $E$ such that  $7$-parity conjecture holds unconditionally for $\{f_n \otimes \sigma\}_k$. 
\medskip

\item[(Ex-2)]  Let $p=3$ and consider two  elliptic curves $E_1= 272c1$ and $E_2= 275b1$ defined over $\Q$. Note that $E_1$ and $E_2$ are $3$-ordinary and from \cite[Table 3-272C1 \& 3-275B1]{dd2}, we have ($S_\BK(E_i/K_\infty)^\vee$ is $\Lambda_K$-torsion and)  $\mu(S_\BK(E_1/K_\infty)^\vee)=0= \mu(S_\BK(E_2/K_\infty)^\vee)$ for $K=\Q(\mu_3)$.  Let $m>1$ be a $3$-power free integer and set $F:= K(\sqrt[3]{m})$.  So by the Kida's formula \cite[Corollary 3.4]{Hachimori_Matsuno}, $S_\BK(E_i/F_\infty)^\vee$ is $\Lambda_F$-torsion and $\mu(S_\BK(E_i/F_\infty)^\vee)=0$. Hence for any irreducible orthogonal Artin representation $\sigma$ of $\Gal(F/\Q)$, by Theorem \ref{P_parity_as_DD}, the $3$-parity conjecture is true for $E_i\otimes \sigma$ over $\Q$. 
Therefore,  the $3$-parity conjecture is true for $(E_1 \bigoplus E_2)\otimes\sigma$ over $\Q$.  Now as explained in Example \ref{example_isomorphic_residual_but_not_symplectic_similitude}, there exists an  infinite family of principally polarised $3$-ordinary abelian varieties $B_n$ over $\Q$, up to twist, such that $B_n[3]\cong E_1[3]\bigoplus E_2[3]$ as $G_\Q$-module. 

Consequently, by Theorem \ref{relative_p_parity_over_Q}, the $3$-parity conjecture holds true (unconditionally) for $B_n\otimes \sigma$ over $\Q$  for the infinite family of abelian varieties.

\medskip 

\item[(Ex-3)] In \cite{dd2}, for prime $p=3,\ 5,\ 7$, the authors provide a list of examples of $p$-ordinary elliptic curves such that ($S_\BK(E/K_\infty)^\vee$ is $\Lambda_K$-torsion and)  $\mu(S_\BK(E/K_\infty)^\vee)=0$ for $K=\Q(\mu_p)$. Now for any $p$-power extension $L$ of $\Q(\mu_p)$, we have $S_\BK(E/L_\infty)^\vee$ is $\Lambda_L$-torsion and $\mu(S_\BK(E/L_\infty)^\vee)=0$ by Kida's formula \cite[Corollary 3.4]{Hachimori_Matsuno}. 

Also in \cite{cdlss}, the authors provided a list of examples of cuspidal newforms $f\in S_{2k}(\Gamma_0(N))$ of higher weight ($k>1$) with $\mu(S_\BK(A_f/K_\infty)^\vee)=0$ for $K=\Q(\mu_p)$.

\medskip

\item[(Ex-4)] Consider an odd prime $p$. Let $K/\Q$ be a Galois number field such that $G:= \Gal(K/\Q) \cong \Z/2r\Z$ with $\gcd(2r,p)=1$. Let $\sigma\in G$ be a generator of $G$ and set $c= \sigma^r$. 
Fix a lift $\tilde{\sigma}$ of $\sigma$ to $G_\Q$ and an element $\psi\in \Hom(G_K,\C^\times)$ of finite order prime to $p$ such that
\begin{equation*}
\psi^c=\psi^{-1},\quad
\mathrm{Stab}_G(\bar\psi)=1,\quad \psi(\tilde \sigma^{2r})=1,
\end{equation*}
where $\bar \psi$ denotes the  residual character of $\psi$.
Now consider a family of characters $\phi_\alpha \in \Hom(G_K,\C^\times)$ such that 
\begin{equation}\label{condition_on_phi}
\phi_\alpha^c=\phi_\alpha^{-1},\quad
\ord(\phi_\alpha)=p^{a_\alpha}.
\end{equation}
For each $\alpha$, consider the character $\eta_\alpha=\psi\phi_\alpha$ and set $\rho_\alpha=\mathrm{Ind}_{G_K}^{G_\mathbf Q}\ \eta_\alpha$. Then $\dim\rho_\alpha=2r$.  The following isomorphism 
\begin{align*}
T:\mathrm{Ind}_{G_K}^{G_\mathbf Q}\ \eta_\alpha^{c} & \lra \mathrm{Ind}_{G_K}^{G_\mathbf Q}\ \eta_\alpha\\
    g &\longmapsto (T\cdot g)(x):= g(cx)
\end{align*}
with the equality $\eta_\alpha^c =\eta_\alpha^{-1}$ implies that $\rho_\alpha^\vee\simeq\rho_\alpha$. As $\ord(\phi_\alpha)=p^{a_\alpha}$, we have
\begin{equation*}
\bar\rho_\alpha\simeq \mathrm{Ind}_{G_K}^{G_\mathbf Q}\ \bar\psi \quad \text{for each $\alpha$}.
\end{equation*}
Thus, we obtain that for each $\alpha$, the  representation $\rho_\alpha$ is orthogonal (see \cite[Proposition 4.1(3)]{PR}) and moreover, for any two $\alpha$ and $\beta$, $\bar\rho_\alpha \cong \bar\rho_\beta$. Also, our assumption $\mathrm{Stab}_G (\bar\psi)=1$ implies that $\bar\rho_\alpha$ is an absolutely irreducible representation of degree $2r$.

 Let $V$ be a symplectic self-dual $p$-adic representation of $G_\Q$ satisfying \hyperref[Hyp:Mov]{(Mov)}. Also assume that $V$  satisfies the Panchishkin condition at $p$ and is pure (of weight $-1$). Let $F$, $K_\alpha$, $K_\beta$ be the trivialising extension of $\psi$, $\phi_\alpha$, $\phi_\beta$, respectively and set $L:= FK_\alpha K_\beta$. 
Then under the assumption $S_\BK(A/L_\infty)^\vee$ is $\Lambda_L$-torsion and $\mu(S_\BK(A/L_\infty)^\vee)=0$, by Theorem \ref{relative_p_parity_over_Q}, the $p$-parity conjecture for $V\otimes \rho_\alpha$ over $\Q$ is equivalent to the $p$-parity conjecture for $V\otimes \rho_\beta$ over $\Q$.

Note that we are not assuming our characters $\eta_\alpha/\eta_\beta$ to be unramified at the prime above $p$. Thus our example is not covered by \cite[Theorem 4.16]{Ne}. 
The results in \cite{Dokchister-Dokchister_parity-conjecture} also do not consider  Artin representations of this type.

\end{itemize}

\begin{small}
	
\end{small}

\begin{thebibliography}{999999}
 \bibitem[BC]{BC}
O.~Brinon and B.~Conrad,
\emph{CMI Summer School Notes on $p$-adic Hodge Theory}, \url{https://math.stanford.edu/~conrad/papers/notes.pdf}. 
\bibitem[BCGP]{BCGP}
G.~Boxer, F.~Calegari, T.~Gee, and V.~Pilloni,
{Abelian surfaces over totally real fields are potentially modular},
{\it Publ. Math. Inst. Hautes \'Etudes Sci.} \textbf{134} (2021), 153--501.



\bibitem[Ben]{Be}D. Benois, An introduction to $p$-adic Hodge theory, in {\it Perfectoid spaces}, 69--219,  Infosys Sci. Found. Ser. Math. Sci.,  Springer, Singapore.

\bibitem[Ber]{Ber} L. Berger, Repr\'esentations $p$-adiques et \'equations diff\'erentielles, {\it Invent. Math.} {\bf 148} (2002), no.~2, 219--284.  

\bibitem[BHKT]{bhkt}G. B\"ockle,  M. Harris, C. Khare, J. Thorne, $\hat G$-local systems on smooth projective curves are potentially automorphic, {\it Acta Math.} {\bf 223} (2019), no.~1, 1--111.

\bibitem[BK]{bk}S.~J. Bloch and K. Kato, $L$-functions and Tamagawa numbers of motives, in {\it The Grothendieck Festschrift, Vol.\ I}, 333--400, Progr. Math., 86, Birkh\"auser Boston, Boston, MA.

\bibitem[BKNO]{bkno}  {A. Burungale, S. Kobayashi,  K. Nakamura and K. Ota, A local sign decomposition for symplectic self-dual Galois representations of rank two, \url{https://arxiv.org/abs/2508.17776}.} 

\bibitem[BS]{BS}B.~J. Birch and N.~M. Stephens, The parity of the rank of the Mordell-Weil group, {\it Topology} {\bf 5} (1966), 295--299.


 \bibitem[Ca]{Castella_Heegner_points} {F. Castella, On the $p$-adic variation of Heegner points, {\it J. Inst. Math. Jussieu} {\bf 19}(6) (2020), 2127--2164.}

\bibitem[CDLSS]{cdlss}J.~H. Coates, T.  Dokchitser, Z. Liang, W.  Stein, R. Sujatha,  Non-commutative Iwasawa theory for modular forms, {\it Proc. Lond. Math. Soc. (3)} {\bf 107} (2013), no.~3, 481--516.

\bibitem[CE]{CE} H.~Cartan and S.~Eilenberg, \emph{Homological Algebra}, Princeton University Press, Princeton, 1956.

\bibitem[Ce]{Ce}K. \v Cesnavi\v cius, The $p$-parity conjecture for elliptic curves with a $p$-isogeny, {\it J. Reine Angew. Math.} {\bf 719} (2016), 45--73.

\bibitem[CFKS]{CFKS}J.~H. Coates,  T. Fukaya, K. Kato, R. Sujatha, Root numbers, Selmer groups, and non-commutative Iwasawa theory, {\it J. Algebraic Geom.} {\bf 19} (2010), no.~1, 19--97.

\bibitem[CG]{CG} {J.~H. Coates and R. Greenberg, Kummer theory for abelian varieties over local fields, {\it Invent. Math.} {\bf 124} (1996), 129--174.}



\bibitem[CNT]{CNT}
L.~Clozel, J.~Newton and J.~A.~Thorne,
\emph{Non-abelian base change for symmetric power liftings
of holomorphic modular forms},
arXiv:2312.01774v2, 2024.




\bibitem[DD1]{Dokchister-Dokchister_parity-conjecture} {T. Dokchitser and V. Dokchitser, Regulator constants and the parity conjecture, {\it Invent. Math.} {\bf 178}(1) (2009), 23--71.}

 \bibitem[DD2]{dd2}T. Dokchitser and V. Dokchitser, Computations in non-commutative Iwasawa theory, {\it Proc. Lond. Math. Soc. (3)} {\bf 94} (2007), no.~1, 211--272.

 \bibitem[DD3]{DD3}T. Dokchitser and V. Dokchitser, On the Birch-Swinnerton-Dyer quotients modulo squares, {\it Ann. of Math. (2)} {\bf 172} (2010), no.~1, 567--596.

\bibitem[De1]{Deligne_73} {P. Deligne, Les constantes des \'equations fonctionnelles des fonctions $L$, in {\it Modular functions of one variable, II}, pp. 501--597, Lecture Notes in Math., Vol. 349, Springer.}

\bibitem[De2]{Deligne_79} {P. Deligne, Valeurs de fonctions $L$\ et p\'eriodes d'int\'egrales, in {\it Automorphic forms, representations and $L$-functions, Part 2}, pp. 313--346, Proc. Sympos. Pure Math., XXXIII, Amer. Math. Soc.}

\bibitem[Dr]{dr} M.~J. Drinen, Iwasawa $\mu$-invariants of elliptic curves and their symmetric powers, {\it J. Number Theory} {\bf 102} (2003), no.~2, 191--213. 

\bibitem[EPW]{epw}{M. Emerton, R.~J. Pollack and T.~A. Weston, Variation of Iwasawa invariants in Hida families, {\it Invent. Math.} {\bf 163}(3) (2006), 523--580.}


\bibitem[Fl]{Flach90}{M. Flach, {A generalisation of the Cassels-Tate pairing}, {\it J. Reine Angew. Math.} {\bf 412} (1990), 113--127}


\bibitem[FO]{FO} J.-M. Fontaine, Y. Ouyang, Theory of p-adic galois representations (2008). \url{https://www.math.u-psud.fr/~fontaine/galoisrep.pdf}.

\bibitem[Fo1]{Fontaine_94} {J.-M. Fontaine, Le corps des p\'eriodes $p$-adiques, {\it Ast\'erisque} No. 223 (1994), 59--111.}
\bibitem[Fo2]{Fo2}{J.-M. Fontaine, Repr\'esentations $p$-adiques semi-stables, {\it Ast\'erisque} No. 223 (1994), 113--184.}

\bibitem[FP]{fp}J.-M. Fontaine and B. Perrin-Riou, Autour des conjectures de Bloch et Kato: cohomologie galoisienne et valeurs de fonctions $L$, in {\it Motives}, 599--706, Proc. Sympos. Pure Math., 55, Part 1.



\bibitem[Gr1]{Greenberg1997}  {R. Greenberg, Iwasawa theory for elliptic curves, in {\it Arithmetic theory of elliptic curves (Cetraro, 1997)}, 51--144, Lecture Notes in Math., 1716, Springer, Berlin.}


\bibitem[Gr2]{Greenberg_projectivemodule_modularrepresentation}{R. Greenberg, Iwasawa theory, projective modules, and modular representations, {\it Mem. Amer. Math. Soc.} {\bf 211}(992) (2011), vi+185 pp.}

\bibitem[Gr3]{Gr3}R. Greenberg, Trivial zeros of $p$-adic $L$-functions, in {\it $p$-adic monodromy and the Birch and Swinnerton-Dyer conjecture (Boston, MA, 1991)}, 149--174, {\it Contemp. Math.}, 165, Amer. Math. Soc., Providence.



\bibitem[Gr4]{Greenberg_Iwasawa_p_adic_representation} {R. Greenberg, Iwasawa theory for $p$-adic representations, in {\it Algebraic number theory}, 97--137, Adv. Stud. Pure Math., 17.}

\bibitem[Gr5]{Gr5} R. Greenberg, Iwasawa theory and $p$-adic deformations of motives, in {\it Motives (Seattle, WA, 1991)}, 193--223, Proc. Sympos. Pure Math., 55, Part 2, Amer. Math. Soc., Providence, RI.

\bibitem[Gr6]{Gr6}R. Greenberg, On the Birch and Swinnerton-Dyer conjecture, {\it Invent. Math.} {\bf 72} (1983), no.~2, 241--265.

\bibitem[GS]{Greenberg_Stevens} {R. Greenberg and G. Stevens, $p$-adic $L$-functions and $p$-adic periods of modular forms, {\it Invent. Math.} {\bf 111}(2) (1993), 407--447.}


\bibitem[Gu1]{Guo_Li} {L. Guo, On a generalization of Tate dualities with application to Iwasawa theory, {\it Compositio Math.} {\bf 85}(2) (1993), 125--161.}

\bibitem[Gu2]{Gu2}L. Guo, General Selmer groups and critical values of Hecke $L$-functions, {\it Math. Ann.} {\bf 297} (1993), no.~2, 221--233.

 \bibitem[GV]{Greenberg_Vatsal} {R. Greenberg and V. Vatsal, On the Iwasawa invariants of elliptic curves, {\it Invent. Math.} {\bf 142}(1) (2000),  17--63.}


\bibitem[Hi1]{Hida_p_adic_measure_1} {H. Hida, A $p$-adic measure attached to the zeta functions associated with two elliptic modular forms. I, {\it Invent. Math.} {\bf 79}(1) (1985), 159--195.}
                
\bibitem[Hi2]{Hida_Galois_representation_into_GL2_attached_ordinarycuspforms} {H. Hida, Galois representations into ${\rm GL}_2({\bf Z}_p[[X]])$ attached to ordinary cusp forms, {\it Invent. Math.} {\bf 85}(3) (1986), 545--613.}

 \bibitem[HM]{Hachimori_Matsuno}{Y. Hachimori and K. Matsuno, An analogue of Kida's formula for the Selmer groups of elliptic curves, {\it J. Algebraic Geom.} {\bf 8}(3) (1999), 581--601.}


\bibitem[JN]{JN}C. Johansson and J.~J.~M. Newton, Parallel weight 2 points on Hilbert modular eigenvarieties and the parity conjecture, {\it Forum Math. Sigma} {\bf 7} (2019), Paper No. e27, 36 pp.


\bibitem[Ka]{Kato_p_adic_hodge_theory}{K. Kato, $p$-adic Hodge theory and values of zeta functions of modular forms, {\it Ast\'erisque} No. 295 (2004), ix, 117--290.}

\bibitem[Ki]{Ki}B.~D. Kim, The parity conjecture for elliptic curves at supersingular reduction primes, {\it Compos. Math.} {\bf 143} (2007), no.~1, 47--72.

\bibitem[Kn]{K}
F. Knop,
\emph{Classification of multiplicity free symplectic representations},
{\it J. Algebra} \textbf{301} (2006), 531--553.

\bibitem[LMFDB]{lmfdb}\textit{ The L-functions and modular forms database, https://www.lmfdb.org, 2026}

\bibitem[LS]{Lim-Sujatha_JNT_2018} {M.~F. Lim and R. Sujatha, Fine Selmer groups of congruent Galois representations, {\it J. Number Theory} {\bf 187} (2018), 66--91.}

\bibitem[LY]{ly} Z. Lin, M. Yin, On Greenberg's conjecture for rational elliptic curves at Eisenstein primes, \url{https://arxiv.org/abs/2609.40282}.


\bibitem[MM]{MM}M.~R. Murty and V.~P. Murty, Base change and the Birch-Swinnerton-Dyer conjecture, in {\it A tribute to Emil Grosswald: number theory and related analysis}, 481--494, Contemp. Math., 143, Amer. Math. Soc., Providence.

\bibitem[Mo]{Mo}P. Monsky, Generalizing the Birch-Stephens theorem. I. Modular curves, {\it Math. Z.} {\bf 221} (1996), no.~3, 415--420.

\bibitem[MR1]{MR}B. Mazur and K. Rubin, Finding large Selmer rank via an arithmetic theory of local constants, {\it Ann. of Math. (2)} {\bf 166} (2007), no.~2, 579--612.

\bibitem[MR2]{MR1}
B.~Mazur and K.~Rubin,
{Kolyvagin Systems},
{\it Mem. Amer. Math. Soc.} \textbf{168} (2004), no.~799.


\bibitem[Ne1]{Ne}J. Nekov\'a\v r, Compatibility of arithmetic and algebraic local constants (the case $\ell\ne p$), {\it Compos. Math.} {\bf 151} (2015), no.~9, 1626--1646.
\bibitem[Ne2]{Ne2}J. Nekov\'a\v r, Selmer complexes, {\it Ast\'erisque} No. 310 (2006), viii+559 pp.   

                
\bibitem[Ne3]{Ne3}J. Nekov\'a\v r, On the parity of ranks of Selmer groups. III, {\it Doc. Math.} {\bf 12} (2007), 243--274.

\bibitem[Ne4]{Ne4}J. Nekov\'a\v r, Some consequences of a formula of Mazur and Rubin for arithmetic local constants, {\it Algebra Number Theory} {\bf 7} (2013), no.~5, 1101--1120.

\bibitem[Ne5]{Ne5}J. Nekov\'a\v r, On the parity of ranks of Selmer groups. IV, {\it Compos. Math.} {\bf 145} (2009), no.~6, 1351--1359.

\bibitem[Ne6]{Ne6}J. Nekov\'a\v r, Compatibility of arithmetic and algebraic local constants, II: the tame abelian potentially Barsotti-Tate case, {\it Proc. Lond. Math. Soc. (3)} {\bf 116} (2018), no.~2, 378--427.



\bibitem[PR]{PR}D. Prasad and D. Ramakrishnan, Self-dual representations of division algebras and Weil groups: a contrast, Amer. J. Math. {\bf 134} (2012), no.~3, 749--767.








 \bibitem[Roh1]{Rohrlich_weil_deligne_group}{D.~E. Rohrlich, Elliptic curves and the Weil-Deligne group, in {\it Elliptic curves and related topics}, 125--157, CRM Proc. Lecture Notes, 4, Amer. Math. Soc., Providence, RI.}

 \bibitem[Roh2]{Ro3}D.~E. Rohrlich, Inductivity of the global root number, {\it Acta Arith.} {\bf 159} (2013), no.~3, 245--256.

 \bibitem[Roc]{Roc}T. de~La~Rochefoucauld, Invariance of the parity conjecture for $p$-Selmer groups of elliptic curves in a $D_{2p^n}$-extension, {\it Bull. Soc. Math. France} {\bf 139}(4) (2011), 571--592.


\bibitem[Si]{S}
J.~H. Silverman,
\emph{The Arithmetic of Elliptic Curves},
2nd ed.,
GTM vol.~106,
Springer, Dordrecht, 2009.

 \bibitem[Su]{Sujatha_elliptic_curve_iwasawa_mu0_conjecture} {R. Sujatha, Elliptic curves and Iwasawa's $\mu=0$ conjecture, in {\it Quadratic forms, linear algebraic groups, and cohomology}, 125--135, Dev. Math., 18, Springer, New York.}

 \bibitem[Tat]{Tate} {J.~T. Tate, Number theoretic background, in {\it Automorphic forms, representations and $L$-functions, Part 2}, pp. 3--26, Proc. Sympos. Pure Math., XXXIII, Amer. Math. Soc., Providence, RI.}

 \bibitem[Tay]{taylor}{R.~L. Taylor, {Galois representations}, {\it  Ann. Fac. Sci. Toulouse Math. (6)} {\bf 13}(1) (2004), 73--119.}

 \bibitem[Wa]{wa}H. Wang, Arithmetic level raising on triple product of Shimura curves and Gross-Schoen diagonal cycles, I: Ramified case, {\it Algebra Number Theory} {\bf 16} (2022), no.~10, 2289--2338.

\bibitem[Wi]{Wiles_ordinary_lambda_adic_representation} {A.~J. Wiles, On ordinary $\lambda$-adic representations associated to modular forms, {\it Invent. Math. }{\bf 94}(3) (1988),  529--573.}

\bibitem[Zh]{Shouwu_Zhang_Heights_Of_Heggner_points} {S.-W. Zhang, Heights of Heegner points on Shimura curves, {\it Ann. of Math. (2)} {\bf 153}(1) (2001), 27--147.}

 


               
 
 



 
 










 
 
   






 















\end{thebibliography}
\end{document}